\documentclass[arxiv,reqno,twoside,a4paper,12pt]{amsart}
\usepackage{amssymb}
\usepackage{amsthm}
\usepackage{graphicx}
\usepackage[margin=3cm]{geometry}

\usepackage{enumitem}

\usepackage[colorlinks=true,linkcolor=blue,citecolor=blue]{hyperref}

\usepackage[dvipsnames]{xcolor}
\usepackage{mathrsfs, tikz}
\usetikzlibrary{matrix,calc}
\usetikzlibrary{positioning}

\usepackage[breakable]{tcolorbox}
\tcbuselibrary{theorems}

\usepackage[color=yellow!30, linecolor=orange, colorinlistoftodos]{todonotes}
\usepackage[breakable]{tcolorbox}
\tcbuselibrary{theorems}
\newtcolorbox{bluebox}[1][]%
{left=0mm, right=0mm, bottom=0mm, top=0mm, sharp corners, boxrule=.8pt, before skip=\topsep, after skip=\topsep, colback=cyan!5, colframe=cyan, coltitle=black, fonttitle=\bfseries, title=#1, breakable}
\newtcolorbox{yellowbox}[1][]%
{left=0mm, right=0mm, bottom=0mm, top=0mm, sharp corners, boxrule=.8pt, before skip=\topsep, after skip=\topsep, colback=yellow!5, colframe=cyan, coltitle=black, fonttitle=\bfseries, title=#1, breakable}
\tcbset{tcbox raise base, sharp corners, left=0mm,right=0mm,top=0mm,bottom=0mm,nobeforeafter, boxrule =.5pt} 
\tcbset{tcbox raise base, sharp corners, left=0mm,right=0mm,top=0mm,bottom=0mm,nobeforeafter, boxrule =.5pt}

\newcommand{\define}[1]{\emph{#1}}
\DeclareMathOperator{\diag}{diag}
\DeclareMathOperator*{\Span}{span}
\DeclareMathOperator*{\range}{ran}
\DeclareMathOperator*{\im}{Im}
\newtheorem{theorem}{Theorem}[section]
\newtheorem{lemma}[theorem]{Lemma}
\newtheorem{proposition}[theorem]{Proposition}
\newtheorem*{proposition*}{Proposition}
\newtheorem{corollary}[theorem]{Corollary}
\newtheorem{example}[theorem]{Example}
\newtheorem{remark}[theorem]{Remark}
\newtheorem*{remark*}{Remark}

\numberwithin{equation}{section}
\theoremstyle{definition}
\newtheorem{definition}[theorem]{Definition}
\newtheorem*{notation*}{Notation}

\begin{document}

\title[Difference equations with operator coefficients]{Scattering theory 
for difference equations with operator coefficients}

\author[Sher]{David Sher}
\author[Silva]{Luis Silva}
\author[Vertman]{Boris Vertman}
\author[Winklmeier]{Monika Winklmeier}

\thanks{2000 Mathematics Subject Classification 39Axx,  47B39.}

\date{}
\dedicatory{}

\begin{abstract}
  We consider a second order difference equation with operator-valued
  coefficients. More precisely, we study either compact or trace class
  perturbations of the discrete Laplacian in the Hilbert space of bi-infinite
  square-summable sequence with entries in a fixed Hilbert space. We discuss
  its continuous and discrete spectrum, as well as properties of the
  associated scattering matrix.
\end{abstract}

\maketitle
\tableofcontents

\section {Introduction and statement of the main results} 

\subsection {Statement of the problem}\label{sec:statement-problem}
Let $H$ be a separable Hilbert space with the scalar product and norm denoted
by $\langle \cdot , \cdot \rangle$ and $\| \cdot \|$, respectively.  We
consider the Hilbert space $\mathcal{H}$ of bi-infinite square-summable
sequences with entries in $H$, namely
\begin{equation*}
  \mathcal H := \ell^2(\mathbb Z,H)\,, \quad
  \left\langle \frac{}{} \{v_{n}\}_{n\in\mathbb Z},\,
  \{w_{n}\}_{n\in\mathbb Z} \right\rangle_{\mathcal H} :=
  \sum_{n=-\infty}^\infty \langle v_n , w_n \rangle\,.
\end{equation*}
We denote by $\mathscr{B}(H)$ the space of bounded operators defined
on the whole space $H$.
Let us write 
$$
\ell^2(\mathbb Z,\mathscr{B}(H))
$$
for the space
of all sequences $T=\{T_n\}_{n\in\mathbb Z} \subset \mathscr{B}(H)$ such that
$$
 \|T\|^2:=
\sum_{n\in\mathbb Z} \|T_n\|^2 < \infty.
$$
Let the sequences $\{A_{n}\}_{n\in\mathbb Z}$ and $\{B_{n}\}_{n\in\mathbb Z}$
be such that $A_{n}, B_{n}\in\mathscr{B}(H)$ for all $n\in\mathbb{Z}$.
 Consider the difference expression
\begin{equation}
  \label{eq:main-difference}
    (\tau u)_n := A_{n-1}u_{n-1} + B_nu_n + A_nu_{n+1}\,.
\end{equation}
\emph{A priori}, this does not define a bounded operator on $\mathcal H$. Assuming
however that the sequences $\{A_{n}\}_{n\in\mathbb Z}$ and
$\{B_{n}\}_{n\in\mathbb Z}$ are uniformly bounded with respect to the
standard $\mathscr{B}(H)$ norm, the difference expression $\tau$ does define
a bounded operator $\mathcal J$ on $\mathcal H = \ell^2(\mathbb Z,H)$, which
is called a Jacobi operator (see Subsection~\ref{sec:pert-discr-lapl}). In this work,
we will actually need more than the uniform boundedness condition, the
so-called \emph{moment} conditions, which we summarize in the
following definition.

\begin{definition}
  \label{assumption:introduction}
  Let $\{A_{n}\}_{n\in\mathbb Z}, \{B_{n}\}_{n\in\mathbb Z} \subset \mathscr{B}(H)$ be
  sequences of self-adjoint operators.  We say that for any $k \in \mathbb{N}_0$
  the ``\define{$k$-th moment condition}'' holds,  if
     \begin{equation}\label{eq:starassumption}
	\sum_{n=-\infty}^{\infty}|n|^k(\|I-A_n\| + \|B_n\|) < \infty.
     \end{equation}
For instance, if the exponent $k=0, 1$ or $2$, we speak of the zero-th,  
first or second moment condition, respectively.  Note that the $k$-th moment
condition implies the $j$-th moment condition for any $j \in \mathbb{N}_0$ with $j \leq k$.
We say that the ``\define{exponential moment condition}'' is satisfied if for
     some $\varepsilon > 0$
     \begin{equation}\label{eq:starassumption3}
	\sum_{n=-\infty}^{\infty}e^{\varepsilon |n|}(\|I-A_n\| + \|B_n\|) < \infty.
     \end{equation}
\end{definition}
     
Note that exponential moment condition implies the $k$-th moment
condition for any $k \in \mathbb{N}$. 
Under any of these moment conditions, the Jacobi operator $\mathcal J$ may be considered as a
``small'' perturbation of the so-called \define{discrete free Laplacian}
$\mathcal J_{0}$ (see Section~\ref{sec:perturbed-laplacian}). In this paper
we study some of the main objects of stationary scattering theory for the
operator $\mathcal J$ with respect to $\mathcal J_{0}$.
\subsection{Historical background on discrete Schr\"odinger
  operators}\label{sec:background}
Second order difference equations, like their continuous counterparts
(Sturm-Liouville equations), arise in many applications in physics and
other sciences. In the particular case of Jacobi matrices
(\emph{i.\,e.} symmetric infinite tridiagonal matrices arising from the
difference expression $\tau$ with scalar coefficients), the
corresponding operators play a central role in the theory of simple
self-adjoint operators and the classical moment problem
\cite{akhiezer,MR1627806,teschl}.

Jacobi operators also have been used
directly to model systems with immediate neighbors interactions. In
quantum mechanics, a particular case of a Jacobi operator, namely the
discrete Schrödinger operator, corresponds to the so-called
\emph{tight binding} approximation \cite{MR883643}. The discrete
Schrödinger operator is given by a Jacobi matrix whose off diagonals
are all equal to $1$. For discrete Schrödinger
operators, the main diagonal corresponds to the potential. As in the
case of (continuous) one dimensional Schrödinger operators, two situations are
considered: the half-line problem (sequences on $\mathbb{N}$) and the
whole line problem (sequences on $\mathbb{Z}$).

Most of the spectral and scattering theory for difference equations
has been developed following the unfolding of the theory for
differential equations which for historical reasons began its
development earlier. The approach to the stationary scattering theory
for difference equations is therefore based on the classical approach
for Schrödinger operators (see \cite{MR985100} for a review and
references).

The work in \cite{MR0133609} seems to have been the first rigorous
attempt to solve the scattering problem for difference equations with
scalar coefficients (another early work but for lattices is
\cite{MR0198301}). Formal solutions to this problem on the whole real
line were presented in \cite{MR332065,MR332066,MR332067,MR332023},
where, in a first approximation, \cite{MR332067} deals with the case
with a zero potential sequence and off-diagonals rapidly converging to
$1$. \cite{guseinov76half,guseinov1976,guseinov1977} provide complete
solutions to the direct and inverse scattering problem for Jacobi
operators, on the whole and half line, being a finite first moment
perturbation of the discrete free Laplacian.

As in the case of equations with scalar coefficients, the theory of
difference equations with matrix coefficients is preceded by the
theory of Schrödinger operators with matrix coefficients (see
\cite{MR954906,MR1855088,MR1861473} and more recent contributions in
\cite{MR2894582,MR3221247}). The case of discrete Schrödinger
operators with matrix coefficients is studied in
\cite{MR4192212,MR4385984}.

There are not so many works devoted to the case of differential or
difference equations with operator coefficients
(with operators acting on a space which may be infinite dimensional).
The scattering theory
for difference operators with operator coefficients is studied in
\cite{MR0422923,MR1191552,mutlu}.

\subsection {Statements of the main results}\label{sec:statements-results}

The novelty of our work is that we study the stationary scattering
theory for operators associated with \eqref{eq:main-difference} with
operator valued coefficients $A_n$ and $B_n$. We will always assume
that the operators $A_n$ and $B_n$ are self-adjoint. Our results are
fourfold. After discussing the spectral properties of $\mathcal J$ as
a perturbation of $\mathcal J_{0}$, we establish existence of Jost
solutions and their properties, then we discuss their Wronskians and
conditions under which these solutions form a fundamental system.
Next we study the scattering matrix and its continuity and analyticity
properties. Finally, we discuss (non-)accumulation of eigenvalues of
the Jacobi operator. \medskip

Consider $\tau$ defined in \eqref{eq:main-difference} and
\begin{itemize}
   \item the \define{formal eigenvalue equation}
\begin{equation}\label{eq:eigenvalue}
    \tau u = \lambda u, \quad u=\{u_n\}_{n\in \mathbb Z} \subset H
\end{equation}
\item the \define{formal operator eigenvalue equation}
\begin{equation}\label{eq:eigenvalue-operator}
    \tau U = \lambda U, \quad U=\{U_{n}\}_{n\in\mathbb Z} \subset \mathscr{B}(H).
\end{equation}
\end{itemize}
A \define{formal vector solution} $u$ is a sequence in $H$ which
satisfies \eqref{eq:eigenvalue}, but does not necessarily lie in
$\ell^2(\mathbb Z,H)$. A \define{formal operator solution} $U$ is a
sequence in $\mathscr{B}(H)$ which satisfies
\eqref{eq:eigenvalue-operator}, but does not necessarily lie in
$\ell^2(\mathbb Z,\mathscr{B}(H))$. In
Section~\ref{sec:pert-discr-lapl} we will associate with $\tau$ a
bounded selfadjoint operator $\mathcal J$ on
$\mathcal H = \ell^2(\mathbb Z, H)$ under the assumption of the first
moment condition. We call $\lambda$ an \define{eigenvalue} if there
exists a non-zero $u\in \mathcal H$ such that \eqref{eq:eigenvalue}
holds. Such vector $u$ is an \define{eigenfunction of} $\mathcal J$.
Clearly, every eigenfunction $u$ is also a formal solution of
\eqref{eq:eigenvalue}. Any formal operator solution $U$ to
\eqref{eq:eigenvalue-operator} defines a formal vector solution by
setting $u=\{U_n v\}_{n\in\mathbb Z}$ for any $v\in H$. Analogously, any operator solution
$U\in \ell^2(\mathbb Z,\mathscr{B}(H))$ defines an eigenfunction $u$
of $\mathcal J$ by setting $u=\{U_n v\}_{n\in\mathbb Z}$ as before.

Let us take $z\in \mathbb C$ with $0 < |z|\le 1$ such that
$\lambda = z+z^{-1}$. It turns out that it is often more convenient to
work with $z$ rather than with $\lambda$, see
Example~\ref{ex:feeJacobi}. Since $U$ and $u$ are functions of $z$, we
use the notation $U(z)$ and $u(z)$, respectively. The analytic
properties of the solutions to \eqref{eq:eigenvalue-operator} as
functions of $z$ are crucial for this work. We will say more about the
relation between vector and operator solutions below in
Appendix~\ref{vector-operator-section}.

\subsubsection *{\underline{I. Existence and properties of the Jost solutions}} \ \medskip

\noindent Our first main result is the following theorem, which is a
combination of statements obtained in Theorem \ref{thm:jostsolutions} and
Proposition \ref{Jost-infinite-moment-prop}.
In what follows, $I$ denotes the identity operator.
\begin{theorem}
       \begin{enumerate}
      \item\label{item:JostFirstMomentCondition} 
Suppose that the $k$-th moment condition in Definition
      \ref{assumption:introduction} holds for some $k \geq 1$ and that $A_{n}$ is invertible with
      $A_{n}^{-1}\in\mathscr{B}(H)$ for any $n\in\mathbb{Z}$. If $z\in\mathbb C$ with
      $0<|z|\leq 1$, then there exist (formal) operator solutions $U^{\pm}(z)$ of
      \eqref{eq:eigenvalue} with $\lambda=z+z^{-1}$, called the \emph{Jost
      solutions} for the operator $\mathcal J$,  such that 
      \begin{equation}
	 \label{jost-estimate-second}
	 \begin{split}
	    U^+_n(z) &= z^n(I+o(n^{-k+1})),\ n\to\infty, \\
	    U^-_n(z) &= z^{-n}(I+o(n^{-k+1})),\ n\to-\infty.
	 \end{split}
      \end{equation}
      \item 
      If in addition the exponential moment condition in Definition~\ref{assumption:introduction} holds,
      then the Jost solutions of (\ref{item:JostFirstMomentCondition})
      are in fact well-defined on the larger punctured ball $0<|z| < e^{\varepsilon/2}$
      and satisfy the estimates \eqref{jost-estimate-second} for any $k\in\mathbb N_0$ there.
   \end{enumerate}
\end{theorem}

\subsubsection *{\underline{II. Jost solutions as a fundamental system and
    their Wronskians}} \ \medskip

\noindent Consider the Wronskian as introduced below in Definition
\ref{def:Wronskian}. Namely, for sequences $U = \{U_n\}_{n\in\mathbb Z}$ and
$V = \{V_n\}_{n\in\mathbb Z}$ of bounded operators on $H$, their
\emph{Wronskian} is a sequence of operators on $H$, defined as
   \begin{align*}
      W_n(U,V) := U_{n-1} A_{n-1} V_n  - U_{n} A_{n-1} V_{n-1}.
   \end{align*}
   We will show in Corollary \ref{cor:wronskianofsols} that for sequences
   $U(z) = \{U_n(z)\}_{n\in\mathbb Z}$ and
   $V(z) = \{V_n(z)\}_{n\in\mathbb Z}$, which are formal operator
   solutions to \eqref{eq:eigenvalue-operator} with $\lambda = z + z^{-1}$, the operator
   $$
W_{n}(U(\overline z)^*,V(z)) \in \mathscr{B}(H)
$$ 
does not depend on $n \in\mathbb Z$ and  we can thus omit the subscript. Now we
   can state our second main result, which combines Theorems
   \ref{thm:fundamentalsystem}, \ref{thm:fundamentalsystem-finite} and
   \ref{them:fundamentalsystem2}.
   \begin{theorem}
   Let $A_{n}$ be invertible with
  $A_{n}^{-1}\in\mathscr{B}(H)$ for any $n\in\mathbb{Z}$.
  Assume that 
  \begin{enumerate}[label={(\alph*)}]
     \item either the first moment condition holds and $|z|=1$,
     \item or the exponential moment condition holds and
     $e^{-\varepsilon / 2} < |z| < e^{\varepsilon / 2}$.
  \end{enumerate}
Consider any of the following pairs of Jost solutions
\begin{equation*}
  \{R(z),\widehat{R}(z)\}:=
  \begin{cases}
    \{U^{\pm}(z), U^{\pm}(z^{-1})\}, \\
    \{U^{\pm}(z), U^{\mp}(z^{-1})\}, \\
    \{U^{\pm}(z), U^{\mp}(z)\}.
  \end{cases}
\end{equation*}
  If the inverse $W(R(\bar z)^*, \widehat{R}(z))^{-1}$ exists and is in $\mathscr{B}(H)$,
  then the pair $\{R(z),\widehat{R}(z)\}$ forms a fundamental system of
  formal solutions, i.\,e., for a formal operator solution $V(z) = \{V_n(z)\}_{n\in\mathbb Z}$ to
  \eqref{eq:eigenvalue-operator} with $\lambda = z+z^{-1}$, there
  exists a unique pair of operator coefficients $P^{\pm}(z)$ and $Q^{\pm}(z)$ in $\mathscr{B}(H)$
  (independent of $n$) such that for any $n$
  \begin{equation*}
     V_{n}(z)=R_{n}(z) P(z) + \widehat{R}_{n}(z)Q(z).
  \end{equation*}
The same statement holds for formal vector solutions to
\eqref{eq:eigenvalue} with $P(z), Q(z) \in H$.

In case $\{R(z),\widehat{R}(z)\} = \{U^{\pm}(z), U^{\mp}(z)\}$ with the first moment condition,
the statement holds in a larger regime $|z| \leq 1$ if either $\dim H < \infty$ or if $B_n, A_n-I$
are compact for all $n \in \mathbb{Z}$. 
\end{theorem}

In Proposition \ref{prop:wronsks} below we establish that
\begin{equation}
    W( U^{\pm}(\bar z)^*, U^{\pm}(z^{-1})) = \pm(z^{-1}-z)I\,,\quad z \neq \pm 1
\end{equation}
for $|z| = 1$ when the first moment condition holds and for
$e^{-\varepsilon/2} < |z| < e^{\varepsilon/2}$ when the exponential moment
condition holds. Consequently, $\{U^{\pm}(z), U^{\pm}(z^{-1})\}$ form a
fundamental system for $z \neq \pm 1$ and $|z| = 1$ under the first moment, and for
$e^{-\varepsilon/2} < |z| < e^{\varepsilon/2}$ under the exponential moment
condition.

\subsubsection *{\underline{III. Scattering matrix and its continuity}} \ \medskip

\noindent The scattering matrix $\mathtt S(z)$ is a $2\times 2$ operator-valued matrix for which
\begin{equation}
   \begin{bmatrix}U^{-}(z),\ U^{+}(z^{-1})
   \end{bmatrix}
   = \begin{bmatrix}U^{+}(z),\ U^{-}(z^{-1})
   \end{bmatrix} \mathtt S(z).
\end{equation}
It is constructed explicitly in Proposition \ref{eq:scatteringmatrixdef} and
exists under the second moment condition \emph{a priori} for $|z|=1$ and
$z \neq \pm 1$. Our third main result is its continuity in the operator norm at the points of
exclusion $z = \pm 1$, as asserted below in Theorem \ref{thm:contscat}.

\begin{theorem}
  Suppose that $A_{n}$ is invertible with $A_{n}^{-1}\in\mathscr{B}(H)$ for
  any $n\in\mathbb{Z}$. If the second moment condition in Definition
  \ref{assumption:introduction} holds, and the Wronskian $W(U^{+}(z_0)^*, U^{-}(z_0))$
  has closed range for $z_0= \pm 1$, then the scattering matrix
  $\mathtt S(z)$, defined initially for $|z|=1$ but $z\neq \pm 1$, has a
  continuous extension to the entire unit circle $|z|=1$.
\end{theorem}

\subsubsection *{\underline{IV. Accumulation of eigenvalues of the Jacobi
    operator}} \ \medskip

Our next main results are concerned with finiteness of the discrete spectrum of $\mathcal J$.
We shall assume for the discussion of the spectrum that
all $B_n$ and $A_n-I$ are compact and the zero-th moment condition holds
$$
\sum\limits_{n\in\mathbb Z} \|A_n - I\| + \|B_n\| < \infty.
$$
Proposition \ref{prop:Jspectrum} 
asserts in this case that the essential spectrum of $\mathcal J$ as an operator on 
$\ell^2(\mathbb Z,H)$ is 
 \begin{align*}
      \sigma_{ess}(\mathcal J) = [-2,2].
  \end{align*}
The operator $\mathcal J$ is bounded and selfadjoint and hence its discrete spectrum
$\sigma_{d}(\mathcal J) = \sigma (\mathcal J) \backslash \sigma_{ess}(\mathcal J)$
is bounded and can only accumulate at the essential spectrum $\sigma_{ess}(\mathcal J)= [-2,2]$.
Our next main result rules out accumulation of the discrete spectrum under certain conditions, 
particularly in the case $\dim H = \infty$.

\begin{theorem}  Suppose that $A_{n}$ is invertible with $A_{n}^{-1}\in\mathscr{B}(H)$ 
and the operators $A_n-I$ and $B_n$ on $H$ are compact for all $n \in \mathbb Z$.   Assume
 the third moment condition of Definition \ref{assumption:introduction} holds and 
     $$
     W(U^+(1)^*, U^-(1))\ \textup{and}\  W(U^+(-1)^*, U^-(-1)) \  \textup{have closed range.}
     $$
Then the discrete spectrum of $\mathcal J$ 
has no accumulation points and hence is finite.  
\end{theorem}

This result is proved below in see Theorem \ref{matrix-case-non-acc}.

\begin{remark} In \cite{MR4385984}, the authors have proven non-accumulation of the spectrum in the case 
$\dim H < \infty$,  assuming only the first moment condition. 
Our theorem here requires stronger conditions,  but allows $H$ to be infinite-dimensional.
\end{remark}

Our final main result is concerned with the case where only the zero-th moment condition is
imposed.  In this case we cannot expect to have non-accumulation of
eigenvalues in general.  However, we have control on the rate of accumulation.

\begin{theorem}
   Assume that the operators $(A_n-I)$ and $B_n$ on $H$ are trace class for
   all $n \in \mathbb Z$ and that
   \begin{align*}
     2 \sum_{n\in\mathbb Z}  \|A_n - I\|_1 + \sum_{n\in\mathbb Z}  \|B_n\|_1  < \infty,
   \end{align*}
   where $\|\cdot \|_1$ denotes the trace-norm in $H$.
   Consider any $R<1$ and enumerate the finite set of eigenvalues
   $\{ \lambda_j = z_j + z_j^{-1}\}_{j =1}^{J(R)}$ of $\mathcal J$, with $|z_j| < R$, counted
   with their multiplicities.  Then
\begin{equation}
   \begin{split}
      \prod_{j =1}^{J(R)} R |z_j|^{-1} &\leq
      \exp \left( \left| \frac{R}{1-R^2} \right| 
      \left( 2 \sum_{n\in\mathbb Z}  \|A_n - I\|_1 + \sum_{n\in\mathbb Z}  \|B_n\|_1 \right) \right), \\
      J(R) &\leq \frac{1}{\ln R} \left| \frac{R}{1-R^2} \right| 
      \left( 2 \sum_{n\in\mathbb Z}  \|A_n - I\|_1 + \sum_{n\in\mathbb Z}  \|B_n\|_1 \right).
   \end{split}
\end{equation}
\end{theorem}
\noindent This result is proved below in Theorem \ref{resonance-case-non-acc}.
\bigskip

\noindent \emph{Acknowledgements.}
The authors are grateful to the Banff International Research Station and Casa
Matemática Oaxaca for organizing a workshop on "Analytic and Geometric
Aspects of Spectral Theory" (22w5149), where the project was initiated. We are also
grateful to Gerald Teschl and Abdon Choque for insightful discussions. 
The first author was supported by an URC grant from DePaul University.
The second author is
supported by CONAHCyT CF-2019 304005. Part of this work was done during a
visit by the second author to the University of Oldenburg. He is grateful for their kind
hospitality.
The last author gratefully acknowledges the support of Facultad de Ciencias, Universidad de los Andes, Bogot\'a, INV-2023-162-2824.

\section {Perturbed discrete Laplacian and its continuous spectrum}\label{sec:perturbed-laplacian}

\subsection {Compact, Hilbert-Schmidt and trace class operators}\label{trace-class-subsection }\medskip

Let $p\ge 1$. Recall that a compact operator $K$ acting on a Hilbert space
$H$ is said to belong to the $p$-Schatten-von Neumann class if
$\sum_{n=1}^\infty s_{n}^p < \infty$ \cite{birman-solomjak}. Here, the $s_n$
are the \define{singular values} of $K$, that is, the square roots of the
eigenvalues of the nonnegative operator $K^*K$. The set of all such operators
is denoted by $\mathscr{S}_p(H)$. Then
$$
\|K\|_p := \left(\sum_{n=1}^\infty s_{n}^p\right)^{1/p}
$$ 
defines a norm on $\mathscr{S}_p(H)$. Note that
$\mathscr{S}_q(H) \subset \mathscr{S}_p(H)$ and that $\|K\| \le \|K\|_p \le \|K\|_{q}$ if $1 \le q \le p$. 
Moreover, the Schatten-von Neumann classes
are two-sided ideals in $\mathscr{B}(H)$. The operators in $\mathscr{S}_1(H)$
are called \define{trace class operators}, while operators in
$\mathscr{S}_2(H)$ are called \define{Hilbert-Schmidt operators}.
Recall that $\mathcal{H} = \ell^2(\mathbb Z, H)$.

\begin{definition}
   \label{def:diagoperator}
   Let $\{C_n\}_{n\in\mathbb Z}  \subset \mathscr{B}(H)$ such that $\sup\{
   \|C_n\| : n\in\mathbb Z\} < \infty$.  Define the block diagonal operator
   \begin{align*}
      \mathcal C := \diag[ \{C_n\}_{n\in\mathbb Z}] : \mathcal H \to \mathcal H,\qquad
      (\mathcal C u)_n = C_n u_n
      \quad \text{ for } u = \{u_n\}_{n\in\mathbb Z}\in \mathcal H.
   \end{align*}
\end{definition}
The next lemma shows that $\mathcal C$ defines a bounded operator on
$\mathcal{H}$.
\begin{lemma}
   \label{lem:diagnorm}
   Under the condition in Definition~\ref{def:diagoperator}, the operator
   $\mathcal C$ is a bounded operator on the Hilbert space $\mathcal{H}$ with
   the operator norm $\|\mathcal C\| = \sup\{ \|C_n\| : n\in\mathbb Z\}$.
\end{lemma}
\begin{proof}
  Let $c := \sup\{ \|C_n\| : n\in\mathbb Z\}$. Then for
  $u = \{u_n\}_{n\in\mathbb Z}\in \mathcal H$ we have that
   \begin{align*}
      \| \{ C_n u_n \}_{n\in\mathbb Z} \|^2
      = \sum_{n\in\mathbb Z} \| C_n u_n \|^2
      \le c^2 \sum_{n\in\mathbb Z} \| u_n \|^2
      = c^2 \| u \|^2.
   \end{align*}
   This implies that $\|\mathcal C \| \le c$. To show the other inequality,
   note that for every $\epsilon > 0$ there are $n_\epsilon\in\mathbb Z$ and
   $x_\epsilon\in H$ such that $\|x_\epsilon \| = 1$ and
   $\|C_{n_\epsilon} x_\epsilon \| \ge c - \epsilon$. If we define
   $u_\epsilon \in \mathcal H$ such that $(u_\epsilon)_{n} = x_\epsilon$ if
   $n = n_\epsilon$ and $0$ else, then clearly $\|u_\epsilon\| = 1$ and
   $\|\mathcal C u_\epsilon \| = \| C_{n_\epsilon} x_\epsilon\| \ge c -
   \epsilon$.
\end{proof}
We now identify conditions which guarantee that $\mathcal C$ is compact,
Hilbert-Schmidt or trace class on $\mathcal H$.

\begin{lemma} \label{lem:diagcompact} Let
  $\{C_n\}_{n\in\mathbb Z} \subset \mathscr{B}(H)$ be a 
sequence of compact operators on $H$ and let
  $\mathcal C := \diag[ \{C_n\}_{n\in\mathbb Z}] $
be the corresponding block diagonal operator.
\begin{enumerate}
\item If\ $\sum\limits_{n\in\mathbb Z} \|C_n\| < \infty$, then $\mathcal C$ is
  compact on $\mathcal H$.
\item If\ $\sum\limits_{n\in\mathbb Z} \|C_n\|_p < \infty$ for some%
\footnote{In particular $C_n \in \mathscr{S}_p(H)$ for all $n \in \mathbb Z$.}
$p\ge 1$, then $\mathcal C \in \mathscr{S}_p(\mathcal H)$. Moreover
$$
\|\mathcal C\|_p \le \sum_{n\in\mathbb Z} \|C_n\|_p.
$$
\end{enumerate}
\end{lemma}
\begin{proof}
  Let $\mathcal C_n$ be the diagonal operator whose $n$-th component is equal
  to $C_n$ and all others are equal to $0$.
   \begin{enumerate}
   \item Clearly each $\mathcal C_n$ is a compact operator on $\mathcal H$
     and $\|\mathcal C_n\| = \|C_n\|$. Moreover,
     $\mathcal C = \sum_{n\in\mathbb Z} \mathcal C_n$ where the sum converges
     in the operator norm by assumption. Hence $\mathcal C$ is the norm limit
     of compact operators, therefore it is compact.
   \item It is easy to see that each $\mathcal C_n$ belongs to the
     $p$-Schatten-von Neumann class on $\mathcal H$ and that
     $\|\mathcal C_n\|_p = \|C_n\|_p$. Hence, by assumption,
     $\mathcal C = \sum_{n\in\mathbb Z} \mathcal C_p$ converges in the $p$-th
     Schatten norm and therefore it is itself in $\mathscr{S}_p(\mathcal H)$.
     \qedhere
   \end{enumerate}
\end{proof}

\subsection {Perturbed discrete Laplacian as a difference
  operator}\label{sec:pert-discr-lapl}
Before defining the main operators to be considered in this work, let us
introduce the so-called \define{left shift operator}:
\begin{align}
  \label{eq:def-shift}
      \mathcal S: \mathcal H \to \mathcal H,\qquad
      (\mathcal Su)_n = u_{n+1}
      \quad \text{ for } u = (u_n)_{n\in\mathbb Z}\in \mathcal H.
\end{align}

  \begin{remark}\label{Snorm}
     Clearly, $\mathcal S$ is bounded with $\|\mathcal S\| = 1$ and its
     adjoint operator is the block right shift operator, given by
     \begin{align*}
	\mathcal S^*: \mathcal H \to \mathcal H,\qquad
	(\mathcal S^*u)_n = u_{n-1}
	\quad \text{ for } u = (u_n)_{n\in\mathbb Z}\in \mathcal H.
     \end{align*}
     Moreover,  $\mathcal S\mathcal S^* = \mathcal S^*\mathcal S = I$
     (the identity operator on $H$),
     hence $\mathcal S$ is a unitary operator on $\mathcal H$. Note that
     in the case of infinite but not bi-infinite sequences 
     {\upshape(}that is $\ell^2(\mathbb N, H)$ instead of $\ell^2(\mathbb Z, H)${\upshape)}
     the shift
     operator $\mathcal{S}$ is not unitary.
  \end{remark}

  Let us now define the main object of our considerations in this paper.
  \begin{definition}
     \label{def:block-jacobi-operator}
     Let
     $\{A_{n}\}_{n\in\mathbb Z}, \{B_{n}\}_{n\in\mathbb Z} \subset
     \mathscr{B}(H)$ be sequences of bounded self-adjoint operators on $H$, such
     that $\sup\{\|A_n\|, \|B_n\| : n\in\mathbb Z\}<\infty$. Then the diagonal
     operators $\mathcal A := \diag[ \{A_n\}_{n\in\mathbb Z} ]$ and
     $\mathcal B := \diag[ \{B_n\}_{n\in\mathbb Z} ]$ define bounded
     self-adjoint operators on $\mathcal H$ by Lemma \ref{lem:diagnorm}. The
     \define{(block) Jacobi operator} $\mathcal J$ is a self-adjoint bounded
     operator on $\mathcal H$, defined by
     \begin{equation}
	\label{eq:def-block-jacobi-operator}
	\mathcal J : = \mathcal S^* \mathcal A + \mathcal B + \mathcal A \mathcal S.
     \end{equation}
  \end{definition}
  The operator $\mathcal J$ is the closure of the operator generated by the
  difference expression $\tau$ given in \eqref{eq:main-difference} acting on
  finite sequences $u=\{u_{n}\}_{n\in\mathbb{Z}}$ with $u_{n}\in H$ for
  any $n\in\mathbb{Z}$. We henceforth say that $\mathcal J$ is the operator
  generated by $\tau$ or, equivalently, that the bi-infinite block matrix
  \begin{align}
  \label{eq:matrix-block-jacobi}
   \begin{tikzpicture}[baseline=(current bounding box.center),
      every left delimiter/.style={xshift=.5em},
      every right delimiter/.style={xshift=-.5em},
      ampersand replacement=\&,
      ]
      \matrix (m) [matrix of math nodes,nodes in empty cells, right
      delimiter={)},left delimiter={(} ]{
        \qquad \&             \&  \& \&  \&  \& \&   \\[2ex]
        \&  A_{-2}     \& B_{-1}\& A_{-1}\&                \&   \\
        \&  \&  A_{-1} \& B_{0} \& A_{0} \&                \&   \\
        \&  \&         \& A_{0} \& B_{1} \& A_{1}          \&   \\
        \&  \&         \&       \& A_{1} \& B_{2} \& A_{2} \&   \\[2ex]
        \&  \&         \&       \&       \&       \&       \& \qquad  \\
      } ; \draw[thick, dotted] ($(m-2-2.north
      west)+(0ex,0ex)$) -- ++ (-2ex,2ex); \draw[thick, dotted]
      ($(m-2-3.north
      west)+(0ex,0ex)$) -- ++ (-2ex,2ex); \draw[thick, dotted]
      ($(m-2-4.north west)+(0ex,0ex)$) -- ++ (-2ex,2ex);
      \draw[thick, dotted] ($(m-5-5.south east)+(0ex,0ex)$)-- ++ (2ex,-2ex);
      \draw[thick, dotted] ($(m-5-6.south east)+(0ex,0ex)$)-- ++ (2ex,-2ex);
      \draw[thick, dotted] ($(m-5-7.south east)+(0ex,0ex)$)-- ++ (2ex,-2ex);
   \end{tikzpicture}
   \end{align}
   is the block matrix representation of $\mathcal{J}$. The operator
   $\mathcal J$ is related to the block matrix
   \eqref{eq:matrix-block-jacobi} in the same way that a closed
   symmetric operator is related to its matrix representation with
   respect to some basis (see \cite[Sec.47]{MR1255973}). In this
   sense, we allow ourselves to refer to
   \eqref{eq:matrix-block-jacobi} as a matrix representation of
   $\mathcal J$, where, instead of a basis one has a decomposition of
   the space $\mathcal H=\oplus_{n\in \mathbb{Z}}H_{n}$, where
   $H_{n}=H$ for all $n\in\mathbb{Z}$ (see \cite[Chap.\,3 Sec.\,5]{birman-solomjak}).

Alongside the Jacobi operator, we consider the \define{discrete free block Laplacian}
(or free block Jacobi operator)
   \begin{equation}
     \label{eq:def-block-free-laplacian}
     \mathcal{J}_{0}:=\mathcal S^* +\mathcal S\,.
   \end{equation}
The block matrix representation of $\mathcal J_{0}$ is
  \begin{align}
   \label{eq:matrix-block-free-jacobi}
   \begin{tikzpicture}[baseline=(current bounding box.center),
      every left delimiter/.style={xshift=.5em},
      every right delimiter/.style={xshift=-.5em},
      ampersand replacement=\&,
      ]
      \matrix (m) [matrix of math nodes,nodes in empty cells, right
      delimiter={)},left delimiter={(},
      column sep =1.5ex
      ]{
        \&             \&  \& \&  \&  \& \&   \\[2ex]
        \&  I \& 0 \& I \&   \&   \&   \&   \\
        \&    \& I \& 0 \& I \&   \&   \&   \\
        \&    \&   \& I \& 0 \& I \&   \&   \\
        \&    \&   \&   \& I \& 0 \& I \&   \\[2ex]
        \&    \&   \&   \&   \&   \&   \&  \\
      } ; \draw[thick, dotted] ($(m-2-2.north
      west)+(0ex,0ex)$) -- ++ (-2ex,2ex); \draw[thick, dotted]
      ($(m-2-3.north
      west)+(0ex,0ex)$) -- ++ (-2ex,2ex); \draw[thick, dotted]
      ($(m-2-4.north west)+(0ex,0ex)$) -- ++ (-2ex,2ex);
      \draw[thick, dotted] ($(m-5-5.south east)+(0ex,0ex)$)-- ++ (2ex,-2ex);
      \draw[thick, dotted] ($(m-5-6.south east)+(0ex,0ex)$)-- ++ (2ex,-2ex);
      \draw[thick, dotted] ($(m-5-7.south east)+(0ex,0ex)$)-- ++ (2ex,-2ex);
   \end{tikzpicture}.
\end{align}

Note that the block matrix representation for $\mathcal S$ is a matrix
where there are zeros everywhere except on the diagonal above the main
diagonal, where the block entries are equal to the identity operator $I$.

Formal solutions to \eqref{eq:eigenvalue-operator} for the difference
expression generated by $\mathcal J_{0}$  are given in the following example.
\begin{example}
   \label{ex:feeJacobi}
   Let $A_{n}=I$, $B_{n}=0$ for all
   $n\in\mathbb{Z}$. This defines the free Jacobi operator $\mathcal J_0$.
   A direct computation
   shows that the sequences $\{z^nI\}_{n=-\infty}^{\infty}$ and
   $\{z^{-n}I\}_{n=-\infty}^{\infty}$ are solutions to
   \eqref{eq:eigenvalue-operator} for $\mathcal J_0$ with $\lambda= z + z^{-1}$. 

   Note that the
   \define{Zhukovsky function} $\lambda(z) = z + z^{-1}$ extends to a bijective mapping
   of the open unit disk onto
   $\overline{\mathbb C}\setminus[-2,2]$, where $\overline{\mathbb C}$
   is the one point compactification at infinity of the complex plane. 
   Moreover, the sets $\{ |z|=1, \im(z) > 0 \}$ and $\{ |z|=1, \im(z) < 0 \}$ are mapped bijectively onto $(-2,2)$ 
   and $\lambda(z)=\pm 2$ if and only if $z=\pm 1$.
\end{example}

\begin{lemma}
  \label{lem:normfreeJacobi}
  If $\mathcal J$ and $\mathcal J_{0}$ are the operators given by
  \eqref{eq:def-block-jacobi-operator} and
  \eqref{eq:def-block-free-laplacian}, respectively, then
  \begin{equation}
     \|\mathcal J\|\le 2\sup_{n\in\mathbb Z}\{\|A_n\|\} +  \sup_{n\in\mathbb Z}\{\|B_n\|\} 
     \qquad\text{and}\qquad
     \|\mathcal J_0\| = 2.
  \end{equation}
\end{lemma}
\begin{proof}
   We obtain from \eqref{eq:def-block-jacobi-operator} and Remark \ref{Snorm}
   that
   \begin{align*}
      \|\mathcal J \| &\le \|\mathcal S^* \mathcal A \| + \|\mathcal B \| + \|\mathcal A \mathcal S \|
      \le 2 \|\mathcal A \| + \|\mathcal B \|
      = 2\sup\{\|A_n\|\} + \sup\{\|B_n\|\}.
   \end{align*}
   This readily implies that $\|\mathcal J_0\| \le 2$. To show equality,
fix a normalized vector $e\in H$.
For $q\in\mathbb C$ with
$|q| < 1$ we define $x_q = \{(x_q)_n\}_{n \in \mathbb Z} \in \mathcal H$ by
$(x_q)_n := q^{|n|} e$. One computes $\|x_q\|^2 = \frac{1+|q|^2}{1-|q|^2}$.
Moreover
\begin{equation*}
(\mathcal J_0 u_q)_n = \left\{
\begin{split}
  &\, 2 q e, &n=0, \\
  &\, (q^{n-1} + q^{n+1}) e, &n > 0, \\
  &\, (q^{-n-1} + q^{-n+1}) e, &n < 0.
\end{split} \right.
\end{equation*}
Hence $\|\mathcal J_0 x_q\|^2 = \frac{2+8|q|^2-2|q|^4}{1-|q|^2}$ and therefore
$$
\|\mathcal J_0\| \ge \frac{ \|\mathcal J_0 x_q\| }{\|x_q\|} =
\left(\frac{2+8|q|^2-2|q|^4}{1+|q|^2}\right)^{1/2}.
$$
Letting $|q|$ tend to $1$, we see that $\|\mathcal J_0\|\ge 2$.
\qedhere
\end{proof}

\begin{remark}When one or both of the sequences of operators $\{ A_{n} \}_{n\in\mathbb Z}$ and
   $\{B_{n \}_{n\in\mathbb Z}}$ are not uniformly bounded the notion of matrix representation is no
  longer so straightforward. In this case, one should show that $\mathcal J$
  is the minimal closed operator (see the details in
  \cite[Sec.\,47]{MR1255973}) generated by the matrix. \end{remark}

\subsection {Spectrum of the free Jacobi operator}

Our next goal is to find the spectrum of $\mathcal J_0$.
To this end,  we first investigate the spectrum of the shift operator $\mathcal S$.
We use conventions as specified in Appendix~\ref{spectrum-convention-section}.

\begin{lemma}
If $\mathcal S$ is the shift operator defined in \eqref{eq:def-shift}, then
\begin{align*}
  \sigma( \mathcal S ) = \sigma_{cont}( \mathcal S ) 
  = \{ |z| = 1 \},
  \qquad
  \sigma_p( \mathcal S ) = \sigma_{pp}( \mathcal S ) = \varnothing\,,
\end{align*}
where
$\sigma_{cont}(\mathcal S) = \sigma_{sc}(\mathcal S) \cup \sigma_{ac}(\mathcal S)$
stands for the continuous spectrum%
\footnote{This is in general different from continuous spectrum as defined for instance in \cite{KolmogorovFomin} or \cite{TaylorLay}.
There the continuous spectrum consists of all $\lambda\in\mathbb C$ such that 
$\range(\mathcal S-\lambda) \neq \overline{\range(\mathcal S-\lambda)} = H$.
See Appendix~\ref{spectrum-convention-section}.
}
and $\sigma_{pp}(\mathcal S)$ is the pure
point spectrum of $\mathcal S$, see for instance, \cite[\S\, X.1]{kato}.
\end{lemma}
\begin{proof}
  Since $\mathcal S$ is a unitary operator, it follows that
  $\sigma(\mathcal S)\subseteq \{|z| = 1\}$. It is easy to see that
  $\sigma_p(\mathcal S)=\varnothing$, hence also
  $\sigma_{pp}(\mathcal S)=\varnothing$ and
  $\sigma( \mathcal S ) = \sigma_{cont}( \mathcal S )$.
  Let $z\in\mathbb C$ with $|z|=1$. We will show that $\mathcal S - z$ is
  not surjective, which means that $z$ cannot be in the resolvent set. 
  Let $e \in H$ be a vector of norm $1$. Define
  $y = \{y_n\}_{n \in \mathbb Z} \in \mathcal H$ by
  $y_n := \frac{z^n}{1+|n|} e$. Consider the equation
  $(\mathcal S - z)x = y$. It can be written as a recursive relation
  $x_{n+1} - z x_n = y_n$. 
  Note that each $x_n$ can be written in a unique way as $x_n = \xi_n e + v_n$ where $\xi_n, \eta_n\in\mathbb C$ and $v_n\in \Span\{e\}^\perp$.
  Since $y_n\in\Span\{e\}$ for every $n\in\mathbb Z$, it follows that $v_{n+1} - z v_n = 0$. 
  Therefore $v_n=0$ for all $n\in\mathbb Z$ because otherwise $z$ would be an eigenvalue of $\mathcal S$.
  Hence we obtain the recursion 
  $\xi_{n+1} - z \xi_n = \frac{z^n}{1+|n|}$ for the coefficients $\xi_n$.
  For any given $\xi_0\in \mathbb C$, this recursion
  has the unique solution which can be shown by induction to be
  \begin{align*}
     \xi_n = \left\{
     \begin{array}{ll}
	\xi_0, &n=0, \\
	z^n \xi_0 + z^{n-1} \sum\limits_{k=1}^n \frac{1}{k}, &n \geq 1, \\
	z^n \xi_0 - z^{n-1} \sum\limits_{k=2}^{|n|+1} \frac{1}{k}, &n \leq -1.
     \end{array} \right.
  \end{align*}
  Clearly, $x = \left\{ \xi_n e \right\}_{n\in\mathbb N} \notin \mathcal H = \ell^2(\mathbb Z,H)$, hence $\mathcal S-z$ is
  not surjective as a bounded operator on $\mathcal H$ and therefore
  $\sigma( \mathcal S ) = \{ |z| = 1 \}$.
\end{proof}

\begin{proposition}
   \label{prop:J0spectrum}
   Let $\mathcal J_0$ be the free Jacobi operator.
   Then
   \begin{align*}
     \sigma( \mathcal J_0 ) = \sigma_{cont}( \mathcal{J}_{0} ) =
     \sigma_{ess}( \mathcal J_{0}) = [-2,\, 2],
     \qquad
     \sigma_p( \mathcal J_{0} ) = \sigma_{pp}( \mathcal J_{0} ) = \varnothing.
   \end{align*}
\end{proposition}
\begin{proof}
  Let $\lambda(z) = z + z^{-1}$ 
  for $z\in\mathbb C\setminus\{ 0\}$ be the Zhukovsky function.
  Then $\lambda$ is a
  continuous function on $\sigma(\mathcal S) = \{ |z| = 1 \}$ and
  $\mathcal J_0 = \mathcal S + \mathcal S^{-1} = \lambda(\mathcal S)$. Hence the
  spectral mapping theorem for bounded normal operators shows that
   \begin{align*}
      \sigma(\mathcal J_0)
      = \lambda(\sigma(\mathcal S))
      = \{ z+z^{-1} : z\in \sigma(\mathcal S) \}
      = \{ z+ \overline z : |z|=1 \} = [-2,2].
   \end{align*}
   Since also $\sigma_p(\mathcal J_0) = \lambda(\sigma_p(\mathcal S)) = \varnothing$,
   the claim is proved.
\end{proof}
\begin{remark}\label{absolutely-cts-remark}
   Note that the proposition above gives an alternative proof for $\|\mathcal J_0\|=2$, see Lemma~\ref{lem:normfreeJacobi}.
  Since $\mathcal J_0$ is a bounded self-adjoint operator, we have that 
  $\|\mathcal J_0\| = \sup\{ |\lambda| : \lambda\in\sigma(\mathcal J_0)\} = 2$. 
  Moreover, along the same lines of reasoning given in
  \cite[Sec.\,XIII.6]{MR529429}, one can show that the spectrum of
  $\mathcal J_{0}$ is purely absolutely continuous, i.e.  the singular spectrum is empty.
\end{remark}

\begin{remark}
  If $\mathcal B = 0$ in \eqref{eq:def-block-jacobi-operator}, then the block
  Jacobi operator $\mathcal J$ is unitarily equivalent to $-\mathcal J$; in
  particular its spectrum is symmetric with respect to $0$. This can be seen
  if we define $\mathcal V = \diag(\dots, I, -I, I, -I, \dots)$. Clearly,
  $\mathcal V$ is a unitary operator on the Hilbert space $\mathcal H$ with
  $\mathcal V^{-1}=\mathcal V$. Moreover, $\mathcal V$ commutes with any
  block diagonal operator, in particular with $\mathcal A$, and
  $\mathcal V \mathcal S = - \mathcal S \mathcal V$. Hence
  $\mathcal V^{-1} \mathcal J \mathcal V = -\mathcal J$.
\end{remark}

It is possible to give a bit more precise information on the spectrum of
$\mathcal J_0$. To explain this let us consider a generic and useful
construction for illustrations.
\begin{example}
  \label{ex:diagonal-blocks}
  Suppose that $H_{n}=\ell^2(\mathbb Z)$ for all $n\in\mathbb{Z}$ and as before let 
  \begin{equation*}
    \mathcal H=\bigoplus_{n\in\mathbb{Z}} H_{n}\,.
  \end{equation*}
  Consider
  the sequence of sequences $\{e_{n}(k)\}_{k,n\in\mathbb{Z}}$, with
  $k,n\in \mathbb{Z}$, such that, for any fixed $n \in \mathbb{Z}$,
  $\{e_{n}(k)\}_{k\in\mathbb{Z}}$ is the canonical basis in
  $H_{n}=\ell^2(\mathbb Z)$. Let the
  bounded operators $A_{n}$ and $B_{n}$ be such that their matrix
  representations for any $n\in\mathbb{Z}$ (see
  \cite[Sec.\,47]{MR1255973}) with respect to the basis
  $\{e_{n}(k)\}_{k\in\mathbb{Z}}$ are infinite diagonal matrices, that
  is $A_{n}=\diag_{k\in\mathbb{Z}}\{a_{n}(k)\}$,
  $B_{n}=\diag_{k\in\mathbb{Z}}\{b_{n}(k)\}$. Then, under the
  decomposition of the space
  \begin{equation*}
    \mathcal{H}=\dots\oplus\Span_{n\in\mathbb{Z}}\{e_{n}(-1)\}
    \oplus\Span_{n\in\mathbb{Z}}\{e_{n}(0)\}\oplus\Span_{n\in\mathbb{Z}}
    \{e_{n}(1)\}\oplus\dots\,,
  \end{equation*}
  the operator $\mathcal{J}$ becomes an infinite orthogonal sum of
  scalar Jacobi operators $J_{k}$, that is,
  \begin{equation}
    \label{eq:orthogonal-decomposition-matrix}
    \mathcal J=\bigoplus_{k\in\mathbb{Z}}J_{k}\,,
  \end{equation}
  where the main diagonal of the matrix representation of
  $J_{k}$ is the sequence $\{b_{n}(k)\}_{n\in\mathbb{Z}}$
  and its off-diagonals are given by $\{a_{n}(k)\}_{n\in\mathbb{Z}}$.
  To verify that the infinite orthogonal sum in the r.h.s of
  \eqref{eq:orthogonal-decomposition-matrix} is indeed the operator
  $\mathcal J$, one refers to the theory of direct integrals of Hilbert
  spaces \cite[Chap.\,7]{birman-solomjak}.

  A realization of the construction above is the free block Jacobi
  operator $\mathcal{J}_{0}$ which turns out to be an infinite
  orthogonal sum of identical instances of the free discrete scalar
  Laplacian. Thus, taking into account the spectral characterization
  of the scalar discrete free Laplacian \cite[Sec.1.3]{teschl}, one concludes
  that the spectrum of $\mathcal{J}_{0}$ is the set $[-2,2]$, it is
  purely absolutely continuous (we already knew that) and it has
  infinite multiplicity (for the notion of multiplicity of the
  continuous spectrum, see \cite[Sec.\,71]{MR1255973}).
\end{example}
\subsection {Trace class and Hilbert-Schmidt perturbations}
\label{sec:trace-class-hilbert}
\begin{proposition}
\label{prop:Jspectrum}
Let $\mathcal J$ and $\mathcal J_{0}$ be the operators given by
  \eqref{eq:def-block-jacobi-operator} and
  \eqref{eq:def-block-free-laplacian}, respectively.
If all $B_n$ and $A_n-I$ are compact, then we have the following.
\begin{enumerate}
\item If\; $\sum\limits_{n\in\mathbb Z} \|A_n - I\| + \|B_n\| < \infty$, then
  $\mathcal J$ is a compact perturbation of $\mathcal J_0$ and
   \begin{align*}
      \sigma_{ess}(\mathcal J) = \sigma_{ess}(\mathcal J_0) = [-2,2].
   \end{align*}

 \item If\; $\sum\limits_{n\in\mathbb Z} \|A_n - I\|_1 + \|B_n\|_1 < \infty$,
   then $\mathcal J$ is a trace-class perturbation of $\mathcal J_0$ and
   \begin{align*}
      \sigma_{cont}(\mathcal J) = \sigma_{ac}(\mathcal J_0) = [-2,2],
   \end{align*}
where presence of $\sigma_{sc}(\mathcal J)$ is not ruled out \emph{a priori}.
\end{enumerate}
\end{proposition}
\begin{proof}
  By Definition~\ref{def:block-jacobi-operator}, we can write
   \begin{align*}
     \mathcal J = \mathcal J_0 +  \mathcal S^* (\mathcal A - \mathcal I) +
     \mathcal B + (\mathcal A - \mathcal I) \mathcal S
   \end{align*}
   where $\mathcal I$ denotes the identity operator on $\mathcal H$. Recall
   that $\sigma_{ess}(\mathcal J_0) = \sigma_{cont}(\mathcal J_0) = [-2,2]$ by
   Proposition~\ref{prop:J0spectrum}.  
   Moreover,  as noted in the Remark~\ref{absolutely-cts-remark}, $\sigma_{sc}(\mathcal J_0) = \varnothing$ and hence
$\sigma_{ac}(\mathcal J_0) = \sigma_{cont}(\mathcal J_0) = [-2,2]$.
\begin{enumerate}
\item Under our assumptions, Lemma~\ref{lem:diagcompact} shows that
  $\mathcal B$ and $\mathcal A -\mathcal I$ are compact. Since $\mathcal S$
  is bounded, it follows that $\mathcal J - \mathcal J_0$ is compact. Hence,
  by Weyl's theorem, $\sigma_{ess}(\mathcal J) = \sigma_{ess}(\mathcal J_0)= [-2,2]$.

\item Again, Lemma~\ref{lem:diagcompact} shows that $\mathcal B$ and
  $\mathcal A -\mathcal I$ are trace class operators. Since
  $\mathcal S$ is bounded and the trace class operators are an ideal
  in $\mathscr{B}(H)$, it follows that $\mathcal J - \mathcal J_0$ is
  a trace class operator. Hence, by the Kato-Rosenblum theorem
  \cite{katorosenblum-kato, katorosenblum-rosenblum},
  $\sigma_{ac}(\mathcal J) = \sigma_{ac}(\mathcal J_0) = [-2,2]$.
  Finally, note that, from the first statement, there is no singular
  continuous spectrum outside the interval $[-2,2]$. Thus, the second
  statement follows. \qedhere
\end{enumerate}
\end{proof}

We conclude with two examples concerning eigenvalues outside and
inside the absolutely continuous spectrum. To this end, assume that
the operators $A_{n}$ and $B_{n}$ are diagonal for each
$n\in\mathbb{Z}$ as in Example~\ref{ex:diagonal-blocks} and consider
the notation introduced there.

\begin{example}
  Assume that  $a_{n}(k)=a_{n}$ and $b_{n}(k)=b_{n}$ for all $k\in\mathbb{Z}$
  (that is, $A_{n}$ and $B_{n}$ are pseudoscalars). If $a_{n}>0,\
  b_{n}\in\mathbb{R}$ for all
  $n\in\mathbb Z$ and
  \begin{equation*}
  \sum_{n\in\mathbb Z} |n|(|a_{n}-1|+|b_{n}|)<+\infty,
  \end{equation*}
  (see  \cite[Eq.\,2]{guseinov1976}), then
  $\sigma_{ess}(\mathcal J)=\sigma(J_{k})$ for all $k\in\mathbb{Z}$,
  where $\sigma(\mathcal J)\cap [-2,2]$ is absolutely continuous
  spectrum with infinite multiplicity and, additionally, there is a
  finite number of eigenvalues, each of infinite multiplicity. Note
  that in this case, Proposition \ref{prop:Jspectrum} does not apply
  because the requirement of compactness for $A_{n}-I$ and $B_{n}$ is
  not met.
\end{example}

\begin{example}
  \label{ex:other-diagonal}
  Another case, where the hypothesis of
  Proposition~\ref{prop:Jspectrum} is not met is the following. In view
  of \cite[Thm.\,3]{guseinov1976},  see also \cite[Thm. 10.12]{teschl} one can find sequences
  $\{a_{n}(k)\}_{n\in\mathbb{Z}}$ and $\{b_{n}(k)\}_{n\in\mathbb{Z}}$
  such that each set $\sigma(J_{k})\setminus[-2,2]$ (see the notation
  of Example~\ref{ex:diagonal-blocks}) contains only a finite number
  of simple eigenvalues, while $\sigma_{ac}(J_{k})=[-2,2]$. Hence, the
  spectrum of $\mathcal J$ outside $[-2,2]$ is the countable union of
  finite sets (each point of this union is a simple eigenvalue).
  Clearly, the sequences $\{a_{n}(k)\}_{n\in\mathbb{Z}}$ and
  $\{b_{n}(k)\}_{n\in\mathbb{Z}}$ can be chosen in such a way that
  each eigenvalue of $\mathcal J$ has arbitrary multiplicity. Also, on
  the basis of \cite[Thm.\,3]{guseinov1976} and
  \cite[Sec.\,10.4]{teschl}, one establishes that there are operators
  $J_{k}$, $k\in\mathbb{Z}$ such that the eigenvalues of $\mathcal J$
  have accumulation points. Moreover, examples of operators
  $\mathcal J$ with eigenvalues embedded in the absolutely continuous
  spectrum can be obtained using
  \eqref{eq:orthogonal-decomposition-matrix}. Indeed, by
  \cite{MR1190777,MR1173094} there are scalar discrete Schr\"odinger
  operators $J_{k}$ having embedded simple eigenvalues in $(-2,2)$.
  The embedded eigenvalues of the so constructed operator $\mathcal J$
  may be of any multiplicity. We emphasize that in this particular
  construction given by \eqref{eq:orthogonal-decomposition-matrix}, if
  one imposes the
  conditions of Definition~\ref{assumption:introduction}, then even
  the first moment condition excludes the existence of embedded eigenvalues.
\end{example}

\section {Construction of the Jost solutions for the perturbed
  Laplacian}\label{sec:jost-solutions}

This section is devoted to the proof of Theorem~\ref{thm:jostsolutions} below, stated partly e.g.
in \cite[Theorem 3]{mutlu} and based on arguments of \cite{ABC,
  guseinov1977}. This result is central for all the subsequent arguments in
this paper, and since none of the references provides a complete proof, we
give one here.  Note that in contrast to the previous section,  we do not ask for
$B_n$ and $A_n-I$ to be compact.  In fact such a condition is only used in the analysis
of the spectrum and thus is absent from the next three sections. \medskip

Note that all notions of continuity and analyticity of operator valued functions are always understood with respect to the operator norm.
Moreover, symbols like $o(1)$ or $o(n^{-k})$ refer to sequences of operator valued sequences $\{F_n(z)\}_{n\in\mathbb Z}$ such that 
$\|F_n(z)\| \to 0$ and 
$n^k\|F_n(z)\| \to 0$, respectively, as $n\to\infty$.

\begin{theorem}\label{thm:jostsolutions}
  Suppose that the $k$-th moment condition in Definition
   \ref{assumption:introduction} holds with some integer $k \geq 1$ and that $A_{n}^{-1}\in\mathscr{B}(H)$ for
   all $n\in\mathbb{Z}$. If $z\in\mathbb C$ with $0<|z|\leq 1$,
   then there exist formal solutions $U^{\pm}(z)$ of
   \eqref{eq:eigenvalue-operator} with $\lambda=z+z^{-1}$, called the
   \define{Jost solutions} for the operator $\mathcal J$, given by the
   formulae
   \begin{align}\label{eq:jostuplus}
      U^+_n(z) &= T_nz^n \Big(I + \sum_{m=1}^{\infty}K_{n,m}z^m \Big), \\
      \label{eq:jostuminus}
      U^-_n(z) &= R_nz^{-n}\Big(I + \sum_{m=1}^{\infty}M_{n,m}z^{m}\Big),
   \end{align}
   where the individual terms are defined explicitly by recursive formulae,
   outsourced to the Remark \ref{explicitU}. The series \eqref{eq:jostuplus} and
   \eqref{eq:jostuminus} are convergent and satisfy\footnote{As an analogy, note
   that if $z=e^{i\theta}$, then $U^+(e^{i\theta})\sim e^{i n\theta}$ as
   $n\to\infty$ and $U^-(e^{i\theta})\sim e^{-i n\theta}$ as $n\to -\infty$.
   There is a clear connection to incoming and outgoing plane waves.}
   \begin{equation}
      \label{eq:jostplusasymp}
      \begin{split}
	 U^+_n(z) &= z^n(I+o(n^{-k+1})),\ n\to\infty, \\
	 U^-_n(z) &= z^{-n}(I+o(n^{-k+1})),\ n\to-\infty.
      \end{split}
   \end{equation}
   Moreover, $U^{\pm}(z)$ are holomorphic in $z$ in the punctured
   disk $\{0<|z|<1\}$ and continuous up to the boundary $|z|=1$.
\end{theorem}

\begin{remark}\label{explicitU}
The terms in \eqref{eq:jostuplus} and \eqref{eq:jostuminus} are given explicitly by
\begin{equation}\label{eq:tnkn}
 T_n:=\prod_{p=n}^{\infty}A_p^{-1} = A_{n}^{-1} A_{n+1}^{-1} \cdots , \quad
 R_n:= \, {\vphantom{\prod}}'\!\prod_{p=-\infty}^{n-1}A_p^{-1} = \cdots A_{n-2}^{-1} A_{n-1}^{-1}.
\end{equation}
Moreover, $K_{n,m},$ satisfy the following recursive relations (cf. 
for example the recursive formulae in \cite[Theorem 2.1]{ABC})
\begin{equation}\label{eq:knrecur}
\begin{split}
&K_{n,1}:=-\sum_{p=n+1}^{\infty}T_p^{-1}B_pT_p, \\
  &K_{n,2}:=-\sum_{p=n+1}^{\infty}T_p^{-1}B_pT_pK_{p,1} +
  \sum_{p=n+1}^{\infty}T_p^{-1}(I-A_p^2)T_p, \\
  &K_{n,m+2}:=-\sum_{p=n+1}^{\infty}T_p^{-1}B_pT_pK_{p,m+1} +
  \sum_{p=n+1}^{\infty}T_p^{-1}(I-A_p^2)T_pK_{p+1,m} + K_{n+1,m}.
\end{split}
\end{equation}
Similarly,  $M_{n,m},$ satisfy the following recursive relations
\begin{equation}\label{eq:knrecur2}
\begin{split}
  &M_{n,1}:=-\sum_{p=-\infty}^{n-1}R_p^{-1}B_pR_p, \\
  &M_{n,2}:=-\sum_{p=-\infty}^{n-1}R_p^{-1}B_pR_pM_{p,1} +
  \sum_{p=-\infty}^{n-1}R_p^{-1}(I-A_{p-1}^2)R_p, \\
  &M_{n,m+2}:=-\sum_{p=-\infty}^{n-1}R_p^{-1}B_pR_pM_{p,m+1} +
  \sum_{p=-\infty}^{n-1}R_p^{-1}(I-A_{p-1}^2)R_pM_{p+1,m} + M_{n+1,m}.
\end{split}
\end{equation}
\noindent The convergence of the infinite products in \eqref{eq:tnkn} is shown in Proposition \ref{prop:help}.
The convergence of the series in \eqref{eq:knrecur} and \eqref{eq:knrecur2}, and estimates for their norms are shown in Proposition~\ref{prop:weaklemma}.
\end{remark}

The rest of this section is devoted to the proof of
Theorem~\ref{thm:jostsolutions}. First observe that the invertibility of
$A_{n}$ for all $n\in\mathbb{Z}$ together with the 
\eqref{eq:tnkn} and \eqref{eq:knrecur}, \eqref{eq:knrecur2} guarantee that \eqref{eq:jostuplus} is
a formal solution of \eqref{eq:eigenvalue}. Indeed, this can be seen by
plugging \eqref{eq:jostuplus} into \eqref{eq:eigenvalue} and equating the
coefficients of each power of $z$. Equating the coefficients of the $z^{n-1}$
terms yields $A_{n-1}T_{n-1}=T_n$, which is true given the definition of $T_n$ in \eqref{eq:tnkn}. For
the $z^n$ terms we obtain
\[T_nK_{n-1,1}+B_nT_n = T_nK_{n,1},\]
which follows from the $m=1$ case of
\eqref{eq:knrecur}. Continuing this analysis we see that the equality of
the $z^{n+j}$ coefficients, for any $j\geq 0$, follows from the $m=j+1$
case of \eqref{eq:knrecur}. Thus the proof of Theorem~\ref{thm:jostsolutions}
centers on proving convergence of the infinite products in \eqref{eq:tnkn}
and the series in \eqref{eq:knrecur}
and their asymptotic behaviour.
\medskip

We write out the arguments for $U^+(z)$, breaking the proof into several
steps. In Proposition~\ref{prop:help} we will show that the operators $T_n$
from \eqref{eq:tnkn} are well-defined and boundedly invertible, hence the
operators $K_{n,m}$ are also well defined. We give bounds for them in
Lemma~\ref{lem:josthelp} and Proposition~\ref{prop:weaklemma} which show that
the series expressions \eqref{eq:jostuplus} and \eqref{eq:jostuminus} for the
Jost solutions converge in norm. The proof for $U^-(z)$ is essentially
identical.  

\begin{proposition}\label{prop:help}
  Assume the $k$-th moment condition in Definition
  \ref{assumption:introduction} with some $k\geq 0$
  and that $A_n^{-1}\in\mathscr B(H)$ for all $n\in\mathbb Z$.
 Then the operators $T_n$ and $R_n$ from \eqref{eq:tnkn} are well-defined, belong to 
  $\mathscr{B}(H)$ and are boundedly invertible, that is,
  $T_{n}^{-1}, R_{n}^{-1}\in\mathscr{B}(H)$. 
  The norms of the operators and their inverses are bounded uniformly in $n \in \mathbb{Z}$ and moreover
  \begin{equation}\label{TR-asymptotics}
      \begin{split}
	 T_n(z) &= I+o(n^{-k}),\ n\to\infty, \\
	 R_n(z) &= I+o(n^{-k}),\ n\to-\infty.
      \end{split}
   \end{equation}
\end{proposition}
\begin{proof}
  Note that by any moment condition, the operators $A_n$ are 
  uniformly bounded and tend to $I$ in the operator norm.
  Therefore their inverses are uniformly bounded too. Indeed, 
  for $|n|$ large enough such that $\|I - A_n\| < 1$ we see that 
  \begin{align*}
  &A_n^{-1} =  (I - (I-A_n))^{-1} = \sum_{j=0}^\infty (I - A_n)^j,  \\
  &\Rightarrow \quad \|A_n^{-1}\| \le  \sum_{j=0}^\infty \|I - A_n\|^j = \frac{1}{1 - \| I-A_n \|}. 
  \end{align*}
  Consequently,  already by the weakest ($0$-th order) moment condition the series
   $$
   \sum_{n=-\infty}^{\infty} \|I - A_n\| < \infty, \quad
   \sum_{n=-\infty}^{\infty} \|I - A_n^{-1}\| \le
   \sum_{n=-\infty}^{\infty} \|A_n^{-1}\|\, \|I - A_n\| < \infty
   $$
   converge, hence also the infinite products
   \begin{align*}
   &T_n^{-1} = \, {\vphantom{\prod}}' \prod_{p=n}^\infty A_p := \cdots A_{p+1} A_p
   := \lim_{m\to \infty} \, {\vphantom{\prod}}'\prod_{p=n}^m A_p,  \\
   &T_n = \prod_{p=n}^\infty A_p^{-1} = A_p^{-1} A_{p+1}^{-1}\cdots
   := \lim_{m\to \infty} \prod_{p=n}^m A_p^{-1}, 
   \end{align*}
   converge by \cite[Theorem 2.3]{welstead} and are uniformly bounded.  
   A similar argument works for $R_n$ and its inverse.  For the proof of \eqref{TR-asymptotics}, note
   by a telescope argument
   \begin{equation}\label{telescope}\begin{split}
   D_{n,m} := \prod_{p=n}^m A_p^{-1},  \quad T_n &= \sum_{m=n}^\infty (D_{n,m+1}-D_{n,m}) + A_n^{-1}
   \\ &= \sum_{m=n}^\infty D_{n,m} (A^{-1}_{m+1}-I)   + A^{-1}_n.
   \end{split}\end{equation}
   Since the norms of $T_n$ and hence also of $D_{n,m}$ are bounded uniformly in $n$ and $m$, 
   by some constant $C>0$,  we conclude (choosing $C \geq 1$)
   $$
   \| T_n - I \| \leq C \sum_{m=n}^\infty \| A^{-1}_{m}-I \| \leq C \cdot \sup_p \| A^{-1}_p \| \cdot 
   \sum_{m=n}^\infty \| A_{m}-I \|.
   $$
   Now,  given the $k$-th moment condition, we can estimate as $n \to \infty$
   $$
   \| T_n - I \| \leq C \cdot \sup_p \| A^{-1}_p \| \cdot n^{-k} 
   \sum_{m=n}^\infty m^k \| A_{m}-I \| = o(n^{-k}).
   $$
   This proves \eqref{TR-asymptotics} for $T_n$.  A similar argument works for $R_n$.
\end{proof}

Assume from now on that some $k$-th moment condition in Definition
  \ref{assumption:introduction} with $k\geq 1$ holds.  Then,  by the uniformity of norm bounds
  as explained in the proof of Proposition~\ref{prop:help},  there exists a constant
 $\mathscr{C}>1$ such that for all $j$,
 \begin{align}\label{C-constant}
 \|T_j\| + \|T_j^{-1}\| + \|I+A_j\| +
    \sum_{n=-\infty}^{\infty}(|n|+1)(\|I-A_n\|+\|B_n\|) \leq \mathscr{C}.
    \end{align}
For showing that the formal expression \eqref{eq:jostuplus} is convergent and hence a true solution of \eqref{eq:eigenvalue-operator}, it is essential to understand how the constants $C_{n,m}$ introduced in Proposition~\ref{prop:weaklemma} below depend on $n,m \in \mathbb Z$.

\begin{proposition}\label{prop:weaklemma}
We have the following estimate for each $m,n \in \mathbb Z$
\begin{equation}\label{eq:knmestimateweak}
\|K_{n,m}\|\leq C_{n,m}\sum_{p=n+\lfloor\frac m2\rfloor}^{\infty}(\|I-A_p\| + \|B_p\|),
\end{equation}
where the constants $C_{n,m}$ may be defined recursively by
$C_{n,1}=\mathscr{C}^2, C_{n,2}=\mathscr{C}^5$ 
and
\[
C_{n,m}= \mathscr{C}^5\prod_{q=n+1}^{n+m-2}\Big(1+\mathscr{C}^3
 \sum_{p=q}^{\infty} (\|1-A_p\|+\|B_p\|)\Big)
 \quad \textrm{ for }m>2
 \]
 where $\mathscr{C}$ is the constant from \eqref{C-constant}.
 Note that the sequence $\{C_{n,m}\}$ is non-increasing in $n$ for each fixed
 $m$, and non-decreasing in $m$ for each fixed $n$.
\end{proposition}
\begin{proof} We argue by induction. For $m=1$,
  \[\|K_{n,1}\|=\|\sum_{p=n+1}^{\infty}T_p^{-1}B_pT_p\|
    \leq\mathscr{C}^2\sum_{p=n}^{\infty}\|B_p\|\]
    as desired. 
    For $m=2$, by the $m=1$ case,
    \begin{align*}
       \|K_{n,2}\|
       &=\Big\|-\sum_{p=n+1}^{\infty}T_p^{-1}B_pT_pK_{p,1}
       + T_p^{-1} (1-A_p^2) T_p \Big\|
       \\
       &\leq
       \sum_{p=n+1}^{\infty} \|T_p^{-1} \|\, \|B_p \|\, \|T_p \|\, \|K_{p,1} \|
       + \|T_p^{-1}\|\, \|1-A_p\|\,\|1+A_p\|\, \|T_p\|
       \\
       &\leq
       \sum_{p=n+1}^{\infty} \mathscr C^5 \|B_p \| + \mathscr C^3 \|1-A_p\|
       \leq
       \mathscr C^5 \sum_{p=n+1}^{\infty}  \|B_p \| + \|1-A_p\|,
    \end{align*}
    again as claimed. Now assume $m>2$. Then
    \[
    \|K_{n,m}\|\leq \|K_{n+1,m-2}\| + \mathscr{C}^2
    \sum_{p=n+1}^{\infty}\|B_p\|\|K_{p,m-1}\| +
    \mathscr{C}^3\sum_{p=n+1}^{\infty}\|I-A_p\|\, \|K_{p+1,m-2}\|.
    \]
By the inductive hypothesis and the monotonicity of $C_{n,m}$, we obtain
\begin{align*}
\|K_{n,m}\| &\leq C_{n+1,m-2}
  \sum_{p=n+\lfloor\frac m2\rfloor}^{\infty}(\|I-A_p\|+\|B_p\|)\\
  &\phantom{\le\ } + \mathscr{C}^2C_{n+1,m-1}\sum_{p=n+1}^{\infty}
  \sum_{q=p+\lfloor\frac {m-1}2\rfloor}^{\infty}\|B_p\|(\|I-A_q\|+\|B_q\|)\\
  &\phantom{\le\ } + \mathscr{C}^3C_{n+1,m-2}\sum_{p=n+1}^{\infty}
  \sum_{q=p+1+\lfloor\frac {m-2}2\rfloor}\|I-A_p\|(\|I-A_q\|+\|B_q\|).
\end{align*}
Using again the monotonicity of $C_{n,m}$, we arrive at the estimate
\begin{align*}
\|K_{n,m}\| &\leq C_{n+1,m-1}
  \sum_{p=n+\lfloor\frac m2\rfloor}^{\infty}(\|I-A_p\|+\|B_p\|)\\
  &\phantom{\le\ }+ \mathscr{C}^3C_{n+1,m-1}\sum_{p=n+1}^{\infty}
  \sum_{q=p-1+\lfloor\frac m2\rfloor}^{\infty}(\|I-A_p\|+\|B_p\|)(\|I-A_q\|+\|B_q\|).
\end{align*}
The second term here is bounded by
\[\mathscr{C}^3C_{n+1,m-1}\sum_{p=n+1}^{\infty}(\|I-A_p\|+\|B_p\|)
  \sum_{q=n+\lfloor\frac m2\rfloor}^{\infty}(\|I-A_q\|+\|B_q\|).\]
Thus we have the desired bound with
\[C_{n,m}=C_{n+1,m-1}(1 + \mathscr{C}^3\sum_{p=n+1}^{\infty}(\|I-A_p\|+\|B_p\|)).\]
Using the inductive hypothesis to identify $C_{n+1,m-1}$, we are done.
\end{proof}

We can estimate the constants $C_{n,m}$ even further to avoid dependence on
$m$.

\begin{corollary}\label{lem:josthelp}
We have the following estimate for each $m,n \in \mathbb Z$
\begin{equation}\label{eq:knmestimate}
\begin{split}
&\|K_{n,m}\|\leq C_n\sum_{p=n+\lfloor\frac m2\rfloor}^{\infty}(\|I-A_p\| + \|B_p\|), \\
&\textup{where }\quad C_n:=\mathscr{C}^5\exp\left[\mathscr{C}^3
  \sum_{p=n+1}^{\infty}(p-n)(\|I-A_p\|+\|B_p\|)\right].
\end{split}
\end{equation}
Note that the sequence $\{C_n\}$ is decreasing in $n$.
\end{corollary}

\begin{proof}
The proof follows from Proposition \ref{prop:weaklemma}, since
$$
C_{n,m} \leq \mathscr{C}^5\prod_{q=n+1}^{\infty} \Big(1+\mathscr{C}^3
\sum_{p=q}^{\infty}(\|1-A_p\|+\|B_p\|) \Big).
$$
The infinite product converges by the following reasoning: using $\log(1+x)<x$
\[\prod_{q=n+1}^{\infty} \Big( 1+\mathscr{C}^3\sum_{p=q}^{\infty}(\|1-A_p\|+\|B_p\|) \Big) 
\leq \exp\Big[\mathscr{C}^3\sum_{q=n+1}^{\infty}\sum_{p=q}^{\infty}(\|I-A_p\|+\|B_p\|)\Big].\]
Rearranging the double series, which may be done since all terms are non-negative, gives a bound
\[C_{n,m}\leq \mathscr{C}^5 \exp \Big[ \mathscr{C}^3
  \sum_{p=n+1}^{\infty}(p-n)(\|I-A_p\|+\|B_p\|) \Big],\] which,  by
\eqref{eq:starassumption} with $k=1$,  is finite for each fixed $n$. This completes the
proof.
\end{proof}

We now finish the proof of Theorem \ref{thm:jostsolutions}. By Corollary
\ref{lem:josthelp},
\[\sum_{m=1}^{\infty}\|K_{n,m}\|\leq C_n\sum_{m=1}^{\infty}
  \sum_{p=n+\lfloor\frac m2\rfloor}^{\infty}(\|I-A_p\| + \|B_p\|).\]
Rearranging the double series and using the fact that the sequence $\{C_n\}$ is
decreasing in $n$, we obtain
\begin{equation}\label{tail}
\begin{split}
\sum_{m=1}^{\infty}\|K_{n,m}\| &\leq
C_n\sum_{q=n}^{\infty}(2(q-n)+1)(\|I-A_q\| + \|B_q\|) \\
&\leq C \sum_{q=n}^{\infty} |q| \cdot (\|I-A_q\| + \|B_q\|),
\end{split}
\end{equation}
for some uniform constant $C>0$. By the first moment condition, the right-hand
side is finite. This is enough to show that the series \eqref{eq:jostuplus}
converges and that one has the required analyticity and continuity
properties. It also proves \eqref{eq:jostplusasymp} for $k=1$ since the right hand side
in \eqref{tail} is the tail of a convergent series. \medskip

Assuming the $k$-th moment condition in Definition
\ref{assumption:introduction} with $k\geq 1$, the asymptotics \eqref{eq:jostplusasymp}
follows from \eqref{TR-asymptotics} once we estimate \eqref{tail} further for positive $n$
\begin{equation}
\begin{split}
\sum_{m=1}^{\infty}\|K_{n,m}\| &\leq C \sum_{q=n}^{\infty} |q| \cdot (\|I-A_q\| + \|B_q\|) \\
&\leq C n^{-k+1} \sum_{q=n}^{\infty} q^k \cdot (\|I-A_q\| + \|B_q\|) = o(n^{-k+1}), \ n \to \infty.
\end{split}
\end{equation}
This completes the proof of
Theorem \ref{thm:jostsolutions} for $U^+(z)$. The argument for $U^-(z)$ is
essentially identical.%
\medskip

In the case that the perturbation
$\mathcal S^*(\mathcal A - \mathcal I) + \mathcal B + (\mathcal A - \mathcal I) \mathcal S$
is exponentially decreasing in $n$ in the sense of Definition~\ref{def:block-jacobi-operator}
(which includes the case where it has compact support), it is possible to extend the domain
of definition of the Jost solutions as functions of $z$ outside the unit disk.

\begin{proposition}\label{Jost-infinite-moment-prop}
Let $A_{n}$ be invertible with $A_{n}^{-1}\in\mathscr{B}(H)$ for
  any $n\in\mathbb{Z}$. If the exponential moment condition in Definition
  \ref{assumption:introduction} holds, then Theorem \ref{thm:jostsolutions}
  in fact holds\footnote{Note that $e^{\varepsilon / 2}> 1$ and hence we
    indeed have a larger domain of convergence than under the first or second
    moment condition.} on the punctured disk $0<|z|<e^{\varepsilon / 2}$.
    In particular, \eqref{eq:jostplusasymp}
    holds for any $k\in\mathbb N_0$.
\end{proposition}
\begin{proof} Fix $n\in \mathbb Z$. By Corollary \ref{lem:josthelp},
we have using $n+\lfloor\frac m2\rfloor>\frac m2 - 1$ the estimate
\begin{align*}
  \|K_{n,m}\|\leq C_n e^{-\varepsilon \left( n+\lfloor\frac m2\rfloor \right)}
  \sum_{p=n+\lfloor\frac m2\rfloor}^{\infty} e^{\varepsilon |p|} (\|I-A_p\| + \|B_p\|)
  \leq C e^{-\varepsilon \left( n+ \frac m2 \right)}.
\end{align*}
Thus, if $|z|<\ e^{\varepsilon / 2}$, the series
$\sum_{m=1}^{\infty}K_{n,m}z^m$ is convergent. The remainder of the proof
proceeds exactly as in the proof of Theorem \ref{thm:jostsolutions}.
\end{proof}

\section {Wronskian and the Jost solutions as a fundamental system}

\subsection {Wronskian of formal operator solutions}\label{sec:wronskian-eigenvalues}

We shall introduce the notion of a \define{Wronskian}, which is central for all
the arguments and computations below.
\begin{definition}
   \label{def:Wronskian}\label{def:Wronskianvector}
   \begin{enumerate}
   \item For sequences $U = \{U_n\}_{n\in\mathbb Z}$ and
     $V = \{V_n\}_{n\in\mathbb Z}$ of bounded operators on $H$, their
     \emph{Wronskian} is defined by
   \begin{align*}
      W_n(U,V) := U_{n-1} A_{n-1} V_n  - U_{n} A_{n-1} V_{n-1}
      \in \mathscr{B}(H).
   \end{align*}
 \item For a sequence $U = \{U_n\}_{n\in\mathbb Z}$ of bounded
   operators on $H$ and $v = \{v_n\}_{n\in\mathbb Z}$ of vectors in
   $H$, we define their \emph{Wronskian} by
   \begin{align*}
      W_n(U,v) := U_{n-1} A_{n-1} v_n  - U_{n} A_{n-1} v_{n-1}
      \in H.
   \end{align*}
   \end{enumerate}
\end{definition}

\begin{proposition}
   \label{prop:Green}
   For all $n\in\mathbb Z$, one has that
   \begin{align}
      \label{eq:Green0}
      &W_{n+1}(U,V) -  W_{n}(U,V) =  U_n (\tau V)_n - (\tau U^*)_n^* V_n, \\
      \label{eq:Green}
      &W_{m+1}(U,V) -  W_{n}(U,V) =  \sum_{j=n}^{m} U_j (\tau V)_j - (\tau U^*)_j^* V_j .
   \end{align}
\end{proposition}
\begin{proof}
  Clearly, \eqref{eq:Green} follows from \eqref{eq:Green0} by summing
  both sides from $n$ to $m$, and \eqref{eq:Green0} follows from a
  straightforward calculation:
  \begin{multline*}
     U_n (\tau V)_n - (\tau U^*)_n^* V_n
     \\ 
     \begin{aligned}[b]
	&= U_n [ A_{n-1} V_{n-1} + B_{n}V_{n} + A_{n} V_{n+1} ]
	- [ U_{n-1} A_{n-1} + U_{n}B_{n} + U_{n+1} A_{n}] V_n
	\\ &=  U_n A_{n-1} V_{n-1} + U_n A_{n} V_{n+1}
	- U_{n-1} A_{n-1}V_n - U_{n+1} A_{n}V_n
	\\ &= W_{n+1}(U,V) -  W_{n}(U,V).
     \end{aligned}
     \qedhere
  \end{multline*}
\end{proof}

We now observe that the Wronskian of formal operator
solutions to \eqref{eq:eigenvalue-operator} is in fact independent of
$n \in \mathbb Z$. To be precise, let
$U(z) = \{U_n(z)\}_{n\in\mathbb Z}$ be such a formal operator solution
to \eqref{eq:eigenvalue-operator} with $\lambda = z + z^{-1}$. This
equation amounts for each $n \in \mathbb Z$ to
\begin{align}
   \label{eq:recursion}
   A_{n-1} U_{n-1}(z) + (B_{n} - \lambda) U_{n}(z) + A_{n} U_{n+1}(z) = 0 .
\end{align}
Taking the adjoint on both sides, we obtain
\begin{align*}
  [ U_{n-1}(z) ]^* A_{n-1} + [ U_{n}(z) ]^* ( B_{n} -
  \overline\lambda)
  + [ U_{n+1}(z) ]^* A_{n} = 0
\end{align*}
which can be written as
$[U(z)]^* (\tau - \overline{\lambda} ) = 0$, if we understand
that $\tau$ acts from the right on
$[U(\lambda)]^* := \{ [U_n(z)]^*\}_{n\in\mathbb Z}$. Consequently we
obtain
\begin{equation*}
  (\tau - \lambda )U(z) = 0\,\Longleftrightarrow[U(\overline z)]^* (\tau - \lambda ) = 0\,.
\end{equation*}

\begin{corollary}\label{cor:wronskianofsols}
  Let the sequences $U(z) = \{U_n(z)\}_{n\in\mathbb Z}$ and
  $V(z) = \{V_n(z)\}_{n\in\mathbb Z}$ be formal operator solutions to
  \eqref{eq:eigenvalue-operator}
  and let $v(z) = \{v_n(z)\}_{n\in\mathbb Z}$ be a formal vector solution to
  \eqref{eq:eigenvalue}.
  Then for $\lambda = z+z^{-1}$ and
  $\lambda' = z'+z'^{-1}$ the following is true.
   \begin{enumerate}
      \item The operator version of the Christoffel-Darboux formula holds
 (\emph{cf}. \cite{akhiezer})
      \begin{align*}
	 W_{n+1}( U(\overline z)^*, V(z')) -  W_{n}(U(\overline z)^*,V(z'))
	 = (\lambda' - \lambda) U_n(\overline z)^* V_n(z').
      \end{align*}
      Similarly, for the vector solution we obtain
      \begin{align*}
	 W_{n+1}( U(\overline z)^*, v(z')) -  W_{n}(U(\overline z)^*,v(z'))
	 = (\lambda' - \lambda) U_n(\overline z)^* v_n(z').
      \end{align*}

      \item The operators
      $W_{n}(U(\overline z)^*,V(z))$ 
      and the vectors $W_{n}(U(\overline z)^*,v(z))$ 
      do not depend on $n \in \mathbb Z$.
   \end{enumerate}
  \end{corollary}

\begin{proof}
   \begin{enumerate}
   \item A straightforward calculation shows for every $n\in\mathbb Z$
      \begin{align*}
	 W_{n+1}( U(\overline z)^*, V(z')) -  W_{n}(U(\overline z)^*,V(z'))
	 &= U_n(\overline z)^* (\tau V(z'))_n - (\tau U(\overline z))_n^* V_n(z')
	 \\ &= \lambda' U_n(\overline z)^* V_n(z')
	 - ( \overline\lambda U_n(\overline z) )^* V_n(z')
	 \\ &= (\lambda' - \lambda) U_n(\overline z)^* V_n(z').
      \end{align*}

      \item The second statement follows if we take $\lambda = \lambda'$.
      \qedhere
   \end{enumerate}
\end{proof}

\begin{notation*}
   If the Wronskian does not depend on $n$,  we shall omit the index $n$.
\end{notation*}

\begin{remark}
  Let $A_{n}$ be invertible with $A_{n}^{-1}\in\mathscr{B}(H)$ for any
  $n\in\mathbb{Z}$. A formal solution $\{U_{n}(z)\}_{n\in\mathbb{Z}}$
  to \eqref{eq:eigenvalue-operator} with $\lambda = z + z^{-1}$ can be
  recursively calculated from \eqref{eq:recursion} if $U_{-1}(z)$ and
  $U_0(z)$ are given (constants with respect to $z$) and
  $U_{0}(z)\neq 0$. Clearly, they are polynomials in
  $\lambda$ with operator coefficients. More precisely,
  if $n\ge 1$, then
   $$
   U_n(z) = \sum_{j=0}^n \lambda^j U_{n,j}, \quad
   [ U_n(\overline z) ]^* = \sum_{j=0}^n \lambda^j U_{n,j}^*,
   $$
   with operators $U_{n,j}\in \mathscr{B}(H)$.
   An analogous statement holds for $n\le -2$ with
   $$
   U_n(z) = \sum_{j=0}^{n+1} \lambda^{-j} U_{n,j}, \quad
   [ U_n(\overline z) ]^* = \sum_{j=0}^{n+1} \lambda^{-j} U_{n,j}^*.
   $$
   While $\{ U_n(\overline z)^*\}_{n\in\mathbb N}$ is a solution of the adjoint equation $V(\tau - \lambda) = 0$, it is in general not a solution of the original equation $(\tau - \lambda) V = 0$.
\end{remark}

\subsection {Fundamental systems of solutions}

\begin{notation*}
If $V = \{V_n\}_{n\in\mathbb Z}$ is a sequence of bounded operators in $H$ and if 
$x\in H$ and $\alpha$ is an operator on $H$, then we set
\begin{equation*}
   V\alpha := \{V_n\alpha\}_{n\in\mathbb Z},
   \qquad
   Vx := \{V_n x\}_{n\in\mathbb Z}.
\end{equation*}
\end{notation*}

\begin{proposition}
   \label{prop:FS:coeff}
   Let $R(z)$ and $\widehat R(z)$ be formal operator solutions of
   \eqref{eq:eigenvalue-operator} for $\lambda = z + z^{-1}$. Assume
   that
   \begin{equation}\label{eq:wronsikanzero}
   \begin{split}
   &W( R(\overline z)^*, \widehat R(z) )^{-1}, \ W( \widehat R(\overline z)^*,  R(z) )^{-1} \in\mathscr{B}(H), \\
      &W( R(\bar z)^*,  R(z) ) = W( \widehat R(\bar z)^*,  \widehat R(z) ) = 0.
      \end{split}
   \end{equation}
   If $V(z) = \{V_n(z)\}_{n\in\mathbb Z}$ is a formal operator
   solution of \eqref{eq:eigenvalue-operator} which can be written as
      \begin{equation}
	 \label{eq:FS:operatorrep}
	 V(z) = R(z) \alpha(z) + \widehat R(z) \beta(z),
      \end{equation}
      then the coefficients $\alpha(z),\, \beta(z)$ belong to $\mathscr{B}(H)$, 
      and are given by
      \begin{align}
	 \label{eq:FS:operatorcoeff}
	 \begin{bmatrix} \alpha(z) \\ \beta(z)
	 \end{bmatrix}
	 & = Z_j(z)
	 \begin{bmatrix}
	    V_j(z) \\ V_{j+1}(z)
	 \end{bmatrix}
	 := 
	 \begin{bmatrix}
	    W( \widehat R(\overline z)^*,  R(z) )^{-1} W(\widehat R(\overline z)^*, V(z))\\ 
	    W( R(\overline z)^*,  \widehat R(z) )^{-1} W(R(\bar z)^*,V(z))
	 \end{bmatrix}
      \end{align}
      where the operator $Z_j(z)$ acts on $\mathscr B(H)\oplus \mathscr B(H)$ as
      \begin{equation}
	 \label{eq:Z:FS}
	 \begin{split}
	    Z_j(z)  &=  
	    \begin{bmatrix}
	       - W( \widehat R(\overline z)^*,  R(z) )^{-1} & 0 \\
	       0 & W( R(\overline z)^*,  \widehat R(z) )^{-1}
	    \end{bmatrix}
	    \begin{bmatrix}
	       \widehat R_{j+1}(\overline z)^* A_j 
	       & -\widehat R_j(\overline z)^*  A_j\\
	       - R_{j+1}(\overline z)^* A_j
	       & R_j(\overline z)^* A_j \vphantom{\widehat R}
	    \end{bmatrix} 
	    \\[1ex] 
	    &=
	    \begin{bmatrix}
	       - W( \widehat R(\overline z)^*,  R(z) )^{-1} \widehat R_{j+1}(\overline z)^* A_j 
	       &  W( \widehat R(\overline z)^*,  R(z) )^{-1}  \widehat R_j(\overline z)^*  A_j\\
	       -W( R(\overline z)^*,  \widehat R(z) )^{-1} R_{j+1}(\overline z)^* A_j
	       & W( R(\overline z)^*,  \widehat R(z) )^{-1} R_j(\overline z)^* A_j
	    \end{bmatrix}.
	 \end{split}
      \end{equation}      
      If $v(z) = \{v_n(z)\}_{n\in\mathbb Z}$ is a formal vector
      solution of \eqref{eq:eigenvalue} which can be written as
      \begin{equation}
        \label{eq:fsvectorrep}
	 v(z) = R(z) \alpha(z) + \widehat R(z) \beta(z)\,,
      \end{equation}
      for vectors $\alpha(z), \beta(z) \in H$, then
      \eqref{eq:FS:operatorcoeff} holds with $V(z)$ replaced by $v(z)$ and the Wronskians
      $W(\widehat R(\overline z)^*, v(z))$ and
      $W(R(\bar z)^*,v(z))$ in \eqref{eq:FS:operatorcoeff}
      understood in the sense of Definition~\ref{def:Wronskianvector} (2).
      
\end{proposition}

\begin{proof}
  Assume that operator coefficients $\alpha(z)$ and $\beta(z)$ as in
  \eqref{eq:FS:operatorrep} exist. Recall that by the recursion
  \eqref{eq:recursion}, the formal solution $V(z)$ is uniquely
  determined if we know $V_j(z)$ and $V_{j+1}(z)$ for any
  $j\in\mathbb Z$. We obtain from \eqref{eq:FS:operatorrep} for any
  $j$ that
   \begin{equation}\label{eq:UpUmmatrixeq}
      \begin{bmatrix} V_j(z)\\ V_{j+1}(z)
      \end{bmatrix} 
      = Y_j(z)\begin{bmatrix} \alpha(z) \\ \beta(z)
      \end{bmatrix},\quad\textrm{ where } 
      Y_j(z)=\begin{bmatrix} 
	 R_j(z) & \widehat R_j(z) \\ R_{j+1}(z) & \widehat R_{j+1}(z)
      \end{bmatrix},
   \end{equation}
   so we only need to show that $Z_j(z)$ is a left inverse of $Y_j(z)$.
   This follows from the straightforward calculation
   \begin{align*}
      Z_j(z)Y_j(z) &=  
      \diag \left(-W( \widehat R(\overline z)^*,  R(z) )^{-1},\ 
      W( R(\overline z)^*,  \widehat R(z) )^{-1}\right)
      \\
      &\phantom{=\ \ }\circ 
      \begin{bmatrix}
	 \widehat R_{j+1}(\bar z)^* A_j & -\widehat R_j(\bar z)^* A_j\\ 
	 - R_{j+1}(\bar z)^* A_j & R_j(\bar z)^* A_j
	 \vphantom{\widehat R}
      \end{bmatrix}
      \begin{bmatrix} 
	 R_j(z) & \widehat R_j(z) \\ 
	 R_{j+1}(z) & \widehat R_{j+1}(z)
      \end{bmatrix}
      \\
      &=  
      \diag \left(-W( \widehat R(\overline z)^*,  R(z) )^{-1},\ 
      W( R(\overline z)^*,  \widehat R(z) )^{-1}\right)
      \\
      &\phantom{=\ \ }\circ 
      \begin{bmatrix}
	 -W( \widehat R(\bar z)^*,  R(z) )
	 &
	 -W( \widehat R(\bar z)^*, \widehat R(z) )
	 \\ 
	 W( R(\bar z)^*,  R(z))
	 & W( R(\bar z)^*, \widehat R(z) )
      \end{bmatrix} 
      \\ & = 
      \begin{bmatrix}
	 I & 0 \\ 0 & I
      \end{bmatrix}
   \end{align*}%
   where we used that the off-diagonal elements are zero by \eqref{eq:wronsikanzero}.
   This concludes the proof of the statement for formal operator solutions. 
   The proof for formal vector solutions is analogous.
\end{proof}

The next corollary provides a condition when a formal operator or
vector solution has in fact a representation as in
\eqref{eq:FS:operatorrep} or \eqref{eq:fsvectorrep}, respectively.

\begin{corollary}
   \label{cor:FS:coeffuniqueness}
   Let $R(z)$ and $\widehat R(z)$ be formal solutions of
   \eqref{eq:eigenvalue-operator} for $\lambda = z + z^{-1}$.
   Suppose  that \eqref{eq:wronsikanzero} holds true. 
   Define $Z_j(z)$ as in
   \eqref{eq:Z:FS}. Let $V(z)$ be any formal operator solution of
   \eqref{eq:eigenvalue-operator} and $v(z)$ any formal vector
   solution of \eqref{eq:eigenvalue}. If $Z_j(z)$ is injective for
   some $j\in\mathbb Z$, then $V(z)$ and $v(z)$ can be represented as in
   \eqref{eq:FS:operatorrep} and \eqref{eq:fsvectorrep}, respectively.
   In particular, if 
   \begin{equation}
      \label{eq:FS:coefffConstInN}
      V_n(z) = R_n(z)\alpha(z) + \widehat R_n(z)\beta(z),
   \end{equation}
   for some $n$ and linear operators $\alpha(z),\, \beta(z)$ on $H$, then
   \eqref{eq:FS:coefffConstInN} is true for all $n\in\mathbb Z$.
   The analogous statement for formal vector solutions $v(z)$ is obtained by replacing $V(z)$ by  $v(z)$.
   
\end{corollary}
\begin{proof}
  First observe that no matter what $\alpha(z)$ and $\beta(z)$ are, by
  right-linearity, the right hand side of \eqref{eq:FS:operatorrep} is
  a formal operator solution of \eqref{eq:eigenvalue-operator}. Define
  $\alpha(z),\, \beta(z)$ by \eqref{eq:FS:operatorcoeff} and set
   \begin{equation}
      \label{eq:widetildeV}
      \widetilde V(z) =  R(z)\alpha(z) + \widehat R(z)\beta(z).
   \end{equation}
   We have to show that $\widetilde V(z) = V(z)$. Let
   $$
   {\mathbf V}_j = \begin{bmatrix} V_j(z) \\ V_{j+1}(z)
   \end{bmatrix} \ \textup{and} \
   \widetilde {\mathbf V}_j = \begin{bmatrix} \widetilde V_j(z) \\ \widetilde V_{j+1}(z)
   \end{bmatrix}.
   $$
   Then, omitting the argument $z$ everywhere, we obtain in the
   notation of \eqref{eq:UpUmmatrixeq}
   \begin{align*}
      Z_j({\mathbf V}_j - \widetilde{\mathbf V}_j) 
      & = Z_j \left( {\mathbf V}_j- Y_j \begin{bmatrix} \alpha \\ \beta
      \end{bmatrix}
      \right)
      = Z_j \left( {\mathbf V}_j- Y_j Z_j {\mathbf V}_j \right)
      = Z_j {\mathbf V}_j - Z_j Y_j Z_j {\mathbf V}_j
      \\
      & = Z_j {\mathbf V}_j- Z_j  {\mathbf V}_j
      = 0.
   \end{align*}
   The injectivity of $Z_j(z)$ now shows that
   ${\mathbf V}_j=\widetilde{\mathbf V}_j$. By the recurrence relation
   \eqref{eq:recursion} it follows that $V(z) = \widetilde V(z)$.
   The proof for formal vector solutions is analogous.
\end{proof}

\begin{definition}
  We call a system of formal operator solutions $R(z)$ and
  $\widehat R(z)$ a \define{fundamental system} if for every formal (operator
  or vector) solution $V(z)$ there exist (operator or vector)
  coefficients $\alpha(z)$ and $\beta(z)$ such that
  $V(z) = R(z)\alpha(z) + \widehat R(z)\beta(z)$.
\end{definition}

Sufficient conditions for formal operator solutions to be a fundamental system are provided in 
Corollary~\ref{cor:FS:coeffuniqueness} above.

\begin{remark}
   \label{rem:FS:finitedim}
   If the Hilbert space $H$ is finite dimensional with $\dim H = d$,
   then $Z_j(z)$ is a $2d\times 2d$ matrix. 
   If additionally the conditions of Proposition~\ref{prop:FS:coeff} are satisfied,  then
   by the proof of Proposition~\ref{prop:FS:coeff},  $Z_j(z)$ 
   admits a right inverse $Y_j(z)$ and hence is surjective.
   Therefore $Z_j(z)$ is bijective and thus 
 the injectivity condition on $Z_j(z)$ in 
   Corollary~\ref{cor:FS:coeffuniqueness} holds automatically, 
   so formal operator solutions $R(z), \widehat R(z)$ satisfying the hypotheses in Proposition~\ref{prop:FS:coeff} are already  a fundamental system.

\end{remark}

\subsection {Jost solutions as a fundamental system}\label{sec:fund-syst-solut}

We assume either the first or the exponential moment condition from Definition
\ref{assumption:introduction} and that $A_{n}^{-1}\in\mathscr{B}(H)$
for all $n\in\mathbb{Z}$.  The assumption of higher (than first) order moment condition does not 
lead to stronger results here.
Let $U^{\pm}(z)$ be the Jost solutions from Section~\ref{sec:jost-solutions}. 
Then also $U^\pm(z^{-1})$ are formal
operator solution of \eqref{eq:eigenvalue-operator}, since the
eigenvalue parameter is $\lambda = z + z^{-1}$.
In our first result below, we note
that $\{R(z), \widehat R(z)\}:=\{U^{\pm}(z), U^{\pm}(z^{-1})\}$
satisfies the conditions of Proposition \ref{prop:FS:coeff} for $|z|=1$
and $z \neq \pm 1$,  which corresponds to $\lambda\in(-2, 2)$.

\begin{proposition}\label{prop:wronsks}
 Let $A_{n}^{-1}\in\mathscr{B}(H)$ for all $n\in\mathbb{Z}$. Assume either
 \begin{enumerate}[label={(\alph*)}]
    \item the first moment condition and consider $|z|=1$,
    \item or the exponential moment condition and consider
    $e^{-\varepsilon / 2} < |z| < e^{\varepsilon / 2}$.
 \end{enumerate}
   Then for $z \neq \pm 1$, the formal operator solutions
   $\{U^{\pm}(z), U^{\pm}(z^{-1})\}$ satisfy the conditions of
   $\{R(z), \widehat R(z)\}$ in Proposition \ref{prop:FS:coeff}. More
   precisely
\begin{align}\label{eq:wronskzzinv}
    &W( U^{\pm}(\bar z)^*, U^{\pm}(z^{-1})) = \pm(z^{-1}-z)I, \\
\label{eq:wronskzzself}
    &W( U^{\pm}(\bar z)^*, U^{\pm}(z)) = 0.
\end{align}
\end{proposition}
\begin{proof}
  By Corollary \ref{cor:wronskianofsols}, both of the Wronskians
  \eqref{eq:wronskzzinv} and \eqref{eq:wronskzzself} are independent
  of $n$. We may therefore find them by computing their limits as
  $n\to\pm\infty$. We begin with the Wronskian \eqref{eq:wronskzzinv}.
  The sets $|z|=1$ and
  $e^{-\varepsilon / 2} < |z| < e^{\varepsilon / 2}$ are preserved
  under taking the inverse. Thus Theorem \ref{thm:jostsolutions}
  applies to both $U^{\pm}(z)$ and $U^{\pm}(z^{-1})$ simultaneously,
  assuming the first or the exponential moment conditions. We
  use the asymptotics \eqref{eq:jostplusasymp} to obtain
  \begin{align*}
     \lim_{n\to\infty}W_n(U^{\pm}(\bar z)^*, U^{\pm}(z^{-1})) 
     &= \, 
     \lim_{n\to\pm\infty}z^{\pm (n-1)}(I + o(1))A_{n-1}(z^{-1})^{\pm n}(I+o(1)) 
     \\
     & \qquad -\, \lim_{n\to\pm\infty} z^{\pm n}(I+o(1))A_{n-1}(z^{-1})^{\pm (n-1)}(I+o(1)).
  \end{align*}
 Simplifying and using the fact that $A_{n}\to I$ as $n\to\infty$, we
 obtain \eqref{eq:wronskzzinv}. 
 We obtain in a similar fashion
   \[\lim_{n\to\pm\infty} W_n(U^{\pm}(\bar z)^*, U^{\pm}(z))
   =\lim_{n\to\pm\infty}z^{\pm(2n-1)} \big(I + o(1)- (I + o(1)) \big) = 0,
   \]
   which allows us to conclude \eqref{eq:wronskzzself}.
\end{proof}

Now we come to the proof that $U^{\pm}(z)$ and $U^{\pm}(z^{-1})$ form
a fundamental system of formal solutions of
\eqref{eq:eigenvalue-operator} for any $z\neq \pm 1$ when either $|z|=1$
or $e^{-\varepsilon / 2} < |z| < e^{\varepsilon / 2}$. This uses
Corollary \ref{cor:FS:coeffuniqueness}.

\begin{theorem}\label{thm:fundamentalsystem}
  Let $z\neq \pm 1$ and $\lambda = z+z^{-1}$ and assume either
   \begin{enumerate}[label={(\alph*)}]
      \item the first moment condition and consider
      $|z|=1$,
      \item or the exponential moment condition and consider
      $e^{-\varepsilon / 2} < |z| < e^{\varepsilon / 2}$.
   \end{enumerate}
   If $A_{n}^{-1}\in\mathscr{B}(H)$ for all $n\in\mathbb{Z}$, then
   $\{U^+(z),\, U^+(z^{-1})\}$ and $\{U^-(z),\, U^-(z^{-1})\}$ are
   fundamental systems of \eqref{eq:eigenvalue-operator}. That is, for
   every formal operator solution $V(z)$ of
   \eqref{eq:eigenvalue-operator} (or formal vector solution $V(z)$ of
   \eqref{eq:eigenvalue}) there exists a unique pair of operator (or
   vector) coefficients $P^{\pm}(z)$ and $Q^{\pm}(z)$ for which
   \begin{equation}\label{eq:lincombequal}
      V(z)=U^{\pm}(z) P^{\pm}(z) + U^{\pm}(z^{-1})Q^{\pm}(z).
   \end{equation}
   These coefficients are given explicitly by
   \begin{equation}\label{eq:betterdefpq}
      \begin{bmatrix} P^{\pm}(z) \\ Q^{\pm}(z)
      \end{bmatrix}
      = \frac{1}{\pm(z^{-1}-z)}
      \begin{bmatrix} -W(U^{\pm}(\bar z^{-1})^*,V(z))\\ W(U^{\pm}(\bar z)^*,V(z))
      \end{bmatrix}\,,
   \end{equation}
   where in the case of formal vector solutions, the Wronskians are
   understood in the sense of Definition~\ref{def:Wronskianvector}~(2).
   \end{theorem}
\begin{proof}
  By Proposition \ref{prop:wronsks},
  $\{R(z), \widehat R(z)\}:=\{U^{\pm}(z), U^{\pm}(z^{-1})\}$ satisfy the
  conditions of Proposition \ref{prop:FS:coeff}. Hence, by Corollary~\ref{cor:FS:coeffuniqueness}, 
  we only need to show that the operator
  matrix $Z_0(z)$ from \eqref{eq:Z:FS} is injective. We shall work
  with $\{R(z), \widehat R(z)\}:=\{U^{+}(z), U^{+}(z^{-1})\}$ and suppress
  the superscript $+$; the proof for 
  $\{R(z), \widehat R(z)\}:=\{U^{-}(z), U^{-}(z^{-1})\}$ is similar. 

  Suppose $[X_0,X_1]^T \in \ker Z_0(z)$. Let $X(z)$ be the formal operator
  solution of \eqref{eq:eigenvalue-operator} with initial data
  $X_0(z)$ and $X_1(z)$. Then by direct computation
  \[0=Z_0(z)\begin{bmatrix}X_0\\ X_1\end{bmatrix} =
    \frac{1}{\pm(z^{-1}-z)}\begin{bmatrix} -W_1(U(\bar
                             z^{-1})^*,X(z))\\
                              W_1(U(\bar z)^*,X(z))\end{bmatrix}.\]
  Since the Wronskians here are independent of the lower index
  $n \in \mathbb Z$ by Corollary \ref{cor:wronskianofsols}, we must
  have for all $n$
   \begin{equation}\label{eq:wronskswithx}
      W_n(U(\bar z^{-1})^*,X(z)) =  W_n(U(\bar z)^*,X(z)) = 0.
   \end{equation}
Writing out these equations gives for all $n$
\begin{align}
  \label{eq:wronskwithx2} 
  U_{n-1}(\bar z^{-1})^*A_{n-1}X_n(z) &= U_{n}(\bar z^{-1})^*A_{n-1}X_{n-1}(z),
  \\
  \label{eq:wronskwithx1} 
  U_{n-1}(\bar z)^*A_{n-1}X_n(z) &= U_{n}(\bar z)^*A_{n-1}X_{n-1}(z),
\end{align}
As a consequence of \eqref{eq:jostplusasymp}, for sufficiently large
$n>0$, $U_{n}(\bar z)^*$ is invertible and thus for sufficiently
large $n$, \eqref{eq:wronskwithx2} gives
\begin{equation*}
   A_{n-1}X_n(z) = (U_{n-1}(\bar z^{-1})^*)^{-1} U_{n}(\bar z^{-1})^*A_{n-1}X_{n-1}(z),
\end{equation*}
Plug this into \eqref{eq:wronskwithx1} and
rearrange to obtain, for sufficiently large $n$
\begin{equation}\label{zero-equation}
\widetilde U_n(z)A_{n-1}X_{n-1}(z)=0,
\end{equation}
where (note that by \eqref{eq:jostplusasymp}, $U_n(\bar z^{-1})^*$ is
also invertible for sufficiently large $n$)
\begin{equation}\label{eq:tildeun}
  \widetilde U_n(z):=U_{n-1}(\bar z)^*
  (U_{n-1}(\bar z^{-1})^*)^{-1}U_n(\bar z^{-1})^* - U_n(\bar z)^*.
\end{equation}
Observe from \eqref{eq:jostplusasymp} that as $n\to\infty,$
\begin{align*}
  \widetilde U_n(z) &=z^{n-1}(I+o(1))z^{n-1}(I+o(1)) z^{-n}(I+o(1)) -
                  z^n(I+o(1))\\
                &= z^{n-2}(1-z^2)(I+o(1)).
\end{align*}
Since $z\neq\pm 1$, we conclude that $\widetilde U_n(z)$ is invertible for sufficiently large $n$.
Since $A_{n-1}$ is also invertible, it follows from \eqref{zero-equation} that
$X_{n-1}(z)=0$ for sufficiently large $n$, and therefore
$X(z)$, being a formal solution of \eqref{eq:eigenvalue-operator}, is
identically zero\footnote{Recall that any pair $(X_{n_0}, X_{n_0+1})$
determines the formal solution $X= \{X_n\}_n$ uniquely by recursion since
$A_{n}^{-1}\in\mathscr{B}(H)$.}. Hence $[X_0,X_1]^T = 0$ and
$\ker Z_0(z) = 0$. Thus $Z_0(z)$ is injective, completing the proof.
\end{proof}

A direct consequence of Theorem \ref{thm:fundamentalsystem},
establishing Jost solutions as a fundamental system, is the fact that
they can be written as linear combinations of each other. This
observation will be particularly useful below in the construction of
transfer and scattering matrices.

\begin{corollary}
  Let $A_{n}^{-1}\in\mathscr{B}(H)$ for all $n\in\mathbb{Z}$. Assume
  either
   \begin{enumerate}[label={(\alph*)}]
   \item the first moment condition and consider
     $|z|=1$,
   \item or the exponential moment condition and consider
     $e^{-\varepsilon / 2} < |z| < e^{\varepsilon / 2}$.
   \end{enumerate}
   Then for $z\neq \pm 1$ we obtain
\begin{equation}\label{eq:alphabetakey}
U^{\mp}(z)=U^{\pm}(z)\alpha^{\pm}(z) + U^{\pm}(z^{-1})\beta^{\pm}(z)
\end{equation}
with the operator coefficients $\alpha^{\pm}(z)$ and $\beta^{\pm}(z)$
given explicitly by\footnote{Note that the Wronskians in
  \eqref{eq:defalpha} and \eqref{eq:defbeta} do not depend on $n$
  since both its arguments are formal solutions of
  \eqref{eq:eigenvalue-operator} for the same $\lambda$, so we omit the
  index $n$.}
\begin{align}\label{eq:defalpha}
  &\alpha^{\pm}(z) = 
\frac{1}{ \mp(z^{-1}-z) }W(U^{\pm}(\bar z^{-1})^*,U^{\mp}(z)),
  \\ \label{eq:defbeta} &\beta^{\pm}(z)=
  \frac{1}{\pm (z^{-1}-z) }W(U^{\pm}(\bar z)^*,U^{\mp}(z)).
\end{align}
\end{corollary}

\begin{proof}
The statement follows directly from Theorem \ref{thm:fundamentalsystem}
with $V(z) = U^{\mp}(z)$.
\end{proof}

Note that the same formulae appear in e.g. \cite[(59)]{MR4385984}, 
albeit with a different notation, namely $M^z_\pm$ for $\alpha^{\pm}(z)$, and 
$N^z_\pm$ for $\beta^{\pm}(z)$. 
They are new in the infinite-dimensional setting.

\subsection {Other combinations of fundamental systems}

We shall now discuss under which conditions the Jost solutions $\{U^{+}(z),\, U^{-}(z)\}$
form a fundamental system. We first consider the easier case, where $H$
is finite dimensional. The next theorem is basically Remark
\ref{rem:FS:finitedim}.

\begin{theorem}\label{thm:fundamentalsystem-finite}
  Let $H$ be finite-dimensional and $A_{n}^{-1}$ be in
  $\mathscr{B}(H)$ for all $n\in\mathbb{Z}$. Assume
  \begin{enumerate}[label={(\alph*)}]
     \item either the first moment condition and consider $0 < |z|\leq 1$,
     \item or the exponential moment condition and consider $0 < |z| < e^{\varepsilon / 2}$.
  \end{enumerate}
   If both $W( U^+(\bar z)^*, U^-(z))^{-1}\in\mathscr{B}(H)$ and
   $W( U^-(\bar z)^*, U^+(z))^{-1}\in\mathscr{B}(H)$, then the Jost solutions
   $\{U^{+}(z),\, U^{-}(z)\}$ form a fundamental system of formal
   solutions to \eqref{eq:eigenvalue-operator}.
   \end{theorem}

\begin{proof}
  Note first that by Theorem \ref{thm:jostsolutions} and Proposition
  \ref{Jost-infinite-moment-prop}, the Jost solutions
  $\{U^{+}(z),\, U^{-}(z)\}$ exist simultaneously for the specified
  domains of $z$. Setting
  $\{ R(z), \widehat R(z) \} := \{ U^{+}(z), U^{-}(z) \} $, invertibility of the
  Wronskians $W( U^+(\bar z)^*, U^-(z))$ or
  $W( U^-(\bar z)^*, U^+(z))$ implies the conditions of Proposition
  \ref{prop:FS:coeff}.  
  Remark~\ref{rem:FS:finitedim} shows that $Z_j(z)$ from \eqref{eq:Z:FS} is a bijection so that 
  the conditions of Corollary~\ref{cor:FS:coeffuniqueness} are satisfied. 
  This proves the statement.
\end{proof}

The case where $H$ is infinite-dimensional is harder since the limit
argument in the proof of Theorem \ref{thm:fundamentalsystem} that
shows injectivity of $Z_0(z)$ for the fundamental systems $\{U^\pm(z),\, U^\pm(z^{-1}) \}$
does not work, as we do not have asymptotics for $U_n^{\pm}(z)$ as $n$ approaches $\mp\infty$.
However we still have the
following result.

\begin{theorem}
   \label{them:fundamentalsystem2}
   Let $A_{n}^{-1}$ be in $\mathscr{B}(H)$ for all
   $n\in\mathbb{Z}$. Assume 
   \begin{enumerate}[label={(\alph*)}]
      \item either the first moment condition and consider $|z|=1$,
      \item or the exponential moment condition and consider
      $e^{-\varepsilon / 2} < |z| < e^{\varepsilon / 2}$,
      \item or the first moment condition,  $|z|\leq 1$, and additionally that
      $B_{n}, A_{n}-I$ are compact for all $n\in\mathbb{Z}$.
   \end{enumerate} 
    Let $V(z)$ be a formal vector solution to \eqref{eq:eigenvalue} or a formal
    operator solution to \eqref{eq:eigenvalue-operator} 
    with $\lambda = z+z^{-1}$ and let $P^\pm(z), Q^\pm(z)$ be as in Theorem~\ref{thm:fundamentalsystem}.
Suppose that 
\begin{align*}
W( U^+(\bar z)^*, U^-(z))^{-1}, 
 \ W( U^-(\bar z)^*, U^+(z))^{-1}\in\mathscr{B}(H)).
\end{align*}
Then   $U^+(z), U^-(z)$ is a fundamental system of formal solutions to
     \eqref{eq:eigenvalue-operator} for $\lambda=z+z^{-1}$ when
     $z\neq \pm 1$. More precisely, there exists a pair of
     coefficients $P(z)$ and $Q(z)$ for which
      \begin{equation}\label{eq:lincombequal3}
	 V(z)=U^{+}(z) P(z) + U^{-}(z)Q(z).
      \end{equation}
      These coefficients are given by
      \begin{align*}\label{eq:XXXMMMM}
	 \begin{bmatrix} P(z) \\ Q(z)
	 \end{bmatrix} 
	 & =
	 \begin{pmatrix}
	    1 & -\alpha^+(z) \beta^{+}(z)^{-1} \\
	    0 & \beta^{+}(z)^{-1}
	 \end{pmatrix}
	 \begin{bmatrix} P^+(z) \\ Q^+(z)
	 \end{bmatrix}
	& = 
	 \begin{pmatrix}
	    0 & \beta^{-}(z)^{-1} \\
	    1 & -\alpha^-(z) \beta^{-}(z)^{-1}
	 \end{pmatrix}
	 \begin{bmatrix} P^-(z) \\ Q^-(z)
	 \end{bmatrix}.
      \end{align*}
\end{theorem}
\begin{proof} Let us first prove the statements under the first two assumptions.
   Assume that $W( U^+(\bar z)^*, U^-(z))^{-1}\in\mathscr{B}(H)$.
   Then the coefficient $\beta^+(z)$ in \eqref{eq:defbeta} is invertible and we can solve for $U^{+}(z^{-1})$ in \eqref{eq:alphabetakey}.
   Therefore \eqref{eq:lincombequal} shows that
   \begin{align*}
      V(z) 
      & = U^{+}(z) P^{+}(z) + U^{+}(z^{-1})Q^{+}(z) 
      \\
      & = U^{+}(z) P^{+}(z) + \left[ U^{-}(z) - U^{+}(z)\alpha^+(z) \right] \beta^{+}(z)^{-1} Q^+(z)
      \\
      & = U^{+}(z) \left[ P^{+}(z) - \alpha^+(z) \beta^{+}(z)^{-1} Q^+(z) \right] + U^{-}(z) \beta^{+}(z)^{-1} Q^+(z).
   \end{align*}

   The proof in the case when $W( U^-(\bar z)^*, U^+(z))^{-1}\in\mathscr B(H)$ is
   analogous. One uses that $\beta^-(z)$ is invertible. \medskip
   
   Let us now prove the statements under the third assumption.
   In fact, we prove only the first statement, the second statement is argued verbatim.
By induction from the recursive formulae \eqref{eq:knrecur}, we conclude that all $K_{n,m}$
are compact.  By the formula \eqref{telescope},  the operator $T_n - I$ is compact as well since the subspace of compact operators is closed in $\mathscr{B}(H)$.  
Similarly,  $R_n - I$ is compact.
We conclude in view of \eqref{eq:jostuplus} and \eqref{TR-asymptotics} that 
$$
U_n^\pm(z) = z^{\pm n} (1 + K^\pm_n(z)), 
$$
with $K^\pm_n(z)$ being compact operators.  Consequently, using the fact that $A_n-I$ is compact, the Wronskians
$$
W(U^{\mp}(\bar z)^*, U^{\pm}(z)) = \pm(z-z^{-1})(I + K^{\pm}(z)),
$$
where $K^{\pm}(z)$ is compact. By assumption (1) and the second sentence in Proposition \ref{prop:FS:coeff}, these are both invertible, and so
$$
W(U^{\mp}(\bar z)^*, U^{\pm}(z))^{-1} = \pm\frac{1}{z-z^{-1}}(I + \tilde K^{\pm}(z)),
$$
where $\tilde K^{\pm}(z)$ is also compact. Thus the operator $Z_j(z)$ in \eqref{eq:Z:FS}, with $\{R(z), \widehat R(z)\} = \{U^+(z), U^-(z)\}$, is of the form
\begin{equation}\label{ZZZ}
   Z_j(z) = - \frac{1}{z-z^{-1}} 
   \begin{pmatrix}
      z^{-j-1} & - z^{-j} \\
      -z^{j+1} & z^j
   \end{pmatrix} + K_j(z)
\end{equation}
where $K_j(z)$ is a compact operator on $\mathscr B(H) \oplus \mathscr B(H)$.
The determinant of the first summand in \eqref{ZZZ} equals $1$ and hence it
is invertible.  Therefore $Z_j(z)$ is Fredholm with index zero. 
By the proof of Proposition~\ref{prop:FS:coeff},  $Z_j(z)$ 
admits a right inverse $Y_j(z)$ and hence is surjective.  Fredholm index zero implies that
$Z_j(z)$ is also injective and hence by Corollary \ref{cor:FS:coeffuniqueness}
the claim follows.
\end{proof}

Analogously it can be shown that $\{U^+(z), U^-(z^{-1})\}$ is a
fundamental system if $W( U^-(\bar z^{-1})^*, U^+(z))^{-1}\in\mathscr B(H)$
and that $\{ U^-(z), U^+(z^{-1}) \}$ is a fundamental system if
$W( U^+(\bar z)^*, U^-(z^{-1}))^{-1}\in\mathscr B(H)$. 

From \eqref{eq:lincombequal3} we obtain a
relation between $P^+(z), Q^+(z)$ and $P^-(z), Q^-(z)$.
\begin{corollary}
  Let $A_{n}^{-1}\in\mathscr{B}(H)$ for all $n\in\mathbb{Z}$ and
  $P^\pm(z), Q^\pm(z)$ be as in Theorem~\ref{thm:fundamentalsystem}.
  If $W( U^+(\bar z)^*, U^-(z))^{-1}$ and hence also $W( U^-(\bar z)^*, U^+(z))^{-1}$ are
  in $\mathscr{B}(H)$, then
   \begin{align*}
      \begin{bmatrix} P^+(z) \\ Q^+(z)
      \end{bmatrix}
      & =
      \begin{pmatrix}
	 \alpha^+(z) & [ 1 - \alpha^+(z)\alpha^-(z) ] \beta^{-}(z)^{-1} \\
	 \beta^{+}(z)&  -\beta^{+}(z)\alpha^-(z) \beta^{-}(z)^{-1}
      \end{pmatrix}
      \begin{bmatrix} P^-(z) \\ Q^-(z)
      \end{bmatrix},
      \\[2ex]
      \begin{bmatrix} P^-(z) \\ Q^-(z)
      \end{bmatrix}
      & =
      \begin{pmatrix}
	 \alpha^-(z) & [ 1 - \alpha^-(z)\alpha^+(z) ] \beta^{+}(z)^{-1} \\
	 \beta^{-}(z)&  -\beta^{-}(z)\alpha^+(z) \beta^{+}(z)^{-1}
      \end{pmatrix}
      \begin{bmatrix} P^+(z) \\ Q^+(z)
      \end{bmatrix}.
   \end{align*}
\end{corollary}
\begin{proof}
   This follows from \eqref{eq:lincombequal3} and
      \begin{align*}
	 \begin{pmatrix}
	    1 & -\alpha^+(z) \beta^{+}(z)^{-1} \\
	    0 & \beta^{+}(z)^{-1}
	 \end{pmatrix}^{-1}
	 &=
	 \begin{pmatrix}
	    1 & \alpha^+(z) \\
	    0 & \beta^{+}(z)
	 \end{pmatrix} \\
	 \begin{pmatrix}
	    0 & \beta^{-}(z)^{-1} \\
	    1 & -\alpha^-(z) \beta^{-}(z)^{-1}
	 \end{pmatrix}^{-1}
	 &=
	 \begin{pmatrix}
	    \alpha^-(z) & 1 \\
	    \beta^{-}(z) & 0
	 \end{pmatrix}.
	 \qedhere
      \end{align*}
\end{proof}

\section {Construction of the scattering matrix}\label{sec:constr-scatt-matr}

\subsection {Construction of the transfer and scattering matrices}

The \emph{transfer matrix} is defined as the matrix $\mathtt M(z)$,
cf.  \cite[Definition 2]{MR4192212}, whose
entries are operator-valued functions (defined on a subset of $\mathbb{C}$) such
that
\begin{equation}
   \label{eq:transfermatrixproperty}
   \begin{bmatrix}U_n^{-}(z),\ U_n^{-}(z^{-1})
   \end{bmatrix}
   = \begin{bmatrix}U_n^{+}(z),\  U_n^{+}(z^{-1})
   \end{bmatrix} \mathtt M(z)
   \qquad\text{for all }\ n\in\mathbb Z.
\end{equation}
In fact, when $A_{n}^{-1}\in\mathscr{B}(H)$ for all $n\in\mathbb{Z}$ and the
first moment condition of Definition~\ref{assumption:introduction} holds, such a
matrix does exist for $|z|=1, z\neq\pm 1$. Indeed, we may take
\begin{equation}\label{eq:defoftransfermatrix}
  \mathtt M(z) = \begin{bmatrix}\alpha^{+}(z) &
         \beta^{+}(z^{-1})\\ \beta^{+}(z) &
                                            \alpha^{+}(z^{-1})\end{bmatrix}\,,
\end{equation}
where the functions $\alpha^{+}(z)$ and $\beta^{+}(z)$ are given in
\eqref{eq:defalpha} and \eqref{eq:defbeta}, respectively.
\begin{proposition}
  Suppose that $A_{n}^{-1}\in\mathscr{B}(H)$ for all $n\in\mathbb{Z}$ and that
  the first moment condition of Definition~\ref{assumption:introduction}
  holds. For $z\in\mathbb{C}$ such that $|z|=1, z\neq\pm 1$
\label{prop:alphainvertible}
   \begin{enumerate}
      \item\label{item:Minv}
      the matrix $\mathtt M(z)$ is invertible,
      \item\label{item:alphainv}
      the operators $\alpha^{\pm}(z)$ are invertible.
   \end{enumerate}
\end{proposition}
\begin{proof} 
   (\ref{item:Minv})
   We claim that the inverse of $\mathtt M(z)$ is 
   \begin{equation}
      \mathtt M^{-}(z)
      :=\begin{bmatrix}\alpha^{-}(z) & \beta^{-}(z^{-1})\\ 
	 \beta^{-}(z) & \alpha^{-}(z^{-1})
      \end{bmatrix}.
   \end{equation}
   If we multiply \eqref{eq:transfermatrixproperty} on the right by $\mathtt M^{-}(z)$, we obtain, using  
   \eqref{eq:alphabetakey}, that 
   \begin{align*}
      \begin{bmatrix}U_n^{+}(z),\ U_n^{+}(z^{-1})
      \end{bmatrix} 
      &= \begin{bmatrix}U_n^{+}(z),\  U_n^{+}(z^{-1})
      \end{bmatrix} \mathtt M(z) \mathtt M^-(z), \\
      \begin{bmatrix}U_n^{-}(z),\ U_n^{-}(z^{-1})
      \end{bmatrix} 
      &= \begin{bmatrix}U_n^{-}(z),\  U_n^{-}(z^{-1})
      \end{bmatrix} \mathtt M^-(z) \mathtt M(z),
   \end{align*}
   where the second relation is obtained similarly. Since
   $\{ U^+(z), U^+(z^{-1})\}$ is a fundamental system for
   \eqref{eq:eigenvalue-operator} by Theorem~\ref{thm:fundamentalsystem}, we
   conclude by uniqueness of the linear combinations \eqref{eq:lincombequal}
   that $\mathtt M(z) \mathtt M^-(z) = \mathtt M^-(z) \mathtt M(z) = I$.
   \smallskip

   (\ref{item:alphainv})
   First observe that since $W(U^*,V)^*=-W(V^*,U)$ 
   \begin{equation}\label{eq:alphaadjoint}
      \alpha^{\pm}(z)^*
      = \pm \frac{1}{\bar z^{-1}-\bar z} W(U^{\mp}(z)^*,U^{\pm}(\bar z^{-1}))=
      \alpha^{\mp}(\bar z^{-1}) \equiv \, \alpha^{\mp}(z),
   \end{equation}
   where we used $\bar z=z^{-1}$ in the last step. Similarly
   \begin{equation}\label{eq:betaadjoint}
      \beta^{\pm}(z)^*=-\beta^{\mp}(\bar z).
   \end{equation}
   Writing out the top left entry of $\mathtt M^-(z) \mathtt M(z)$ gives
   \[\alpha^{-}(z)\alpha^{+}(z) + \beta^{-}(z^{-1})\beta^{+}(z)=I.\]
   Using \eqref{eq:alphaadjoint} and \eqref{eq:betaadjoint}, we obtain
   for $\bar z=z^{-1}$
   \begin{equation}\label{eq:alphabetarelation}
      \alpha^{\pm}(z)^*\alpha^{\pm}(z) = I + \beta^{\pm}(z)^*\beta^{\pm}(z).
   \end{equation}
   The right-hand side is the identity plus a positive operator and is
   therefore a bijection with bounded inverse.
   Multiplication of the whole equation by its inverse on the
   left shows that $\alpha^{\pm}(z)$ has a left inverse. By
   \eqref{eq:alphaadjoint}, its adjoint also has a left inverse, and thus
   $\alpha^{\pm}(z)$ has a right inverse as well, completing the proof.
\end{proof}
We now construct the \emph{scattering} matrix $\mathtt S(z)$. It is defined,
presuming it exists, as an operator-valued matrix for which

\begin{equation}\label{eq:scatteringmatrixproperty}
   \begin{bmatrix}U_n^{-}(z),\ U_n^{+}(z^{-1})
   \end{bmatrix}
   = \begin{bmatrix}U_n^{+}(z),\ U_n^{-}(z^{-1})
   \end{bmatrix} \mathtt S(z),
   \qquad\text{for all }\ n\in\mathbb Z.
\end{equation}
\begin{proposition}
  If one assumes that $A_{n}^{-1}\in\mathscr{B}(H)$ for all $n\in\mathbb{Z}$ and
  that the first moment condition of Definition~\ref{assumption:introduction}
  holds, then the scattering matrix exists for $|z|=1, z \neq \pm 1$ and is
  given by
\begin{equation}\label{eq:scatteringmatrixdef}
  \mathtt S(z) =
  \begin{bmatrix} (\alpha^{-}(z))^{-1} &
                                         -\beta^{+}(z^{-1})(\alpha^{+}(z^{-1}))^{-1}
    \\ -\beta^{-}(z)(\alpha^{-}(z))^{-1} & (\alpha^{+}(z^{-1}))^{-1}
\end{bmatrix}\,,
\end{equation}
where the functions $\alpha^{+}(z)$ and $\beta^{+}(z)$ are given in
\eqref{eq:defalpha} and \eqref{eq:defbeta}, respectively.
\end{proposition}
\begin{proof}
  Take \eqref{eq:alphabetakey} for $U^+(z)$, multiply on the right by
  $(\alpha^{-}(z))^{-1}$, and solve for $U^{-}(z)$ to get the first row of
  $\mathtt S(z)$. For the second row, take \eqref{eq:alphabetakey} for
  $U^-(z^{-1})$, multiply on the right by $(\alpha^+(z^{-1}))^{-1}$, and
  solve for $U^+(z^{-1})$.
\end{proof}

\begin{remark}
\label{sec:relat-cont-case} (Interpretation and relation to the continuous case)

\begin{enumerate}
\item The transfer matrix expresses the Jost solutions with known asymptotic
behavior as $n\to -\infty$ in terms of the Jost solutions with known
asymptotic behavior as $n\to +\infty$. In other words, it shows how a
solution behaving asymptotically as a ``plane wave'' at $n=\infty$ is
transformed by the transfer matrix into a solution that asymptotically
behaves as a plane wave at $n=-\infty$.

\item For the scattering matrix, we interpret the solutions $U^{+}(z)$ and
$U^{-}(z^{-1})$ as the outgoing ones, moving to the right as $n\to\infty$
and to the left as $n\to -\infty$. Conversely, the solutions $U^{-}(z)$ and
$U^{+}(z^{-1})$ are incoming solutions. If we know the coefficients of the
incoming solutions, then multiplying them by $\mathtt S(z)$ gives the
coefficients of the outgoing solutions.
\end{enumerate}
\end{remark}

\subsection {Continuity of the scattering matrix}\label{sec:cont-scat-matrix}

The main objective of this subsection is to prove continuity of the
scattering matrix $\mathtt S(z)$ at the \emph{a priori} excluded points $z = \pm 1$,
under the second moment condition and the condition that the Wronskian has closed range.
\begin{theorem}\label{thm:contscat}
  Assume that $A_{n}^{-1}\in\mathscr{B}(H)$ for all $n\in\mathbb{Z}$ and
  that the second moment condition of
  Definition~\ref{assumption:introduction} holds. 
  If the Wronskian $W(U^{+}(z_0)^*, U^{-}(z_0))$
  has closed range for $z_0= \pm 1$, then the scattering matrix
  $\mathtt S(z)$, defined initially for $|z|=1$ but $z\neq \pm 1$, has a
  continuous extension to the entire unit circle $|z|=1$.
\end{theorem}

We begin by collecting some facts about the various entries of
$\mathtt S(z)$. All of these hold even if we only assume the first moment
condition.

\begin{proposition}\label{prop:listoffacts}
  Suppose that $A_{n}^{-1}\in\mathscr{B}(H)$ for all $n\in\mathbb{Z}$ and that
  the first moment condition of Definition~\ref{assumption:introduction}
  holds. Then the following statements hold for any $z \in \mathbb C$ with
  $|z|=1$ and $z\neq \pm 1$ and some $N\in\mathbb N$ sufficiently large (independent of $z$).
  Observe that in (\ref{item:invertcpct}) and (\ref{item:continuity}) below, $z=\pm 1$ is included.

  \begin{enumerate}
     \item\label{item:observatione} 
     $U_n^{\mp}(z)=z^{\pm n}\alpha^{\pm}(z) + z^{\mp n}\beta^{\pm}(z) + o(1),\ n \to \pm \infty$.
     \medskip

     \item\label{item:invertcpct} 
     $(U_{\pm n}^{\pm}(z))^{-1}
     \in\mathscr{B}(H)$ with $\|(U_{\pm n}^{\pm}(z))^{-1}\| < 2$\
     for all $n\geq N$ and $|z|=1$.
     \medskip

     \item\label{eq:forbetaalphainverse}
     $\beta^{\pm}(z)(\alpha^{\pm}(z))^{-1} =
     (U_{\pm n}^{\pm}(z^{-1}))^{-1} 
     \Big( U_{\pm n}^{\mp}(z)(\alpha^{\pm}(z))^{-1} - U_{\pm n}^{\pm}(z) \Big)$ for all $n\geq N$.
     \medskip

     \item\label{item:unifbounds}
     $\|(\alpha^{\pm}(z))^{-1}\|$ and $\|\beta^{\pm}(z)(\alpha^{\pm}(z))^{-1}\|$ are uniformly bounded in $z$.\medskip

     \item\label{item:continuity}
     If $(\alpha^{\pm}(z))^{-1}$ has a continuous
     extension to the entire unit circle $\{|z|=1\}$, so does
     $\beta^{\pm}(z)(\alpha^{\pm}(z))^{-1}$.
  \end{enumerate}
\end{proposition}
\begin{proof} The asymptotic formula in (\ref{item:observatione}) follows
  immediately from \eqref{eq:jostplusasymp} and
  \eqref{eq:alphabetakey}. Note
   that the remainder is not necessarily uniform in $z$ as $z\to\pm 1$.
   \medskip

   Item (\ref{item:invertcpct}) holds for each fixed $z$ as an immediate
   consequence of \eqref{eq:jostplusasymp}. However, the condition on the norm
   of the inverse is an open condition, that is, if it is true for some $z$ it
   is true in a small neighborhood of $z$ as well. Since the unit circle
   $\{|z|=1\}$ is compact, we may use finitely many such neighborhoods to cover
   the circle, proving (\ref{item:invertcpct}). \medskip

   Statement \eqref{eq:forbetaalphainverse} is a straightforward consequence of
   \eqref{eq:alphabetakey}, multiplying it on the right by
   $(\alpha^{\pm}(z))^{-1}$ and rearranging. Invertibility of
   $U_{\pm n}^{\pm}(z^{-1})$ follows from (\ref{item:invertcpct}). \medskip

   To get (\ref{item:unifbounds}), start with \eqref{eq:alphabetarelation} and
   observe that for any $v\in H$
   \begin{equation*}
      \left\langle\alpha^{\pm}(z) v,\, \alpha^{\pm}(z)v\right\rangle 
      = \|v\|^2 + \left\langle\beta^{\pm}(z) v,\, \beta^{\pm}(z) v\right\rangle 
      \geq \|v\|^2.
   \end{equation*}

   From this we deduce that $\|\alpha^{\pm}(z)v\|\geq \|v\|$ for all
   $v\in H$, which shows that $\|(\alpha^{\pm}(z))^{-1}\|\leq 1$. The
   boundedness of $\|\beta^{\pm}(z)(\alpha^{\pm}(z))^{-1}\|$ then follows from
   \eqref{eq:forbetaalphainverse} and (\ref{item:invertcpct}). Indeed, if we fix
   any $n\geq N$, it is bounded uniformly by
   \[2\|U_{\pm n}^{\mp}(z)\| + 2\|U_{\pm n}^{\pm}(z)\|,\] which is bounded
   uniformly in $z$ as the Jost solutions are continuous on the closed unit
   disk. Finally (\ref{item:continuity}) is an immediate consequence of
   \eqref{eq:forbetaalphainverse} and, again, continuity of the Jost solutions
   on the closed unit disk.
\end{proof}

In light of part (\ref{item:continuity}) of Proposition
\ref{prop:listoffacts} and the form of the scattering matrix $\mathtt S(z)$,
Theorem \ref{thm:contscat} is an immediate consequence of the following lemma.

\begin{lemma}
   If one assumes that the second moment condition in Definition~\ref{assumption:introduction} holds,
   that $A_{n}^{-1}\in\mathscr{B}(H)$ for all $n\in\mathbb{Z}$
   and that the Wronskian $W(U^{+}(z_0)^*, U^{-}(z_0))$
   has closed range for $z_0= \pm 1$, 
   then $(\alpha^{\pm}(z))^{-1}$ has a continuous extension to $\{|z|=1\}$.
\end{lemma}

\begin{proof} We prove this lemma for $(\alpha^{+}(z))^{-1}$; the proof for
   $(\alpha^{-}(z))^{-1}$ is analogous.  The proof consists of three steps.\medskip
   
   \noindent \textbf{Step 1 $-$ Behaviour of $\alpha^{+}(z)$ at $\pm 1$.} \medskip
 
   \noindent Recall that $\alpha^{+}(z)$ is defined and continuous on $\{|z| = 1, z\neq \pm 1\}$ if the second moment condition holds and on $\{e^{-\epsilon/2} < |z| < e^{\epsilon/2}, z\neq \pm 1\}$ if the exponential moment condition holds.
   Let us fix $z_0=\pm 1$. Define the operators
   $\delta U_n^{\pm}(z)$, for all $z$ with $|z|=1$ and $z\neq z_0$, and all
   $n\geq 1$, by the formula
   \begin{equation}\label{U-exp}
   U_n^{\pm}(z)=U_n^{\pm}(z_0)+(z-z_0)\delta U_n^{\pm}(z).
   \end{equation}
   Then from \eqref{eq:jostuplus} we find for all $z \neq z_0$
   \begin{equation}\label{delta-U-series}
      \delta U_n^+(z) 
      = T_n\frac{z^n-z_0^n}{z-z_0} 
      + T_n \sum_{m=1}^\infty K_{n, m} \frac{z^{m+n} - z_0^{m+n}}{z - z_0} .
   \end{equation}
   First note that $|z^{m+n}-z_0^{m+n}|\leq (m+n)|z-z_0|$. 
   By Corollary \ref{lem:josthelp}, we have for some constant $C>0$
   \begin{equation*}
      \sum_{m=1}^{\infty} (m+n) \|K_{n,m}\|\leq 
      C \sum_{m=1}^{\infty} (m+n) \sum_{p=n+\lfloor\frac m2\rfloor}^{\infty}(\|I-A_p\|+\|B_p\|).
   \end{equation*}
   This infinite sum may be rearranged 
   as before and yields
   \begin{align*}
      \sum_{m=1}^{\infty} (m+n) \|K_{n,m}\|
      & \leq C \sum_{q=n}^{\infty}\frac 12(2(q-n)+1)(2(q-n)+2)(\|I-A_q\|+\|B_q\|)
      \\
      & \phantom{\leq\ } + C n \sum_{q=n}^{\infty} (2(q-n)+1) (\|I-A_q\|+\|B_q\|) 
   \end{align*}
   which is convergent by the second moment condition. Hence
   \eqref{delta-U-series} is well-defined for all $|z|=1$.
   Furthermore, $\delta U_n^{\pm}(z)$ has a
   well-defined limit as $z\to z_0$, which yields a continuous extension 
   \begin{equation}\label{U-exp-delta}
   \delta U_n^+(z_0) := \lim_{z\to z_0}\delta U_n^+(z) = nT_n z_0^{n-1} + \sum_{m= 1}^\infty (m+n) K_{n,m}z_0^{m-1}.
   \end{equation}
   The convergence here is not necessarily uniform in $n \in \mathbb Z$. An
   analogous result holds for $\delta U_n^-(z)$ as $z\to z_0$. Now, using
   \eqref{eq:defalpha}, we write, for $|z|=1, z\neq \pm 1$, and a fixed and sufficiently large $n$, suppressing that subscript,
   \begin{equation}\label{alpha-taylor}
      \alpha^{+}(z) = \frac{1}{z-\bar z} W( U^{+}(z)^*, U^{-}(z)) 
      = \frac{1}{z-\bar z} \Big( W(z_0) + (z-z_0)\delta W(z) \Big)
   \end{equation}
   where
   \begin{align*}
      W(z_0) &:=W(U^{+}(z_0)^*, U^{-}(z_0)), 
      \\
      \delta W(z) &:= W(U^{+}(z_0)^*, \delta U^-(z)) + W(\delta U^+(z)^*, U^{-}(z_0)) 
      \\ &\qquad + (z-z_0)W(\delta U^+(z)^*, \delta U^-(z)).
   \end{align*}
   Observe that $\delta W(z)$ has a continuous extension to $z_0$, with
   \begin{align}\label{deltaW-cts}
      \lim_{z\to z_0}\delta W(z) = W(U^{+}(z_0)^*, \delta U^-(z_0)) + W(\delta U^+(z_0)^*, U^{-}(z_0)).
   \end{align}

   \noindent \textbf{Step 2 $-$ Applying the Schur complement theorem to $\alpha^{+}(z)$.} \medskip
 
   \noindent Let $P$ be the orthogonal projection onto
      $\ker W(z_0)$,  and $P^\perp=I-P$ be the orthogonal projection onto the orthogonal complement
      $\ker W(z_0)^\perp$.  Let $Q$ be the orthogonal projection onto
      $\range W(z_0)^\perp$,  and $Q^\perp=I-Q$ be the orthogonal projection onto the orthogonal complement
      $\range W(z_0)$ (recall that $\range W(z_0)$ is closed by assumption).  By construction
      \begin{align}\label{QWP}
      Q W(z_0) =0, \quad W(z_0) P = 0.
      \end{align}
     We use these orthogonal projections to decompose
      \begin{align*}
      &\alpha^{+}(z) = \begin{bmatrix} A & B\\ C & D\end{bmatrix}: H = \ker W(z_0) \oplus 
      \ker W(z_0)^\perp \to \range W(z_0)^\perp \oplus \range W(z_0), \\
      &\textup{where} \quad A := Q \alpha^{+}(z) P,  \  B:= Q \alpha^{+}(z) P^\perp,  \  
      C = Q^\perp \alpha^{+}(z) P,  \  D := Q^\perp \alpha^{+}(z) P^\perp.
	 \end{align*}
	 The individual blocks follow directly from \eqref{alpha-taylor} and \eqref{QWP}
	 \begin{align*}
	    A(z) &= \frac{z-z_0}{z-\bar z} Q \delta W(z) P,\quad
	    \ \  B(z) = \frac{z-z_0}{z-\bar z}Q\delta W(z)P^\perp,\\
	    C(z) &= \frac{z-z_0}{z-\bar z}Q^{\perp}\delta W(z)P,\quad
	    D(z) = \frac{1}{z-\bar z}Q^\perp W(z_0)P^{\perp} + \frac{z-z_0}{z-\bar z}Q^\perp \delta W(z)P^{\perp}.
	 \end{align*}
	 Note that $(z-z_0)/(z-\bar z)$ is continuous and 
	 tends to $\frac{1}{2}$ as $z\to z_0$ on the unit circle.
	 Hence by \eqref{deltaW-cts}, $A, B, C$ have a continuous
	 extension to $z=z_0$, with the respective limits denoted by write
	 $A(z_0), B(z_0), C(z_0)$.\medskip

	 We first observe that $D$ is invertible on a (punctured) neighborhood of $z_0$ on the unit circle.
	 To see this,  observe that since $\range W(z_0)$ is closed,  it is a Banach space and hence
	 $$
	 Q^\perp W(z_0)P^{\perp}:  \ker W(z_0)^\perp \to \range W(z_0)
	 $$
	 is a bijective,  linear bounded operator between Banach spaces. 
	 By the closed graph theorem,  it admits a bounded inverse\footnote{
	 Our closed range assumption is essential to get a bounded inverse.}.
	 Therefore
	 \begin{equation*}
	    D(z) = \frac{1}{z-\bar z}Q^\perp W(z_0)P^{\perp} 
	    \Big[ I +  (z-z_0) \Big( Q^\perp W(z_0)P^{\perp}  \Big)^{-1}
	    Q^\perp \delta W(z)P^{\perp} \Big].
	 \end{equation*}
	 Hence a Neumann series argument shows that $D(z)$ is invertible for $z$ sufficiently close to $z_0$ with inverse
	 \begin{align*}
	 D(z)^{-1} 
	    &= (z-\bar z)
	    \Big[ I + (z-z_0) \Big( Q^\perp W(z_0)P^{\perp} \Big)^{-1}Q^\perp \delta W(z)P^{\perp} \Big]^{-1}
	   (Q^\perp W(z_0)P^{\perp})^{-1}  \\
	    &= (z - z_0)
	    \frac{z-\bar z}{z- z_0}
	    \Big[ I + (z-z_0) \Big( Q^\perp W(z_0)P^{\perp} \Big)^{-1}Q^\perp \delta W(z)P^{\perp} \Big]^{-1}
	    \\ & \phantom{\ =\ }\times  (Q^\perp W(z_0)P^{\perp})^{-1}. 
	 \end{align*}

	 Hence $D(z)^{-1}$ has a continuous extension to $z_0$ since $\delta W$ is continuous in $z_0$. Even better: it is
	 $(z-z_0)$ times a continuous function down to $z=z_0$. 
	 \medskip

	 Let us write $\Sigma(z) := A(z) - B(z)D^{-1}(z)C(z)$ for the so-called \define{Schur complement}.
	 Then we can factorize
	 %
	 \begin{align}\label{alpha-schur}
	    \begin{bmatrix} I & -BD^{-1}\\  0 & I
	    \end{bmatrix}
	    \alpha^+(z) 
	    \begin{bmatrix} I & 0 \\ -D^{-1}C & I 
	    \end{bmatrix}
	    =
	    \begin{bmatrix} \Sigma(z) & 0\\ 0 & D(z)
	    \end{bmatrix}.
	 \end{align}
	 Since $\alpha^+(z)$ and the triangular matrices on the left hand side of \eqref{alpha-schur} are boundedly invertible for $|z|=1, z\neq \pm 1$, the right hand side is boundedly invertible too and 
	 \begin{align}\label{alpha-schur-inverse}
	    \begin{bmatrix} I & 0 \\ D^{-1}C & I 
	    \end{bmatrix}
	    (\alpha^+(z))^{-1}
	    \begin{bmatrix} I & BD^{-1}\\  0 & I
	    \end{bmatrix}
	    =
	    \begin{bmatrix} (\Sigma(z))^{-1} & 0\\ 0 & (D(z))^{-1}
	    \end{bmatrix}.
	 \end{align}
	 \medskip
	        
   \noindent \textbf{Step 3 $-$ Limit of $(\alpha^{+}(z))^{-1}$ as $z
     \to \pm 1$.} \medskip
 
   \noindent Now we study $A(z)$ as $z \to z_0=\pm 1$. 
   Recall that $A(z), B(z), C(z)$ extend continuously to $z=z_0$, and that $D(z)^{-1}$ is $(z-z_0)$ times a continuous function at $z=z_0$.
   Therefore, the norms of the triangular matrices in \eqref{alpha-schur-inverse} are uniformly bounded by some $\gamma > 0$ for all $|z|=1$
   and we obtain that $\|(\Sigma(z))^{-1} \| \le \gamma^2 \|\alpha^+(z)^{-1}\| \le 2\gamma^2 $ for $|z|=1, z\neq \pm 1$.
   Now we have for $|z|=1, z\neq \pm 1$, that
   \begin{equation*}
      A(z_0) = \Sigma(z_0) 
      = \Sigma(z) \Big( I + (\Sigma(z))^{-1} (\Sigma(z_0) - \Sigma(z)) \Big).
   \end{equation*}
   For $z$ in a small enough punctured neighbourhood of $z_0$ on the unit circle, the second term is a bijection with bounded inverse by an Neumann series argument.
   Hence $A(z_0)$ is a bijection with bounded inverse.
   \medskip

   Finally, 
   as $z \to z_0$ on the unit circle, we obtain from \eqref{alpha-schur-inverse}
   \begin{equation*}
      (\alpha^+(z))^{-1}
      = 
      \begin{bmatrix} I & 0 \\ -D^{-1}C & I 
      \end{bmatrix}
      \begin{bmatrix} (\Sigma(z))^{-1} & 0\\ 0 & (D(z))^{-1}
      \end{bmatrix}
      \begin{bmatrix} I & -BD^{-1}\\  0 & I
      \end{bmatrix}
      \longrightarrow 
      \begin{bmatrix} (A(z_0))^{-1} & 0\\ 0 & 0
      \end{bmatrix}.
   \end{equation*}
\end{proof}

\begin{remark}
Under the assumption that $\range W(z_0)$ is closed, the last formula shows in particular that $\alpha^+(z)$ has a continuous extension to $z_0$ if and only if $P^\perp =0$ and $Q^\perp = 0$, that is, if and only if 
$\ker W(z_0) = \{0\}$ and 
$\range W(z_0) = H$.
\end{remark}

\section {Accumulation of the discrete spectrum}\label{sec:finit-discr-spectr}

The aim of this section is to study the discrete spectrum of $\mathcal J$,
\emph{i.\,e.} the isolated eigenvalues of finite multiplicity of $\mathcal J$. In
particular we will be interested in the question of whether the eigenvalues
accumulate at the edges of the continuous spectrum $[-2,2]$. We begin with
some fundamental observations that are in particular consequences of our
scattering matrix construction.

\subsection {Properties of the discrete spectrum}\label{subsec:finit-discr-spectr}
\begin{proposition}\label{propC}
   Let $A_{n}^{-1}\in\mathscr{B}(H)$ and let $B_{n}, A_{n}-I$ be compact operators for all $n\in\mathbb{Z}$. 
   If the first moment condition in Definition~\ref{assumption:introduction} holds,
   then $\mathcal J$ has no eigenvalues embedded in the interior of its essential spectrum,
   that is,
   $(-2,\, 2) \cap \sigma_{p}(\mathcal J) = \varnothing$.
  \end{proposition}
\begin{proof}
  Recall from Proposition \ref{prop:Jspectrum} that
  $\sigma_{ess}(\mathcal J) = [-2,2]$. 
  Suppose for the sake of a contradiction that $\lambda\in (-2,\, 2)$ is an eigenvalue of $\mathcal J$. 
  Let $\lambda = z + z^{-1}$. 
  If $\lambda\in (-2,\, 2)$, then $|z|=1$ with $z\neq \pm 1$. 
  Since $\lambda$ is an eigenvalue, there
  exists a non-zero vector $v(z)\in\mathcal H = \ell^2(\mathbb Z,H)$ which solves
  \eqref{eq:eigenvalue}. Note that $\|v_n(z)\|\to 0$ as $n\to\infty$ and
  that, by Theorem \ref{thm:jostsolutions}, $|U_n^+(z)|$ is bounded as
  $n\to\infty$. Therefore, using the Wronskian from the second part of
  Definition \ref{def:Wronskianvector} we obtain
   \begin{equation}\label{wronski-zero}
\lim_{n\to\infty} W_n(U^+(\overline z^{\pm 1})^*,v)=0.
\end{equation}
   Since $U^+(z)$ is a formal solution of \eqref{eq:eigenvalue-operator} and
   $v(z)$ is a solution of \eqref{eq:eigenvalue}, Corollary~\ref{cor:wronskianofsols} 
   shows that their Wronskian is independent of $n$
   and therefore must equal zero for all $n \in \mathbb Z$. 
   Due to the fact that $|z|=1$ and $z\neq \pm 1$, one can apply
   Theorem~\ref{thm:fundamentalsystem} (more precisely
   \eqref{eq:betterdefpq}) to show that $v=0$.
 \end{proof}
 \begin{remark}
   \label{rem:not-embedded-eigenvalues}
   It is known \cite{MR0211028} that for discrete Schrödinger operators with scalar
   coefficients
   there are no eigenvalues in $(-2,\,2)$ 
   if the potential sequence $\{q_{n}\}_{n\in\mathbb{N}}$ is
   $o(1/n)$ as $n\to\infty$.
   Note that in the scalar case, the first moment condition guarantees the
   absence of embedded eigenvalues in $(-2,\,2)$. This is also true
   for the case of operator coefficients when the Jacobi operator is
   given by \eqref{eq:orthogonal-decomposition-matrix} (see Example~\ref{ex:other-diagonal})
 \end{remark}
 \begin{proposition}
    \label{prop:edge-points-no-eigenvalues}
    Let $A_{n}^{-1}\in\mathscr{B}(H)$ and
    $B_{n}, A_{n}-I$ be compact for all $n\in\mathbb{Z}$. If
    the second moment condition in Definition \ref{assumption:introduction}
    holds, then the points ${\pm 2}$ are not eigenvalues.
 \end{proposition}
 \begin{proof}
   If $|\lambda|=2$, then
   $z= \pm 1$ and writing out the formula for the Wronskian in the equality
   $W_{n+1}(U^+(\overline z^{\pm 1})^*,v) = 0$ from \eqref{wronski-zero}, we conclude for all $n$
   \begin{equation*}
      U_{n}^+(\pm 1)^*A_nv_{n+1} = U_{n+1}^+(\pm 1)^*A_nv_{n}.
   \end{equation*}
   For sufficiently large $n$, we have invertibility of 
$U_{n}^+(\pm 1)$ by Theorem~\ref{thm:jostsolutions} and therefore
   \begin{equation*}
      v_{n+1} = A_n^{-1}U_{n}^+(\pm 1)^*{}^{-1}U_{n+1}^+(\pm 1)^*A_nv_{n}.
   \end{equation*}
   Again using \eqref{eq:jostplusasymp} for the second moment condition,  
   and also using the second moment condition to 
conclude that $A_n = I + o(n^{-2})$ as $n \to \pm \infty$, we see that, as $n\to\infty$,
   \[v_{n+1} = (I + o(n^{-1}))v_n.\]
   In order to obtain a contradiction, we want to apply this relation recursively. 
   For sufficiently large $n \in \mathbb N$ we have
\[\|v_{n+1}\|\ge \left(1-\frac{1}{2n}\right)\|v_n\|.\]
Applying this inequality recursively,  we obtain that for any sufficiently large $n_0$, 
and for any $N\in\mathbb N$,  using the double factorial from combinatorics
\begin{align*}
\|v_{n_0+N}\| &\ge \prod_{k=0}^{N-1}\left(1-\frac{1}{2(n_0+k)}\right)\|v_{n_0}\| \\
&=\frac{(2n_0-1)(2n_0+1)(2n_0+3)\dots(2n_0+2N-3)}{2n_0(2n_0+2)(2n_0+4)\dots(2n_0+2N-2)}\|v_{n_0}\|\\
&=\frac{(2n_0-2)!!}{(2n_0-3)!!}\frac{(2n_0+2N-3)!!}{(2n_0+2N-2)!!}\|v_{n_0}\|.
\end{align*}
By Stirling's formula, there is a positive constant $C=C(n_0)$ for which
\[\frac{(2n_0+2N-3)!!}{(2n_0+2N-2)!!} = \frac{(2n_0+2N-2)!}{2^{2N}((n_0+N-1)!)^2}\sim\frac{C}{\sqrt N}.\]
But since $\{v_n\}_{n\in\mathbb Z}\in\ell^2(\mathbb Z,H)$, this can only happen if $v_{n_0}=0$. Since this is true for all sufficiently large $n_0$, uniqueness implies that the entire sequence $\{v_n\}=0$, a contradiction.
\end{proof}

We now discuss the relationship between invertibility of the Wronskians and
eigenvalues under the first and the exponential moment conditions.

\begin{lemma}\label{corollaryB-Boris}
Let $A_{n}^{-1}\in\mathscr{B}(H)$ for all $n\in\mathbb{Z}$.  Assume either
\begin{enumerate}[label={(\alph*)}]
   \item $\dim H < +\infty$, the first moment condition and $0 < |z|\leq 1$.
   \item $\dim H = +\infty$, and either 
\begin{itemize}
\item the first (or higher) moment condition,  and $|z|=1,  z \neq \pm 1$,
\item or the exponential moment condition, and 
$e^{-\varepsilon / 2} < |z| \le 1,  z \neq \pm 1$,
\item or the first moment condition, $|z|\leq 1$, $z\neq\pm 1$, and that $B_{n}, A_{n}-I$ are compact for all $n\in\mathbb{Z}$.
\end{itemize}
\end{enumerate}
If $\lambda = z + z^{-1}$ is an eigenvalue of $\mathcal J$, then
$W(U^+(\bar z)^*, U^-(z))$ is not invertible.
\end{lemma}
\begin{proof}
   Consider an eigenvalue $\lambda = z+ z^{-1}$ of $\mathcal J$ with the
   non-trivial eigenfunction $w\in\ell^2(\mathbb Z, H)$. Let us assume to the
   contrary that $W(U^+(\bar z)^*,U^-(z))$ is invertible. Then for
   $\dim H < +\infty$, by Theorem \ref{thm:fundamentalsystem-finite},
   $\{U^+(z), U^-(z)\}$ forms a fundamental system. Hence
   \[
   w=U^+(z) \alpha(z) + U^{-}(z) \beta(z),
   \]
with the coefficients given, due to \eqref{eq:FS:operatorcoeff}, by
\begin{align*}
\alpha(z) &= W( U^-(\overline z)^*,  U^+(z) )^{-1} W(U^-(\overline z)^*, w(z)),\\ 
\beta(z) &= W(U^+(\bar z)^*,U^-(z))^{-1} W( U^+(\bar z)^*,w(z) ).
\end{align*}
Since $w(z)\in\ell^2(\mathbb Z, H)$ and $U_n^-(z) = z^{-n}(I + o(1))$ for
$n\to-\infty$, we conclude
\begin{align*}
W( U^-(\bar z)^*,w(z) ) = \lim_{n\to-\infty} W_n( U^-(\bar z)^*,w(z) ) = 0.
\end{align*}
Therefore $\alpha(z) =0$. Similarly we can show that $\beta(z) =0$, in
contradiction to our assumption that $w\neq 0$. For $\dim H = \infty$ the
proof repeats the same argument and uses Theorem
\ref{them:fundamentalsystem2} instead of Theorem
\ref{thm:fundamentalsystem-finite}.
\end{proof}

The converse,  at least in the region $0 < |z| < 1$, 
is easier and does not require any distinct discussion of finite
and infinite-dimensional $H$.

\begin{lemma}\label{corollaryB-Boris2} 
Suppose $A_{n}^{-1}\in\mathscr{B}(H)$ for all $n\in\mathbb{Z}$.
Assume the first moment condition
  of Definition \ref{assumption:introduction} and $0 < |z| < 1$. If
  $W(U^+(\bar z)^*, U^-(z))$ is not invertible, then $\lambda = z + z^{-1}$ is an
  eigenvalue of $\mathcal J$.
\end{lemma}

\begin{proof}
Let $0 < |z|<1$ and take a non-zero vector $v\in\ker W(U^+(\bar z)^*,U^-(z))$. 
By definition of the Wronskian, we obtain for any $n \in \mathbb Z$
   \[U_{n-1}^+(\bar z)^*A_{n-1}U_n^{-}(z)v = U_n^+(\bar
   z)^*A_{n-1}U_{n-1}^{-}(z)v.\] 
   Recall that $U_n^+(z)$ is invertible for sufficiently large $n$ by Theorem \ref{thm:jostsolutions} 
   and that $A_n^{-1}\in \mathscr B(H)$.
   So we get
  \[U_n^{-}(z)v =
    A_{n-1}^{-1}[U_{n}^+(z)U_{n-1}^+(z)^{-1}]^*A_{n-1}U_{n-1}^{-}(z)v.\]
  Again by Theorem \ref{thm:jostsolutions}, the term in brackets on the
  right-hand side is $z(I+o(1))$ as $n\to + \infty$. Since $A_n$ and
  $A_{n}^{-1}$ are bounded as $n\to\infty$, we see that
  \[U_{n}^{-}(z)v=z(I+o(1))U_{n-1}^{-}(z)v, \ n \to \infty.\] Fixing $n_0 \in \mathbb N$
  sufficiently large so that $\|z(I+o(1))\| \leq q < 1$, we find for all
  $N\in\mathbb N$
  \[ \| U_{N+n_0}^{-}(z)v \| \leq q^N  \| U_{n_0}^{-}(z)v \| ,\] and therefore
  $U_{n}^{-}(z)v$ is exponentially decreasing as $n\to\infty$. By
  \eqref{eq:jostplusasymp} we conclude
  $\{U_n^{-}(z)v\}_n \in\ell^2(\mathbb Z,H)$. Therefore it is indeed an
  eigenfunction of $\mathcal{J}$.
\end{proof}

\subsection {Accumulation of eigenvalues}\label{sec:overv-known-results}

Since $\mathcal J$ is a bounded operator on $\mathcal H$, its eigenvalues
cannot accumulate at infinity and can only accumulate at the essential
spectrum $\sigma_{ess}(\mathcal J)$. 
Proposition \ref{prop:Jspectrum} gives sufficient conditions for
\begin{align*}
      \sigma_{ess}(\mathcal J) = [-2,2].
\end{align*}
Consequently, under these conditions finiteness of the discrete spectrum follows, once we rule out
$\pm 2$ as accumulation points. We study this question in the following two
regimes:

\begin{enumerate}
    \item \textbf{the third moment condition and a closed range assumption},\smallskip
   \item \textbf{the general case}. \smallskip
\end{enumerate}
In the first regimes we obtain the desired non-accumulation.

The case $\dim H=1$ has been studied in
\cite{MR332065,MR332066,MR332023,MR0198301,guseinov76half,guseinov1976,guseinov1977}
(see Teschl \cite{teschl}). It is proved there that the corresponding first
moment condition is sufficient for having finite discrete spectrum.
Non-accumulation for the case $\dim H<+\infty$ has been obtained in
\cite{MR4192212} under the assumption that the perturbation of the discrete free
Laplacian $\mathcal J_0$ is compactly supported. In \cite{MR4385984} the
result is extended to sequences which satisfy the first moment
condition in Definition \ref{assumption:introduction}. We note that the finite-dimensional case
was studied before (see for instance \cite{ABC}), but the first complete
proof of the finiteness of the discrete spectrum seems to be in
\cite{MR4385984}.

We could not find any instance where non-accumulation has been addressed in
the general case with possibly infinite-dimensional $H$.  We study this
situation in Theorem \ref{matrix-case-non-acc}  under additional assumptions.
Without these we could not rule
out non-accumulation and only prove some estimates on the rate of
accumulation in Theorem \ref{resonance-case-non-acc}.

\subsection {Nonaccumulation under the third moment condition with $\dim H = \infty$} 
\label{sec:nonacc-when-h} \ \medskip

Note that Ballesteros et al.
\cite{MR4192212} prove non-accumulation in case $\dim H < \infty$, 
assuming the first moment condition.  The price we pay for 
studying $\dim H = \infty$ is the stronger assumption that the third
moment condition in Definition \ref{assumption:introduction} is satisfied.
We also assume that $W:=W(U^+(1)^*, U^-(1))$ has closed range.  However, since we
assume that $B_{n}, A_{n}-I$ are compact for all $n \in \mathbb{Z}$, then
$W$ is also compact. Thus our closed range assumption is equivalent to requiring
that $W$ is a finite-rank operator.  Since finite-rank operators are dense
in the space of compact operators, this is a reasonable condition.

\begin{theorem}\label{matrix-case-non-acc} 
Assume the third moment condition in 
Definition \ref{assumption:introduction}.  
Assume also  that
$B_{n}, A_{n}-I$ are compact for all $n \in \mathbb{Z}$ and moreover that
$$
W(U^+(\pm 1)^*, U^-(\pm 1)) \ \textup{has closed range.}
$$
Then the eigenvalues $z+z^{-1},\ z \in \mathbb R$,
do not accumulate at $z+z^{-1} = \pm 2$.  In particular, if the Wronskians at both $\pm 1$ have closed range, the discrete spectrum of $\mathcal J$ is finite. 
\end{theorem}

\begin{proof} 
Let $\{\lambda_j\}_j\in\mathbb R\backslash[-2,2]$ be a sequence of 
eigenvalues of $\mathcal J$ and let $\{x_j\} \subset (-1,1)$ such that $\lambda_j = x_j + x_j^{-1}$.

Then $W_j:=W(U^+(x_j)^*, U^-(x_j))$ is not invertible by the last condition in Lemma~\ref{corollaryB-Boris}.
Now let us assume to the contrary that the eigenvalues
accumulate e.g. at $\lambda = 2$.
Let $\lambda_j$ be a sequence of eigenvalues in $(2,\infty)$ with corresponding $x_j\in (-1, 1)$ such that $\lambda_j\to 2$. Then $x_j \to 1$. 
We abbreviate 
\begin{align*}
   W_j &:= W(U^+(x_j)^*, U^-(x_j)),\\
   W   &:= W(U^+(1)^*, U^-(1)).
\end{align*}
The proof now follows in several steps. \medskip

\noindent \textbf{Step 1 $-$ Asymptotic properties of $W_j$ as $x_j \to 1$.}
\medskip

\noindent It follows by induction from the recursive formulae \eqref{eq:knrecur} that all $K_{n,m}$
are compact.  By the formula \eqref{telescope},  the operator $T_n - I$ is compact as well since the subspace of compact operators is closed in $\mathscr{B}(H)$.  
Similarly,  $R_n - I$ is compact.
We conclude in view of \eqref{eq:jostuplus} and \eqref{TR-asymptotics} that 
$$
U_n^\pm(z) = z^{\pm n} (1 + K^\pm_n(z)), 
$$
with $K^\pm_n(z)$ being compact operators with 
$\| K^\pm_n(z) \| = o(1)$ as $n \to \pm \infty$ uniformly in $z$.
Consequently,  we obtain (recall $A_{n}-I$ are compact for all $n \in \mathbb{Z}$)
\begin{equation}\label{WWj}
   \begin{split}
      W_j &= W(U^+(x_j)^*, U^-(x_j)) = \left(\frac{1}{x_j}-x_j\right) I + K_j,  \\
      W   &= W(U^+(1)^*, U^-(1)) = K,
   \end{split}
\end{equation}
with $K_j$ and $K$ being compact for all $j$.  By continuity of 
$U^\pm(z)$ in $z$, we find that $W_j$ converges to $W$ and in particular, 
$K_j$ converges to $K$.  In fact,  we can say even more: 
By \eqref{U-exp}, we can write $U^\pm(x_j) = U^\pm(1) + (x_j - 1)\delta U^\pm(x_j)$ where the limit
$\lim_{x_j\to 1}\delta U^\pm(x_j)$ exists.
Thus we can write for any $n\in\mathbb Z$,
\begin{align*}
   W_j & = W( U^+(x_j)^*, U^-(x_j) )
   \\&=  W\left( U^+(1)^* + (x_j-1)\delta U^+(x_j)^*,\ U^-(1) + (x_j - 1)\delta U^-(x_j) \right)
   \\ &=   W( U^+(1)^*, U^-(1) )
   \\ &\phantom{=\ }+ (x_j-1) \Big[ W_n( U^+(1)^* , \delta U^-(x_j) )
   + W_n( \delta U^+(x_j)^*, U^-(1) ) \Big]
   \\ &\phantom{=\ }+ (x_j-1)^2 W_n( \delta U^+(x_j)^*, \delta U^-(x_j) )
\end{align*}
and we obtain with \eqref{WWj} that
\begin{equation}\label{KKj}
\begin{split}
   K_j &= W_j - \left( \frac{1}{x_j} - x_j \right) I 
   = W_j - \frac{ (x_j-1)(x_j+1) }{x_j}  I 
   \\
   & = K + (x_j - 1)\left[ W_n\left( U^+(1)^* , \delta U^-(x_j) \right) + 
   W_n\left( \delta U^+(x_j)^*, U^-(1) \right) - \frac{ x_j+1 }{x_j}\, I\right]
   \\ & \phantom{=\ }+ (x_j-1)^2 W_n\left( \delta U^+(x_j)^*, \delta U^-(x_j) \right)
 =: K + (x_j-1) K'_j.
\end{split}
\end{equation}
Now for $x_j\neq 1$, we have that $K'_j = (x_j-1)^{-1} (K-K_j)$ is compact because 
$K$ and $K_j$ are compact.
Since $(K'_j)_{j\in\mathbb N}$ converges in in the operator norm, its limit for 
$j\to\infty$ (that is, $x_j\to 1$)
\begin{align}
   \label{KK-limit:}
   K'_\infty := \lim_{j\to\infty} K'_j 
   = W_n( U^+(1)^* , \delta U^-(1) ) + W( \delta U^+(1)^*, U^-(1) ) - 2 I
\end{align}
is compact too.
We can also check explicitly that $K'_\infty$ is indeed compact.
We use that 
$U^\pm(1) = 1 + \text{compact}$ and, by \eqref{U-exp-delta},
$\delta U_n^+(x_j) = n x_j^{n-1} + \text{compact}$ and 
$\delta U_n^-(x_j) = -n x_j^{-n-1} + \text{compact}$.
Hence
\begin{align*}
   K'_\infty 
   &= W_n( ( I, (n)_{n\in\mathbb Z} ) + W_n( (-n)_{n\in\mathbb Z}, I ) - 2 I
   + \text{compact}
   \\
   &= n A_{n-1} - (n-1)A_{n-1} + (-(n-1))A_{n-1} - (-n)A_{n-1} - 2I + \text{compact}
   \\
   &= 2A_{n-1} - 2I + \text{compact}
   \\
   &= 2 (A_{n-1} - I) + \text{compact}
   \\
   &= \text{compact}.
\end{align*}

\noindent \textbf{Step 2 $-$ Spectral projection of $K^*K$.}
\medskip

\noindent The compact operator $K$ has closed range if and only if
it has finite rank, i.e.  we have
\begin{equation}\label{finite-range}
\dim \range W < \infty.
\end{equation}
By the spectral theorem for compact self-adjoint operators,  $\sigma(K^*K) \backslash \{0\}$ 
is discrete and consists of finite-dimensional eigenvalues that may only accumulate at $0$.  However, 
by \eqref{finite-range},  $\dim \range K^*K \leq \dim \range K < \infty$ and hence there are only finitely many 
non-zero eigenvalues of $K^*K$.  In particular,  $0 \in \sigma(K^*K)$ is an isolated point and 
there exists $\varepsilon > 0$ such that 
\begin{align*}
\overline{B_\varepsilon (0)} \cap \sigma(K^*K) = \{0\}.
\end{align*}
We can now define the spectral projection
\begin{align*}
P = \frac{1}{2\pi i} \oint\limits_{\partial B_\varepsilon (0)} (K^* K- \mu)^{-1} d\mu,
\end{align*}
which is in fact the orthogonal projection onto the infinite-dimensional 
subspace $\ker W = \ker K^*K$.  Let us now abbreviate
$$
L_j:= K^*_jK_j +  \left(\frac{1}{x_j}-x_j\right) (K^*_j +K_j), \quad 
R_j := L_j - K^* K,
$$
such that $W^*_jW_j = \left(\frac{1}{x_j}-x_j\right)^2 I + L_j$.
Then for any $\mu \notin \sigma(K^*K)$
\begin{equation}\label{resolvent-formula}
L_j - \mu = \Bigl(I + R_j (K^* K- \mu)^{-1}\Bigr)(K^* K- \mu).
\end{equation}
Since $L_j \to K^* K$,  clearly $R_j \to 0$. 
Moreover, $(K^*K-\mu)^{-1}$ is uniformly bounded on $\partial B_\epsilon(0)$.
Thus we can choose
$j_0 \in \mathbb{N}$ sufficiently large such that $\| (K^* K- \mu)^{-1} \| \|R_j\| < 1$ 
for $j \geq j_0$ and $\mu\in B_\epsilon(0)$.
Hence by applying the Neumann series argument to \eqref{resolvent-formula} we find for all $j \geq j_0$
\begin{align*}
   \partial B_\varepsilon (0) \cap \sigma(L_j) = \{0\}.
\end{align*}
This allows us to define the spectral projection 
\begin{align*}
P_j = \frac{1}{2\pi i} \oint\limits_{\partial B_\varepsilon (0)} (L_j- \mu)^{-1} d\mu.
\end{align*}
We can now estimate the norm of the difference of the kernel projections
\begin{align*}
\frac{P_j- P}{1-x_j} = \, &\frac{1}{2\pi i } \oint\limits_{\partial B_\varepsilon (0)} 
\frac{(L_j- \mu)^{-1} - (K^* K- \mu)^{-1} }{1-x_j} d\mu  \\
 = \, &\frac{1}{2\pi i} \oint\limits_{\partial B_\varepsilon (0)} 
 (L_j- \mu)^{-1} \circ \frac{(K^* K-L_j)}{1-x_j} \circ (K^* K- \mu)^{-1} d\mu
\\ \xrightarrow{j \to \infty} \, &\frac{1}{2\pi i} \oint\limits_{\partial B_\varepsilon (0)} 
 (K^* K- \mu)^{-1} \circ \delta L \circ (K^* K- \mu)^{-1} d\mu =: Q,
\end{align*}
where we used \eqref{KKj} and \eqref{KK-limit} in the last step and abbreviated 
$$
\delta L = K^* K'+K'^*K + 2(K+K^*).
$$
In particular,  we have 
$$
\| P-P_j \| \xrightarrow{j \to \infty} 0.
$$

\noindent \textbf{Step 3 $-$ Constructing a convergent subsequence of normalized $v_j \in \ker W_j$.}
\medskip

\noindent Consider $v_j \in \ker W_j$ with $\| v_j \| = 1$.  
Because of \eqref{WWj} and \eqref{KKj} we have
\begin{align*}
   &\left(\frac{1}{x_j}-x_j\right) v_j + \Bigl(K + (x_j-1) K'_j\Bigr) v_j = 0\\[1ex]
   \Longleftrightarrow\quad
   &\left(\frac{1}{x_j}-x_j\right) v_j + K v_j  = (1-x_j) K'_j v_j.
\end{align*}
Note that $v_j$ lies in the $(x_j - 1/x_j)^2$ eigenspace of $L_j$ and 
the eigenvalue $(x_j - 1/x_j)^2 \in B_\varepsilon (0)$ for $j \in \mathbb{N}$
sufficiently small.  Thus for such $j$ we have $v_j \equiv P_j v_j$. 
Writing $v_j \equiv P_j v_j = P v_j + (P_j-P)v_j$, we obtain 
\begin{align*}
&\left(\frac{1}{x_j}-x_j\right) v_j + K (P_j-P) v_j  = (1-x_j) K'_j v_j\\
\Leftrightarrow  \
& \ v_j = \frac{x_j}{1+x_j} \left( K'_j v_j - K \circ \left( \frac{P_j-P}{1-x_j} \right) v_j \right)
\\ & \quad = \frac{x_j}{1+x_j} \left( K'v_j + (K'_j-K)v_j- K \circ \left( \frac{P_j-P}{1-x_j} \right) v_j \right).
\end{align*}
The first and the last term on the right hand side of the equality above are
compact operators applied to 
bounded sequences,  and the intermediate term converges in norm to zero.  
This shows that $v_j$ admits a convergent subsequence 
$$
v_{j_k} \xrightarrow{k \to \infty} v = 
\frac{1}{2} \left( K' - K \circ Q \right) v.
$$
We shall simply write $v_j$ with $v_j \to v$ for the subsequence.  \medskip

\noindent \textbf{Step 4 $-$ Arriving at a contradiction.}
\medskip

\noindent We may now adapt an argument by Guseinov \cite[Lemma 2.4]{guseinov1977}.
The proof of Lemma \ref{corollaryB-Boris2}
asserts that $U^-(x_j)v_j \in \ell^2(\mathbb Z, H)$ is an eigenfunction of $\mathcal J$ for the eigenvalue $(x_j+x_j^{-1})$.
By \eqref{eq:jostplusasymp} and arguing exactly as in the proof of Lemma \ref{corollaryB-Boris2}, we have in particular
\begin{equation}\label{U}
\begin{split}
U^-_n(x_j)v_j &= x_j^{-n} (I+ o(n^{-2}))v_j, \quad n \to -\infty, \\
U^-_n(x_j)v_j &= x_j (I + o(n^{-2}))U^-_{n-1}(x_j)v_j, \quad n \to \infty,
\end{split}
\end{equation}
By \cite[Theorem 2.3]{welstead} the infinite product
$\, {\vphantom{\prod}}' \prod_{n=n_0}^\infty (I + o(n^{-2})) 
$
exists for any $n_0\in\mathbb Z$ and hence
\begin{equation}\label{UU}
\begin{split}
U^-_{n_0+N}(x_j)v_j = x^{N}_j (I + o(1))U^-_{n_0}(x_j)v_j, \quad n_0 \to \infty,
\end{split}
\end{equation}
where $o(1)$ is of arbitrarily small norm as $n_0 \to \infty$.
Moreover we have 
\begin{equation}\label{UUU}
   U^-_{n_0}(x_j)v_j \to U^-_{n_0}(1)v, 
   \qquad
   \| U^-_{n_0}(1)v \| \neq 0.
\end{equation}
The fact that $\| U^-_{n_0}(1)v \| \neq 0$ follows since otherwise by 
\eqref{UU} we have $U^-_{n}(1)v = U^-_{n_0}(1)v = 0$ for $n \geq n_0$ and thus by uniqueness of 
formal solutions, $U^-_{n}(1)v = 0$ for all $n \in \mathbb Z$. However, by \eqref{eq:jostplusasymp} 
this in turn would imply
$0 = U^-_{n}(1)v = (I+o(1))v \to v$ as $n \to -\infty$ and hence $v = 0$. This is a contradiction to 
$\| v \| = 1$. Using \eqref{UUU}, we may choose now $j_0 \in \mathbb N$ sufficiently large such that for any $i,j \geq j_0$
\begin{equation}\label{UUUU}
\begin{split}
&\langle v_i, v_j \rangle \geq \frac{1}{2} \| v \|^2 = \frac{1}{2}, \\
&\langle U^-_{n_0}(x_i)v_i, U^-_{n_0}(x_j)v_j \rangle \geq \frac{1}{2} \| U^-_{n_0}(1)v \|^2 > 0.
\end{split}
\end{equation}
Using \eqref{U}, we compute for any $i,j \geq j_0$
and some uniform constant $C>0$
\begin{multline*}
   \quad \sum_{n=-\infty}^{-n_0} \langle U^-_{n}(x_i)v_i, U^-_{n}(x_j)v_j \rangle 
   \\
   \begin{aligned}
      &= \sum_{n=-\infty}^{-n_0} (x_ix_j)^{-n}\langle v_i, v_j \rangle 
      + \sum_{n=-\infty}^{-n_0} (x_ix_j)^{-n}\langle o(n^{-1})v_i, o(n^{-1})v_j \rangle 
      \\ 
      &\geq  \sum_{n=-\infty}^{-n_0} (x_ix_j)^{-n}\langle v_i, v_j \rangle 
      - \sum_{n=-\infty}^{-n_0} \frac{C}{n^2}(x_ix_j)^{-n}\|v_i\| \|v_j \| 
      \\
      &\geq \sum_{n=n_0}^{\infty} (x_ix_j)^{n} \left( \frac{1}{2}-\frac{C}{n^2}\right) 
      \geq \frac{(x_ix_j)^{n_0}}{4( 1 - x_ix_j)} \geq \delta_0 > 0,
   \end{aligned}
   \qquad
\end{multline*}
for some positive $\delta_0>0$, independent of $i,j$.
We used $\|v_i\| = \|v_j\|=1$, as well as \eqref{UUUU} in the penultimate estimate.
Moreover,  we assumed $n_0$
to be possibly even larger to ensure $C/n_0^2 < 1/4$.
Similarly, possibly increasing $n_0$ even further, we compute for any $i,j \geq j_0$ using \eqref{UU} and \eqref{UUUU}, 
\begin{align*}
   \sum_{n=n_0}^{\infty} \langle U^-_{n}(x_i)v_i, U^-_{n}(x_j)v_j \rangle 
   & \geq \frac{1}{2} \sum_{n=n_0}^{\infty} (x_ix_j)^{n-n_0} \langle U^-_{n_0}(x_i)v_i, U^-_{n_0}(x_j)v_j  \rangle 
   \\ 
   & \geq \frac{1}{4( 1 - x_ix_j)}  \| U^-_{n_0}(1)v \|^2 \geq \delta_0 > 0.
\end{align*}
On the other hand, noting the eigenfunctions to different eigenvalues of self-adjoint operators are orthogonal, we compute
\begin{align*}
   0 &= \sum_{n=-\infty}^{\infty} \langle U^-_{n}(x_i)v_i, U^-_{n}(x_j)v_j \rangle
   = \sum_{|n| > n_0} \langle U^-_{n}(x_i)v_i, U^-_{n}(x_j)v_j \rangle \\ 
   &\phantom{=\ }+
   \sum_{|n| \leq n_0} \langle U^-_{n}(x_i)v_i, U^-_{n}(x_j)v_j -U^-_{n}(x_i)v_i \rangle + \sum_{|n| \leq n_0} \| U^-_{n}(x_i)v_i \|^2.
\end{align*}
Noting that $U^-_{n}(x_j)v_j -U^-_{n}(x_i)v_i \to 0$ as $i,j \to \infty$ uniformly in $|n| \leq n_0$, we conclude in view of 
\eqref{UUU}
$$
2\delta_0 \leq \sum_{|n| > n_0} \langle U^-_{n}(x_i)v_i, U^-_{n}(x_j)v_j \rangle \leq 0.
$$
This is a contradiction and hence a sequence $x_j \to 1$ may not exist. 
The statement follows.
\end{proof}
%

%

\subsection {Accumulation without the moment conditions}
\label{sec:nonacc-reson-oper}

In this subsection,  rather than proving non-accumulation of eigenvalues at $\pm 1$, we
obtain estimates on the product of eigenvalues. Let us assume here that 
$\mathcal J = \mathcal J_0 + V$ is a trace class perturbation of the
free Jacobi operator $\mathcal J_0$.
As in Proposition~\ref{prop:Jspectrum} we find that
\begin{equation}
   V = \mathcal S^* (\mathcal A- \mathcal I) + \mathcal B + (\mathcal A-\mathcal I) \mathcal S.
\end{equation}
Here,  $\mathcal I$ denotes the block diagonal operator
$\diag[ \{I\}_{n\in\mathbb Z}]$ acting on $\mathcal H = \ell^2(\mathbb Z, H)$.
We shall write $ \| \cdot \|_1$ both for the trace norm on $H$ and $\mathcal H$.
Since $\mathcal S$ is unitary,  we conclude by the observation in Lemma \ref{lem:diagcompact}
 \begin{equation}\label{V-trace-est}
   \| V \|_1 \leq 2 \sum_{n\in\mathbb Z} \|A_n - I\|_1 + \sum_{n\in\mathbb Z} \|B_n\|_1.
   \end{equation}
Hence, to ensure that the perturbation $V$ is trace class, we
assume that the operators $A_n-I$ and $B_n$ on $H$
are trace class for all $n \in \mathbb Z$ and that
\begin{align}\label{no-moment-trace}
2 \sum_{n\in\mathbb Z} \|A_n - I\|_1 + \sum_{n\in\mathbb Z} \|B_n\|_1 < \infty.
\end{align}
The resolvent kernel of $\mathcal J_0$ can be computed explicitly, as in Teschl
\cite[\S 1.3]{teschl}, and is given for $|z| < 1$ and any $n,m \in \mathbb Z$ by
$$
(\mathcal J_0 - (z+z^{-1}))^{-1} (n,m) = \frac{z}{z^2-1} \, z^{|n-m|}.
$$
Using the trace-class property \eqref{no-moment-trace}, we can 
estimate uniformly in $|z| < 1$
\begin{align}\label{JV}
\| (\mathcal J_0 - (z+z^{-1}))^{-1} \circ V \|_1 
\leq  \left| \frac{z}{z^2-1} \right| \| V \|_1.
\end{align}
The Birman-Schwinger principle states that $-1$ is an eigenvalue of
$(\mathcal J_0 - (z+z^{-1}))^{-1} \circ V$ if and only if $z+z^{-1}$ is an
eigenvalue of $\mathcal J$ (counted with multiplicities). Using Simon \cite[Theorem
3.3 and 3.9]{Simon-det} we conclude that the perturbation determinant (cf.
\cite{MR246142})
\begin{align}\label{determinant}
f(z):= \det\nolimits_1 (I + (\mathcal J_0 - (z+z^{-1}))^{-1} \circ V)
\end{align}
is an analytic function on $|z| < 1$, with zeros precisely at those values of
$z$ which correspond to eigenvalues $z+z^{-1}$ of $\mathcal J$. Thus the
analysis of eigenvalues of $\mathcal J$ is transformed into the study of zeros of
analytic functions, which may well accumulate at the boundary. Due to the
estimate \eqref{JV} and \cite[Theorem 3.2]{Simon-det}, we find for any
$R < 1$
\begin{align*}
&f(z) = 1+ O(z^{-1}), z \to 0, \\
&|f(|z|=R)| \leq \exp \left( \left| \frac{R}{1-R^2} \right| \| V \|_1 \right).
\end{align*}
Note that since $f$ is analytic in $|z| < 1$, it has only finitely many zeros in $|z| \leq R$
for any fixed $R< 1$. A direct consequence of \cite[Proposition 2.1]{hulko} 
is an estimate on the number of such eigenvalues.

\begin{theorem}\label{resonance-case-non-acc}
Consider any $R<1$ and the enumeration of the eigenvalues $\{ z_1 + z_1^{-1}, \cdots,  z_{J(R)} + z_{J(R)}^{-1}\}$ of $\mathcal J$,
with $|z_j| < R$,  counted with their multiplicities.  Then 
\begin{equation} 
\begin{split}
\prod_{j =1}^{J(R)} R |z_j|^{-1} &\leq \exp \left( \left| \frac{R}{1-R^2} \right| 
\left(2 \sum_{n\in\mathbb Z} \|A_n - I\|_1 + \sum_{n\in\mathbb Z} \|B_n\|_1\right) \right), \\
J(R) &\leq \frac{1}{\ln R} \left| \frac{R}{1-R^2} \right| 
\left(2 \sum_{n\in\mathbb Z} \|A_n - I\|_1 + \sum_{n\in\mathbb Z} \|B_n\|_1\right) .
\end{split}
\end{equation}
\end{theorem}

\begin{proof}
The first estimate follows directly from the observation in \cite[Proposition 2.1]{hulko}.
The second estimate is a direct consequence of the first, since
\begin{equation*}
\left| \frac{R}{1-R^2} \right| \| V \|_1 \geq \ln \prod_{j =1}^{J(R)} R |z_j|^{-1} 
= \sum_{j =1}^{J(R)} \Bigl( \ln R -\ln |z_j| \Bigr) \geq J(R) \ln R.
\qedhere
\end{equation*}
The statement now follows from \eqref{V-trace-est}.
\end{proof}

%
%

\appendix
\section {Relation between vector and operator solutions}\label{vector-operator-section}

In this section we elaborate on the relation between formal 
solutions to \eqref{eq:eigenvalue} and \eqref{eq:eigenvalue-operator}.
Note first that for a formal solution $U = \{U_n\}_{n\in\mathbb Z}$ of \eqref{eq:eigenvalue-operator},
and any fixed vector $v_0\in H$, the sequence $u := \{U_n v_0\}_{n\in\mathbb Z}$ is a formal solution of \eqref{eq:eigenvalue}.
Moreover, if $U \in \ell^2(\mathbb Z, \mathscr{B}(H))$, then $u \in \ell^2(\mathbb Z, H) = \mathcal H$, 
i.e. the solution is not only a formal solution but even an eigenfunction of $\mathcal J$. 
This immediately implies the following result.

\begin{proposition}[Vector solution constructed from an operator solution]
   Assume that $\lambda$ is an eigenvalue of $\mathcal J$ viewed as a 
   bounded operator on $\ell^2(\mathbb Z, \mathscr{B}(H))$,
   with an eigenvector $U \in \ell^2(\mathbb Z, \mathscr{B}(H))$.
   Then $\lambda$ is also an eigenvalue of $\mathcal J$ viewed as a 
   bounded operator on $\ell^2(\mathbb Z, H)$ and its multiplicity is at least 
   $$\sup\limits_{n\in\mathbb Z} \dim \range (U_n).$$
\end{proposition}

Passing back from vector to operator solutions is more intricate. 

\begin{proposition}[Operator solution constructed from a vector solution.]
   \label{prop:vectoroperator}
   Let $u = \{u_n\}_{n\in\mathbb Z}$ be a formal vector solution of \eqref{eq:eigenvalue}.
   For an auxiliary $v \in H\setminus \{0\}$ let $\widetilde H$ be a closed subspace of $H$ such that 
   $H = \Span\{v\} \oplus \widetilde H$.
   For $n\in\mathbb Z$ we define the linear operators $U = \{U_n\}_{n \in \mathbb Z} \subset \mathscr{B}(H)$ by 
   \begin{equation}
      \label{eq:Udef}
      U_n: H\to H,\quad
      U_n  v = u_n
      \quad\text{ and }\quad
      U_n  w = 0
      \quad\text{for } w\in\widetilde H.
   \end{equation}

   \begin{enumerate}[label={(\arabic*)}]
      \item 
      \label{item:vectoroperator:i}
      Then $U$ is a formal solution of \eqref{eq:eigenvalue-operator}.
      If $u \in \ell^2(\mathbb Z, H)$, then
      $U \in \ell^2(\mathbb Z, \mathscr{B}(H))$. 

    \item\label{item:vectoroperator:ii} Consider a tuple
      $u^{(1)}, \cdots, u^{(k)} \in \ell^2(\mathbb Z, H)$ of solutions to \eqref{eq:eigenvalue}.
     Then any formal solution $U = \{U_n\}_{n \in \mathbb Z} \subset \mathscr{B}(H)$ of
      \eqref{eq:eigenvalue-operator} with
      \begin{equation*}
	 \range \begin{pmatrix} U_0 \\ U_1 
	 \end{pmatrix}
	 = \Span\left\{
	 \begin{pmatrix} u^{(1)}_0 \\ u^{(1)}_1
	 \end{pmatrix},\, \dots,\,
	 \begin{pmatrix} u^{(k)}_0 \\ u^{(k)}_1
	 \end{pmatrix}
	 \right\}
      \end{equation*}
      belongs to $\ell^2(\mathbb Z, \mathscr{B}(H))$.
\end{enumerate}
\end{proposition}

\begin{proof}
   \ref{item:vectoroperator:i}
   By definition of $U$ we obtain
   \begin{align*}
      &\bigl(A_{n-1}U_{n-1} + B_n U_n + A_n U_{n+1}\bigr)v =
      A_{n-1}u_{n-1} + B_n u_n + A_n u_{n+1} = \lambda u_n =  \lambda U_n v,
   \end{align*}
where we use that $u$ is a formal solution to  \eqref{eq:eigenvalue} in the second equality.
By construction we obtain for any $w \in \widetilde H$
   \begin{align*}
      \bigl(A_{n-1}U_{n-1} + B_n U_n + A_n U_{n+1}\bigr)w = 0 = \lambda U_n w.
   \end{align*}
Hence $U$ is indeed a formal solution of the operator equation \eqref{eq:eigenvalue-operator}.
   Note that
   $$
\|U_n\| = \|P\| \|u_n\|\,\|v\|^{-1},
$$ 
where $P$ is the projection of $H$ onto $\Span\{v\}$ along $\widetilde H$.
   Hence $u\in \ell^2(\mathbb Z,H)$ if and only if $U\in \ell^2(\mathbb Z, \mathscr{B}(H))$.

   \ref{item:vectoroperator:ii}
   We set 
   \begin{equation}
      \label{eq:def:widetildeU}
      \widetilde U := 
      \begin{pmatrix} U_0 \\ U_1
      \end{pmatrix}: H \to H\oplus H
      \qquad\text{and}\qquad
      \widetilde u^{(j)} := 
      \begin{pmatrix} u^{(j)}_0 \\ u^{(j)}_1
      \end{pmatrix}.
   \end{equation}
   Without restriction we may assume that the $\widetilde u^{(1)}, \dots, \widetilde u^{(k)}$ are linearly independent.
   We can write $H = \widehat H \oplus\ker\widetilde U$.
   Then $\widetilde U|_{\widehat H}:\widehat H \to \Span\{ \widetilde u^{(1)}, \dots, \widetilde u^{(k)} \}$ is an isomorphism.
   Choose a basis $v_1,\dots, v_k$ of $\widehat H$ and define 
   $$
   \widetilde H^{(j)} = \Span\{v_1,\dots, v_{j-1}, v_{j+1}, \dots, v_k\} \oplus \ker\widetilde U.
   $$
   Then $ U^{(j)}$ as defined in \eqref{eq:Udef} belongs to $\ell^2(\mathbb Z, \mathscr{B}(H))$ and
   the statement follows since $U$ is a linear combination of $\{U^{(j)}\}_{j}$.
\end{proof}

In summary, we showed that $U = \{U_n\}_{n\in\mathbb Z}$ is a formal operator solution of \eqref{eq:eigenvalue-operator}, then if and only if $v\in H$, $\{U_nv\}_{n\in\mathbb Z}$ is a formal vector solution of \eqref{eq:eigenvalue}.
If $U \in \ell^2(\mathbb Z, \mathscr{B}(H))$, then $u \in\mathcal H = \ell^2(\mathbb Z, H)$. 
Conversely, $u = \{u_n\}_{n\in\mathbb Z}$ is a formal vector
solution of \eqref{eq:eigenvalue} if and only if there exists a formal
operator solution $U$ of \eqref{eq:eigenvalue-operator}, such that
$u_n = U_n v$ for some vector $v \in H$ and all $n \in \mathbb Z$. 
If $u \in \mathcal H$, then $U$ can
be chosen in $\ell^2(\mathbb Z, \mathscr{B}(H))$. 
In particular, $\lambda$ is an eigenvalue of $\mathcal J$ on $\mathcal H$ if and only if it is an
eigenvalue of $\mathcal J$ on $\ell^2(\mathbb Z,\mathscr{B}(H))$. 
\smallskip

We have the following relation between the respective multiplicities.

\begin{corollary}
   \label{cor:mutliplicities}
   A complex number $\lambda$ is an eigenvalue of $\mathcal J$ on $\mathcal H$ if and only if it is an
eigenvalue of $\mathcal J$ on $\ell^2(\mathbb Z,\mathscr{B}(H))$.
   Let $k\in\mathbb N\cup \{\infty\}$ be the multiplicity of $\lambda$ as an eigenvalue of $\mathcal J$ on $\mathcal H$ and 
   let $\widetilde k \in\mathbb N\cup \{\infty\}$ be the multiplicity of 
$\lambda$ as an eigenvalue of $\mathcal J$ on $\ell^2(\mathbb Z, \mathscr{B}(H))$. 
   Then the multiplicities $k$ and $\widetilde k$ are related as follows.
   \begin{enumerate}

      \item If $\dim H = d < \infty$, then $\widetilde k = dk$.
      \item If $\dim H = \infty$, then $\widetilde k = \infty$.

   \end{enumerate}
\end{corollary}
\begin{proof}
Let $\lambda$ be an eigenvalue of $\mathcal J$ and let
   $V = \ker (\mathcal J - \lambda) \subseteq \mathcal H$.
   Note that, with the notation from \eqref{eq:def:widetildeU},
   every operator solution of $(\mathcal J - \lambda)U = 0$ is uniquely determined by $\widetilde U$
   and every $v\in V$ is uniquely determined by $\widetilde v$.
   If we set $\widetilde V = \{ \widetilde v : v\in V \}\subseteq H\oplus H$, 
   then clearly $\dim \widetilde V = k$ and 
   Proposition~\ref{prop:vectoroperator} \ref{item:vectoroperator:ii} 
   shows that the set of all operator solutions of $(\mathcal J - \lambda)U = 0$ is isomorphic to the space 
   \begin{align*}
      W
      = \{ \widetilde U: H \to H \oplus H \text{ linear with } 
      \range \widetilde U \subseteq \widetilde V\}
   \end{align*}
   If $\dim H = d < \infty$, then $\dim W = \dim H \cdot \dim \widetilde V = dk$.
   This proves the first statement.
   If $\dim H = \infty$, then clearly $\dim W = \infty$.
   This proves the second statement.
\end{proof}

\begin{corollary}
   If $\dim H =\infty$ and the conditions of Theorem 
   \ref{matrix-case-non-acc}  or Lemma \ref{infinite-moment-case-non-acc} hold, then the eigenvalues of 
   $\mathcal J$ viewed as operator on $\ell^2(\mathbb Z, \mathscr B(H))$ do not accumulate,  but they have infinite multiplicity; therefore $\mathcal J$ as operator on $\ell^2(\mathbb Z, \mathscr B(H))$ has only essential spectrum and its discrete spectrum is empty.
\end{corollary}

\section{Different types of spectrum and various conventions}
\label{spectrum-convention-section}

Let us recall and introduce notation for the various types of spectrum.  
This is partly necessary,  due to different conventions in the literature.
Recall that for a closed linear operator on a Hilbert space $\mathcal H$ its \define{resolvent set} $\rho(T)$ consists of all points $\lambda\in \mathbb C$ such that $T-\lambda$ is a bijection from its domain $\mathcal D(T)$ to $\mathcal H$.
The \define{spectrum} of $T$ is $\sigma(T) := \mathbb C \setminus \rho(T)$.

\subsection{Point,  continuous and residual spectrum.}
There is a disjoint decomposition of the spectrum of $T$ into its 
\define{point, continuous} and
\define{residual spectrum}
\begin{equation*}
\sigma (T) = \sigma_p (T) \ \dot{\cup} \ \sigma_{c}( T ) \ \dot{\cup} \ \sigma_{res}( T )
\end{equation*}
where
\begin{align}
   \sigma_{p} (T) &:= \{ \lambda \in \mathbb C \mid (T-\lambda) 
   \ \textup{is not injective}\},  \\
   \label{eq:sigmacont}
   \sigma_{c} (T) &:= \{ \lambda \in \mathbb C \mid (T-\lambda) 
   \ \textup{is injective,} \range (T-\lambda) \neq \overline{\range(T-\lambda)} = \mathcal H\}, \\
   \sigma_{res} (T) &:= \{ \lambda \in \mathbb C \mid (T-\lambda) 
   \ \textup{is injective and }  \overline{\range(T-\lambda)} \neq \mathcal H\}.
\end{align}
The residual spectrum is empty for self-adjoint operators. Note that
\eqref{eq:sigmacont} differs from the definition of continuous
spectrum given in \cite[Sec.\,3.7]{birman-solomjak}.
\medskip

\subsection{Discrete and essential spectrum.}
For any selfadjoint operator $T$ on the Hilbert space $\mathcal H$ we have
a disjoint decomposition of the spectrum into the 
\define{discrete spectrum} and the \define{essential spectrum}
\begin{equation*}
\sigma (T) = \sigma_d (T) \ \dot{\cup} \ \sigma_{ess}( T )
\end{equation*}
where
\begin{align}
   \sigma_{ess} (T) &:= \{ \lambda \in \mathbb C \mid (T-\lambda) 
   \ \textup{is not Fredholm}\}, \\
   \sigma_d (T) &:= \{ \lambda \in \mathbb C \mid 0 < \dim \ker (T-\lambda) < \infty,  \ 
   \lambda \ \textup{not acc.  point  of } \sigma (T) \}.
\end{align}
\smallskip

\subsection{Pure point, absolutely continuous and singular spectrum.}
For selfadjoint operators $T$
there is also a third decomposition of the spectrum, in general no longer disjoint, into its 
\define{pure point, absolutely continuous} and \define{singular continuous spectrum}
\begin{equation}
   \label{eq:def:abscont}
   \sigma (T) = \sigma_{pp} (T) \ \cup \ \sigma_{ac}( T ) \ \cup \ \sigma_{sc}( T ),
   \qquad  \sigma_{pp} (T) := \overline{ \sigma_{p} (T)}.
\end{equation}
Moreover one defines the \define{singular spectrum} and the \define{continuous spectrum} by
\begin{align}
   \sigma_{sg} (T) &:= \sigma_{pp} (T) \ \cup \ \sigma_{sc}(T),
   \\
   \label{eq:def:sigmacontscatterin}
   \sigma_{cont} (T) &:= \sigma_{ac} (T) \ \cup \ \sigma_{sc}(T)
\end{align}
The definitions in \eqref{eq:def:abscont} are more involved, using spectral measures,
and we refer the reader to \cite[\S\, X.1]{kato} for details.
The set $\sigma_{cont}(T)$ from 
\eqref{eq:def:sigmacontscatterin}
should not be confused with $\sigma_c(T)$ from 
\eqref{eq:sigmacont}.

There are the following relations between the various types of spectrum:
\begin{equation*}
\sigma_{d} (T) \subseteq \sigma_{p} (T) \subseteq \sigma_{pp} (T),  \qquad
\sigma_{ac} (T) \cup \sigma_{sc} (T) \subseteq \sigma_{ess} (T).
\end{equation*}

\subsection{Differences of conventions}
Now let us address some confusing points of conventions in the literature. 
Some of the results in the field of difference equations,  such as \cite{ABC} and
\cite{mutlu}, refer for their notions of different spectra to Glazman \cite{glazman}.
However,  what \cite{glazman} calls discrete spectrum  is actually the point spectrum
$\sigma_{p} (T)$.  What \cite{glazman} calls continuous spectrum is actually the essential spectrum
$\sigma_{ess} (T)$.  

This is important,  since Weyl's theorem asserts invariance of the essential 
spectrum under compact perturbations,  but gives no statement on the continuous spectrum.
In particular,  the perturbation results in 
\cite[Theorem 3.1]{ABC} and \cite[Theorem 8]{mutlu} (for correctness of the latter result, 
the author should in fact have additionally assumed compactness of $A_n$ and $B_n$) are in fact correct.

\bibliographystyle{amsalpha}

\begin{thebibliography}{AKvdM01}

\bibitem[AG93]{MR1255973}
N.~I. Akhiezer and I.~M. Glazman, \emph{Theory of linear operators in {H}ilbert
  space}, Dover Publications Inc., New York, 1993, Translated from the Russian
  and with a preface by Merlynd Nestell, Reprint of the 1961 and 1963
  translations, Two volumes bound as one. \MR{1255973 (94i:47001)}

\bibitem[AK01]{MR1861473}
T.~Aktosun and M.~Klaus, \emph{Small-energy asymptotics for the
  {S}chr\"{o}dinger equation on the line}, vol.~17, 2001, Special issue to
  celebrate Pierre Sabatier's 65th birthday (Montpellier, 2000), pp.~619--632.
  \MR{1861473}

\bibitem[Akh65]{akhiezer}
N.~I. Akhiezer, \emph{The classical moment problem and some related questions
  in analysis}, Translated by N. Kemmer, Hafner Publishing Co., New York, 1965.
  \MR{0184042 (32 \#1518)}

\bibitem[AKvdM01]{MR1855088}
T.~Aktosun, M.~Klaus, and C.~van~der Mee, \emph{Small-energy asymptotics of the
  scattering matrix for the matrix {S}chr\"{o}dinger equation on the line}, J.
  Math. Phys. \textbf{42} (2001), no.~10, 4627--4652. \MR{1855088}

\bibitem[AKW11]{MR2894582}
T.~Aktosun, M.~Klaus, and R.~Weder, \emph{Small-energy analysis for the
  self-adjoint matrix {S}chr\"{o}dinger operator on the half line}, J. Math.
  Phys. \textbf{52} (2011), no.~10, 102101, 24. \MR{2894582}

\bibitem[AKW14]{MR3221247}
\bysame, \emph{Small-energy analysis for the selfadjoint matrix
  {S}chr\"{o}dinger operator on the half line. {II}}, J. Math. Phys.
  \textbf{55} (2014), no.~3, 032103, 25. \MR{3221247}

\bibitem[BAC16]{ABC}
E.~Bairamov, Y.~Aygar, and S.~Cebesoy, \emph{Spectral analysis of a selfadjoint
  matrix-valued discrete operator on the whole axis}, J. Nonlinear Sci. Appl.
  \textbf{9} (2016), no.~6, 4257--4262. \MR{3530129}

\bibitem[BFGSB22]{MR4385984}
M.~Ballesteros, G.~Franco, G.~Garro, and H.~Schulz-Baldes, \emph{Band edge
  limit of the scattering matrix for quasi-one-dimensional discrete
  {S}chr\"{o}dinger operators}, Complex Anal. Oper. Theory \textbf{16} (2022),
  no.~2, Paper No. 23, 31. \MR{4385984}

\bibitem[BFSB21]{MR4192212}
M.~Ballesteros, G.~Franco, and H.~Schulz-Baldes, \emph{Analyticity properties
  of the scattering matrix for matrix {S}chr\"{o}dinger operators on the
  discrete line}, J. Math. Anal. Appl. \textbf{497} (2021), no.~1, Paper No.
  124856, 27. \MR{4192212}

\bibitem[BS87]{birman-solomjak}
M.~Sh. Birman and M.~Z. Solomjak, \emph{Spectral theory of selfadjoint
  operators in {H}ilbert space}, Mathematics and its Applications (Soviet
  Series), D. Reidel Publishing Co., Dordrecht, 1987, Translated from the 1980
  Russian original by S. Khrushch{\"e}v and V. Peller. \MR{1192782 (93g:47001)}

\bibitem[Cas73]{MR332066}
K.~M. Case, \emph{On discrete inverse scattering problems. {II}}, J.
  Mathematical Phys. \textbf{14} (1973), 916--920. \MR{332066}

\bibitem[Cas74]{MR332023}
\bysame, \emph{The discrete inverse scattering problem in one dimension}, J.
  Mathematical Phys. \textbf{15} (1974), 143--146. \MR{332023}

\bibitem[CC73]{MR332067}
K.~M. Case and S.~C. Chiu, \emph{The discrete version of the {M}archenko
  equations in the inverse scattering problem}, J. Mathematical Phys.
  \textbf{14} (1973), 1643--1647. \MR{332067}

\bibitem[CFKS87]{MR883643}
H.~L. Cycon, R.~G. Froese, W.~Kirsch, and B.~Simon, \emph{Schr\"{o}dinger
  operators with application to quantum mechanics and global geometry}, study
  ed., Texts and Monographs in Physics, Springer-Verlag, Berlin, 1987.
  \MR{883643}

\bibitem[CK73]{MR332065}
K.~M. Case and M.~Kac, \emph{A discrete version of the inverse scattering
  problem}, J. Mathematical Phys. \textbf{14} (1973), 594--603. \MR{332065}

\bibitem[CS89]{MR985100}
K.~Chadan and P.~C. Sabatier, \emph{Inverse problems in quantum scattering
  theory}, second ed., Texts and Monographs in Physics, Springer-Verlag, New
  York, 1989, With a foreword by R. G. Newton. \MR{985100}


\bibitem[\`E66]{MR0198301}
M.~S. \`Eskina, \emph{The direct and the inverse scattering problem for a
  partial difference equation}, Soviet Math. Dokl. \textbf{7} (1966), 193--197.
  \MR{0198301}



\bibitem[G66]{glazman}
I.~M.  Glazman, \emph{Direct methods of qualitative spectral analysis of singular
              differential operators}, Translated from the Russian by the IPST staff,  Israel Program for Scientific Translations,  Jerusalem, (1966),  Daniel Davey \& Co., Inc.,  New York\MR{0190800}



\bibitem[GK69]{MR246142}
I.~C. Gohberg and M.~G. Kre\u{\i}n, \emph{Introduction to the theory of linear
  nonselfadjoint operators}, Translations of Mathematical Monographs, vol. Vol.
  18, American Mathematical Society, Providence, RI, 1969, Translated from the
  Russian by A. Feinstein. \MR{246142}

\bibitem[Gus76a]{guseinov76half}
G.~\v{S}. Guse\u{\i}nov, \emph{Determination of an infinite {J}acobi matrix
  from scattering data}, Dokl. Akad. Nauk SSSR \textbf{227} (1976), no.~6,
  1289--1292. \MR{0405160}

\bibitem[Gus76b]{guseinov1976}
\bysame, \emph{The inverse problem of scattering theory for a second order
  difference equation on the whole real line}, Dokl. Akad. Nauk SSSR
  \textbf{231} (1976), no.~5, 1045--1048. \MR{0454437}

\bibitem[Gus77]{guseinov1977}
\bysame, \emph{The scattering problem for an infinite {J}acobi matrix}, Izv.
  Akad. Nauk Armyan. SSR Ser. Mat. \textbf{12} (1977), no.~5, 365--379, 416.
  \MR{559594}

\bibitem[Hul17]{hulko}
Artem Hulko, \emph{On the number of eigenvalues of the discrete one-dimensional
  {S}chr\"{o}dinger operator with a complex potential}, Bull. Math. Sci.
  \textbf{7} (2017), no.~2, 219--227. \MR{3671736}

\bibitem[Kat57]{katorosenblum-kato}
Tosio Kato, \emph{Perturbation of continuous spectra by trace class operators},
  Proc. Japan Acad. \textbf{33} (1957), 260--264. \MR{92133}

\bibitem[Kat95]{kato}
T.~Kato, \emph{Perturbation theory for linear operators}, Classics in
  Mathematics, Springer-Verlag, Berlin, 1995, Reprint of the 1980 edition.
  \MR{1335452}

\bibitem[KF57]{KolmogorovFomin}
A.~N. Kolmogorov and S.~V. Fomin, \emph{Elements of the theory of functions and
  functional analysis. {V}ol. 1. {M}etric and normed spaces}, Graylock Press,
  Rochester, N.Y., 1957, Translated from the first Russian edition by Leo F.
  Boron. \MR{0085462}

\bibitem[Kla88]{MR954906}
M.~Klaus, \emph{Low-energy behaviour of the scattering matrix for the
  {S}chr\"{o}dinger equation on the line}, Inverse Problems \textbf{4} (1988),
  no.~2, 505--512. \MR{954906}

\bibitem[Lja67]{MR0211028}
V.~\`E. Ljance, \emph{A non-selfadjoint difference operator}, Dokl. Akad. Nauk
  SSSR \textbf{173} (1967), 1260--1263. \MR{211028}

\bibitem[MBO92]{MR1191552}
F.~G. Maksudov, \`E.~M. Ba\u{\i}ramov, and R.~U. Orudzheva, \emph{An inverse
  scattering problem for an infinite {J}acobi matrix with operator elements},
  Dokl. Akad. Nauk \textbf{323} (1992), no.~3, 415--419. \MR{1191552}

\bibitem[Mut20]{mutlu}
G.~Mutlu, \emph{Spectral properties of the second order difference equation
  with selfadjoint operator coefficients}, Commun. Fac. Sci. Univ. Ank. Ser.
  A1. Math. Stat. \textbf{69} (2020), no.~1, 88--96. \MR{4034546}

\bibitem[NY92a]{MR1190777}
S.~N. Naboko and S.~I. Yakovlev, \emph{The discrete {S}chr\"{o}dinger operator.
  {A} point spectrum lying in the continuous spectrum}, Algebra i Analiz
  \textbf{4} (1992), no.~3, 183--195. \MR{1190777}

\bibitem[NY92b]{MR1173094}
\bysame, \emph{The point spectrum of a discrete {S}chr\"{o}dinger operator},
  Funktsional. Anal. i Prilozhen. \textbf{26} (1992), no.~2, 85--88.
  \MR{1173094}

\bibitem[Ros57]{katorosenblum-rosenblum}
Marvin Rosenblum, \emph{Perturbation of the continuous spectrum and unitary
  equivalence}, Pacific J. Math. \textbf{7} (1957), 997--1010. \MR{90028}

\bibitem[RS79]{MR529429}
M.~Reed and B.~Simon, \emph{Methods of modern mathematical physics. {III}},
  Academic Press [Harcourt Brace Jovanovich Publishers], New York, 1979,
  Scattering theory. \MR{529429 (80m:81085)}

\bibitem[Sim77]{Simon-det}
Barry Simon, \emph{Notes on infinite determinants of {H}ilbert space
  operators}, Advances in Math. \textbf{24} (1977), no.~3, 244--273.
  \MR{482328}

\bibitem[Sim98]{MR1627806}
B.~Simon, \emph{The classical moment problem as a self-adjoint finite
  difference operator}, Adv. Math. \textbf{137} (1998), no.~1, 82--203.
  \MR{1627806 (2001e:47020)}

\bibitem[Tar61]{MR0133609}
V.~G. Tarnopol'skii, \emph{The dispersion problem for a difference equation},
  Soviet Math. Dokl. \textbf{2} (1961), 135--138. \MR{0133609}

\bibitem[Tar76]{MR0422923}
\bysame, \emph{The scattering problem for a difference equation with operator
  coefficients}, Ukrain. Mat. \v{Z}. \textbf{28} (1976), no.~3, 342--351, 429.
  \MR{0422923}

\bibitem[Tes00]{teschl}
G.~Teschl, \emph{Jacobi operators and completely integrable nonlinear
  lattices}, Mathematical Surveys and Monographs, vol.~72, American
  Mathematical Society, Providence, RI, 2000. \MR{1711536}

\bibitem[TL86]{TaylorLay}
A.~E. Taylor and D.~C. Lay, \emph{Introduction to functional analysis}, second
  ed., Robert E. Krieger Publishing Co., Inc., Melbourne, FL, 1986. \MR{862116}

\bibitem[Wel85]{welstead}
S.~T. Welstead, \emph{Infinite products in a {B}anach algebra}, J. Math. Anal.
  Appl. \textbf{105} (1985), no.~2, 523--532. \MR{778485}

\end{thebibliography}

\def\cprime{$'$} \def\lfhook#1{\setbox0=\hbox{#1}{\ooalign{\hidewidth
  \lower1.5ex\hbox{'}\hidewidth\crcr\unhbox0}}}
\providecommand{\bysame}{\leavevmode\hbox to3em{\hrulefill}\thinspace}
\providecommand{\MR}{\relax\ifhmode\unskip\space\fi MR }
\providecommand{\MRhref}[2]{%
  \href{http://www.ams.org/mathscinet-getitem?mr=#1}{#2}
}
\providecommand{\href}[2]{#2}

\end{document}